\documentclass[12pt]{article}

\usepackage[a4paper]{geometry}
\usepackage{graphicx}

\usepackage{amsmath,amssymb,amsfonts,amsthm}
\usepackage{amscd}
\usepackage{float}
\usepackage{subcaption}
\usepackage[all]{xy}
\usepackage{tikz}
\theoremstyle{definition}
\newtheorem{thm}{Theorem}
\newtheorem*{thm*}{Theorem}
\newtheorem{lemma}[thm]{Lemma}
\newtheorem{claim}[thm]{Claim}
\newtheorem{prop}[thm]{Proposition}

\newtheorem{crl}[thm]{Corollary}
\theoremstyle{definition}

\theoremstyle{definition}
\newtheorem{dfn}[thm]{Definition}
\newtheorem{question}[thm]{Question}

\title{Torelli group action on a colored extended Hatcher complex}

\author{Hiromasa Fuchizawa}
\date{}
\begin{document}

\maketitle

\begin{abstract}

We prove that the Torelli group acts freely on a colored extended Hatcher complex without using the fact that the Torelli group is torsion-free.

\end{abstract}

\section{Introduction}

Let $\Sigma_{g,r}$ 
be an orientable surface of genus $g$ with $r$ boundary components. Assume that Euler characteristic $\chi (\Sigma_{g,r})=2-2g-r<0$.

Recall that a pair of pants is a surface of genus $0$ with $3$ boundary components, and a pants decomposition of a surface
is a decomposition of the surface to some pairs of pants.

In [NS], the authors defined a quilt decomposition of a surface such as definition 1 of this paper, 
which is constructed by adding ``seams" to a pants decomposition.
For example, a quilt decomposition of a surface is as shown in Figure \ref{fig:quiltdecomposition}.

In [F2], we defined a colored quilt decomposition of a surface such as definition 2 of this paper, which is constructed
by adding ``color" to a quilt decomposition. 
For example, a colored quilt decomposition of a surface is as shown in Figure \ref{fig:letterYonsurface} or 
Figure \ref{fig:coloredquiltonsurface}. 
In [F2], we defined a two-dementional cell complex called a colored extended Hatcher complex (definition 3 of this paper),
such that each 0-cell is a colored quilt decomposition of a surface, each 1-cell is a ``move" of four types between
two colored quilt decompositions, and each 2-cell is one of eight types of certain 2-cells.
In [F2], we proved that it is simply-connected (theorem 4 of this paper). 

In [F2], we proved that the Torelli group (or more precisely, a torsion-free subgroup of a mapping class group of a surface)
 acts on a colored extended Hatcher complex freely (theorem 5 of this paper). 
The Torelli group is the subgroup of the mapping class group that acts on the first homology of the surface trivially.
The Torelli group is a normal subgroup of the mapping class group.

In this paper, we prove that the Torelli group acts on a colored extended Hatcher complex freely
without using the fact that the Torelli group is torsion-free.

Theorem 5 suggests that there is a (infinite) presentation of the Torelli group.
This is because, for the 2-cell complex obtained by quotienting the colored extended Hatcher complex by the action of the Torelli
group, the fundamental group is the Torelli group, and all the 2-cells have been described.
However, it is difficult to calculate the presentation by this theorem, because the fundamental domain of the action is too large.

We conclude this paper with two questions that arise naturally from our result.

\section{Torelli group action on a colored extended Hatcher complex}

Let $\Sigma_{g,r}$ 
be an orientable surface of genus $g$ with $r$ boundary components. Assume that Euler characteristic $\chi (\Sigma_{g,r})=2-2g-r<0$.
(We will sometimes omit the subscripts $g$ and $r$ if they are unimportant.)

\begin{figure}[p]
    \centering
    \includegraphics[width=0.4\linewidth]{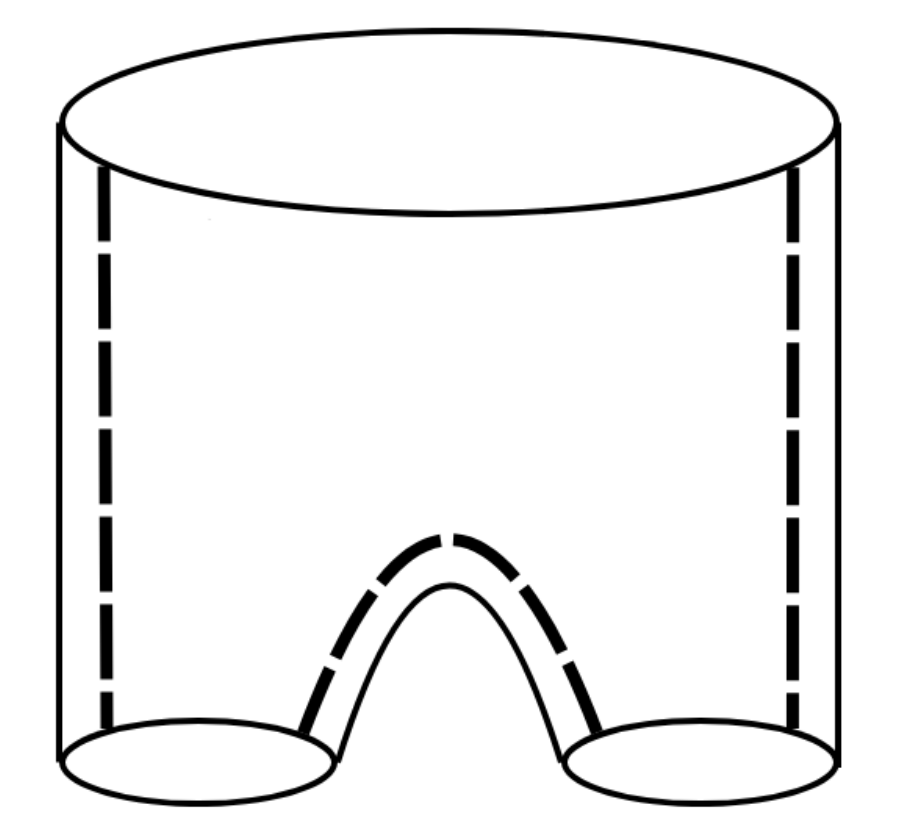}
    \caption{Quilt decomposition of a pair of pants}
    \label{fig:quiltonpants}
\end{figure}

\begin{figure}[p]
    \centering
    \includegraphics[width=0.9\linewidth]{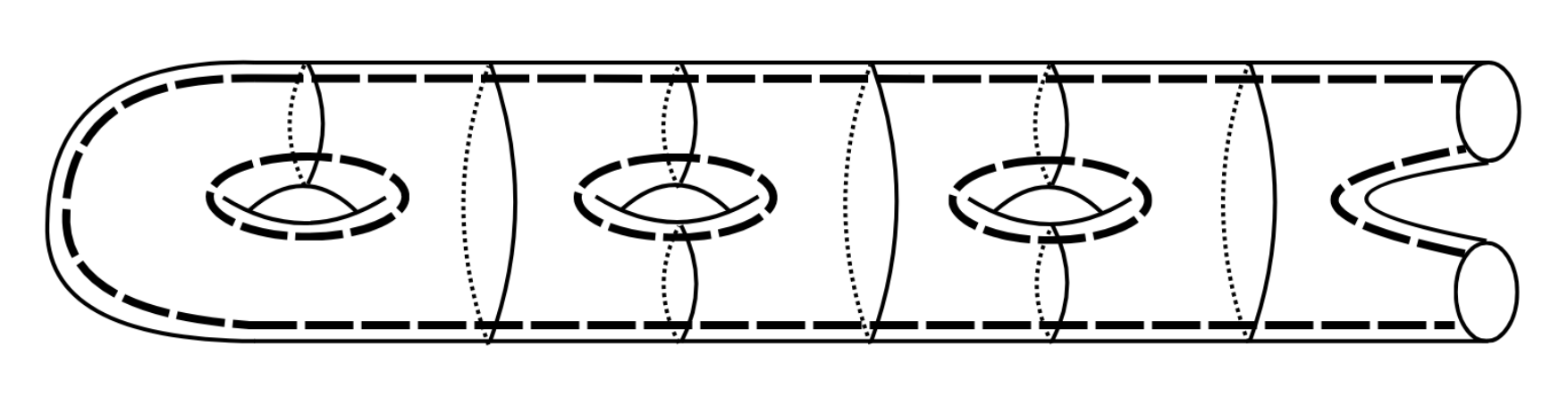}
    \caption{Quilt decomposition of a surface}
    \label{fig:quiltdecomposition}
\end{figure}

\begin{dfn} (a quilt decomposition of a surface) [NS]

Let a {\it quilt decomposition of a pair of pants} be an isotopy class of a collection of 3 seams on the pair of pants (as Figure \ref{fig:quiltonpants}). We call a hexagonal domain bounded by SCCs (simple closed curves) and seams a {\it patch}.

Let a {\it quilt decomposition of a surface} $\Sigma_{g,r}$ be a pants decomposition with each pair of pants having a quilt decomposition of a pair of pants, such that every seam is connected on each SCC of a pants decomposition (Figure \ref{fig:quiltdecomposition}).

\end{dfn}

\begin{figure}[p]
    \centering
    \includegraphics[width=0.4\linewidth]{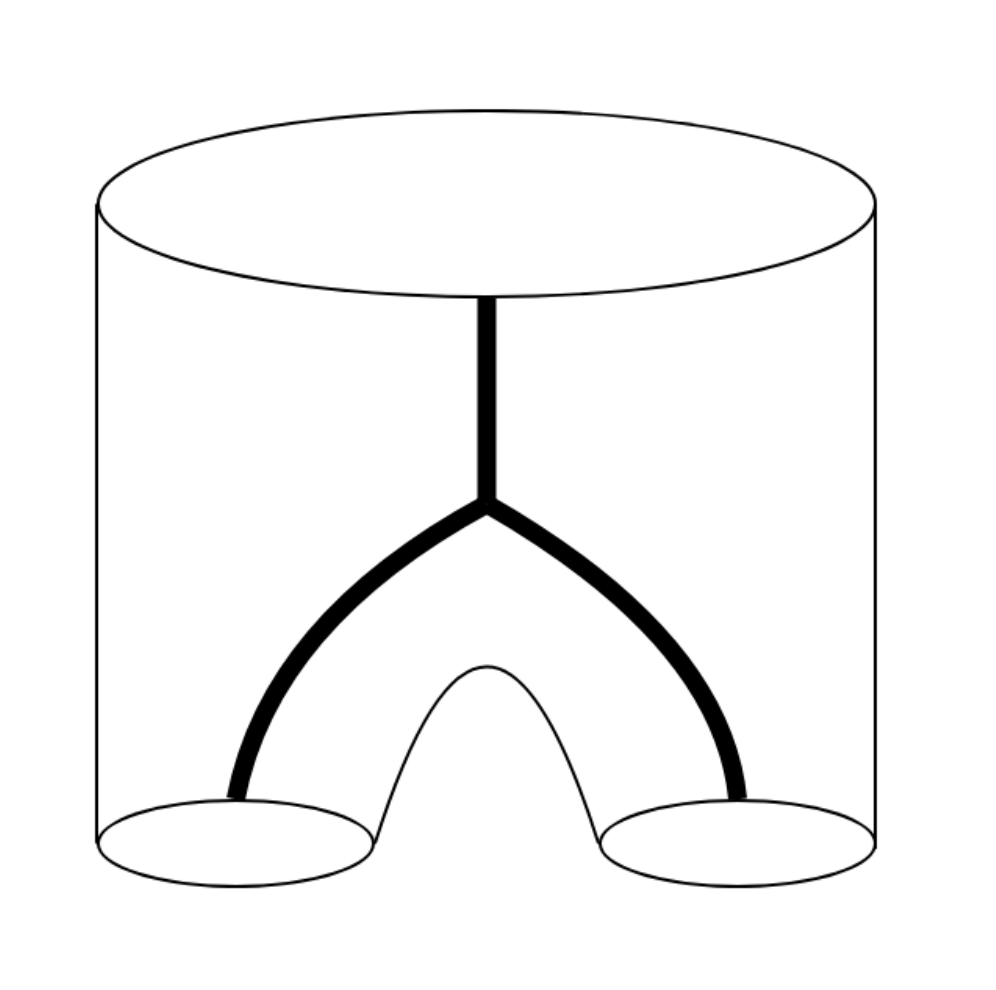}
    \caption{colored quilt decomposition of a pair of pants (I)}
    \label{fig:letterYonpants}
\end{figure}

\begin{figure}[p]
    \centering
    \includegraphics[width=0.9\linewidth]{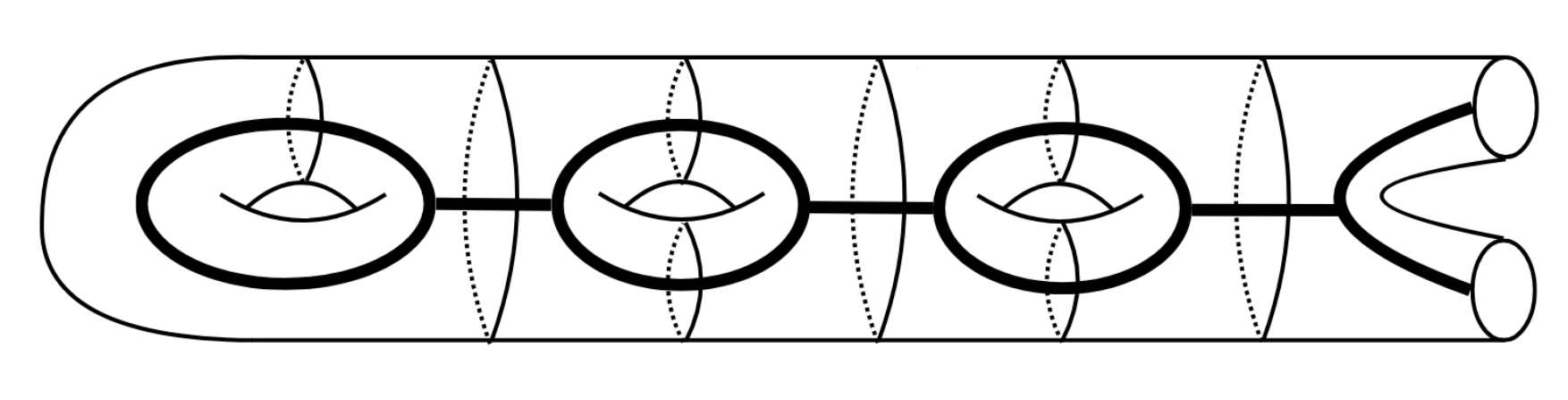}
    \caption{colored quilt decomposition of a surface (I)}
    \label{fig:letterYonsurface}
\end{figure}

\begin{figure}[p]
    \centering
    \includegraphics[width=0.4\linewidth]{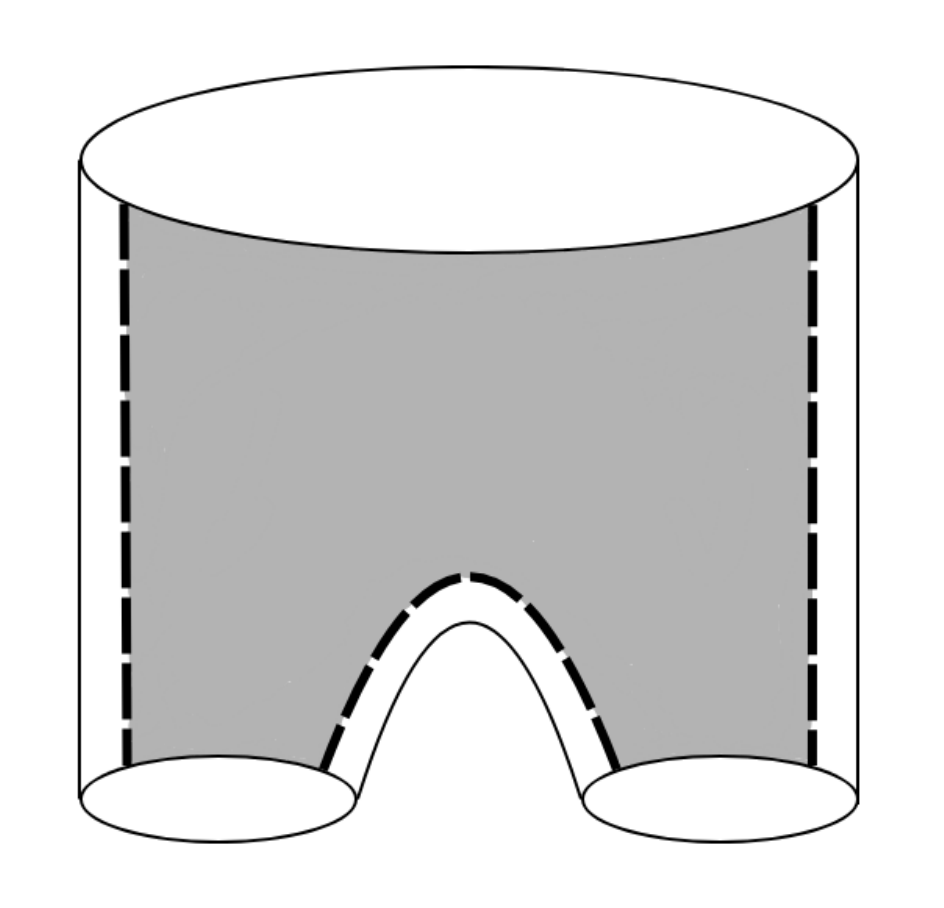}
    \caption{colored quilt decomposition of a pair of pants (II)}
    \label{fig:coloredquiltonpants}
\end{figure}

\begin{figure}[p]
    \centering
    \includegraphics[width=0.9\linewidth]{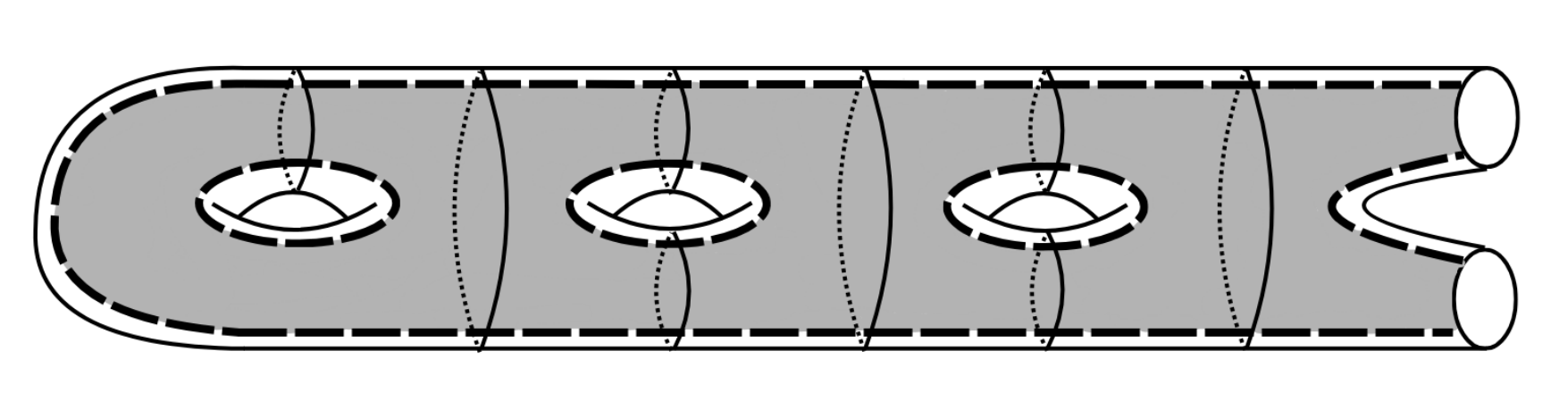}
    \caption{colored quilt decomposition of a surface (II)}
    \label{fig:coloredquiltonsurface}
\end{figure}

\begin{figure}[p]
 \begin{minipage}[p]{0.45\linewidth}
    \centering
    \includegraphics[width=0.75\linewidth]{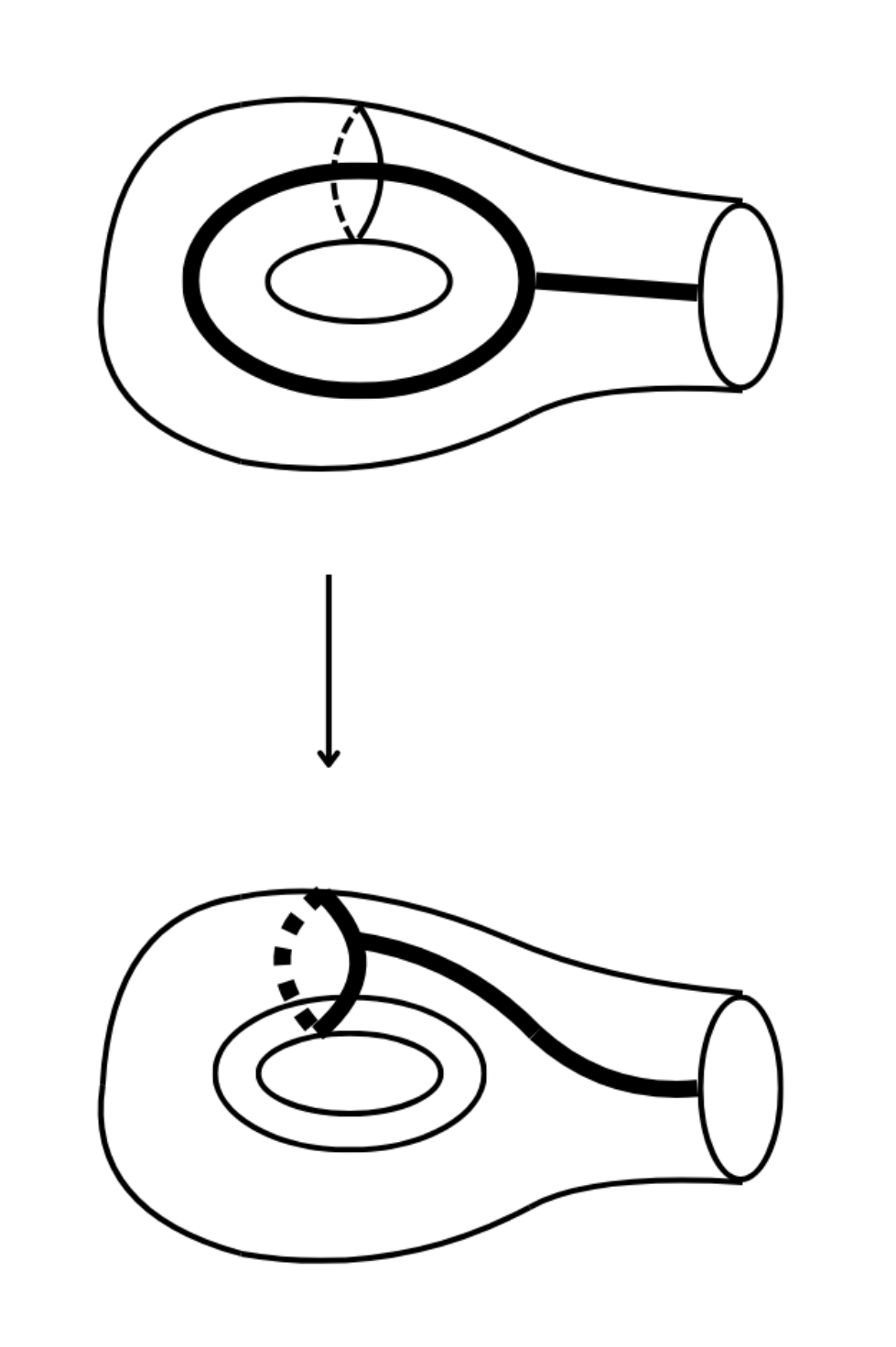}
    \caption{$\overline{S}$-move}
    \label{fig:barS}
 \end{minipage}
 \begin{minipage}[p]{0.45\linewidth}
    \centering
    \includegraphics[width=0.7\linewidth]{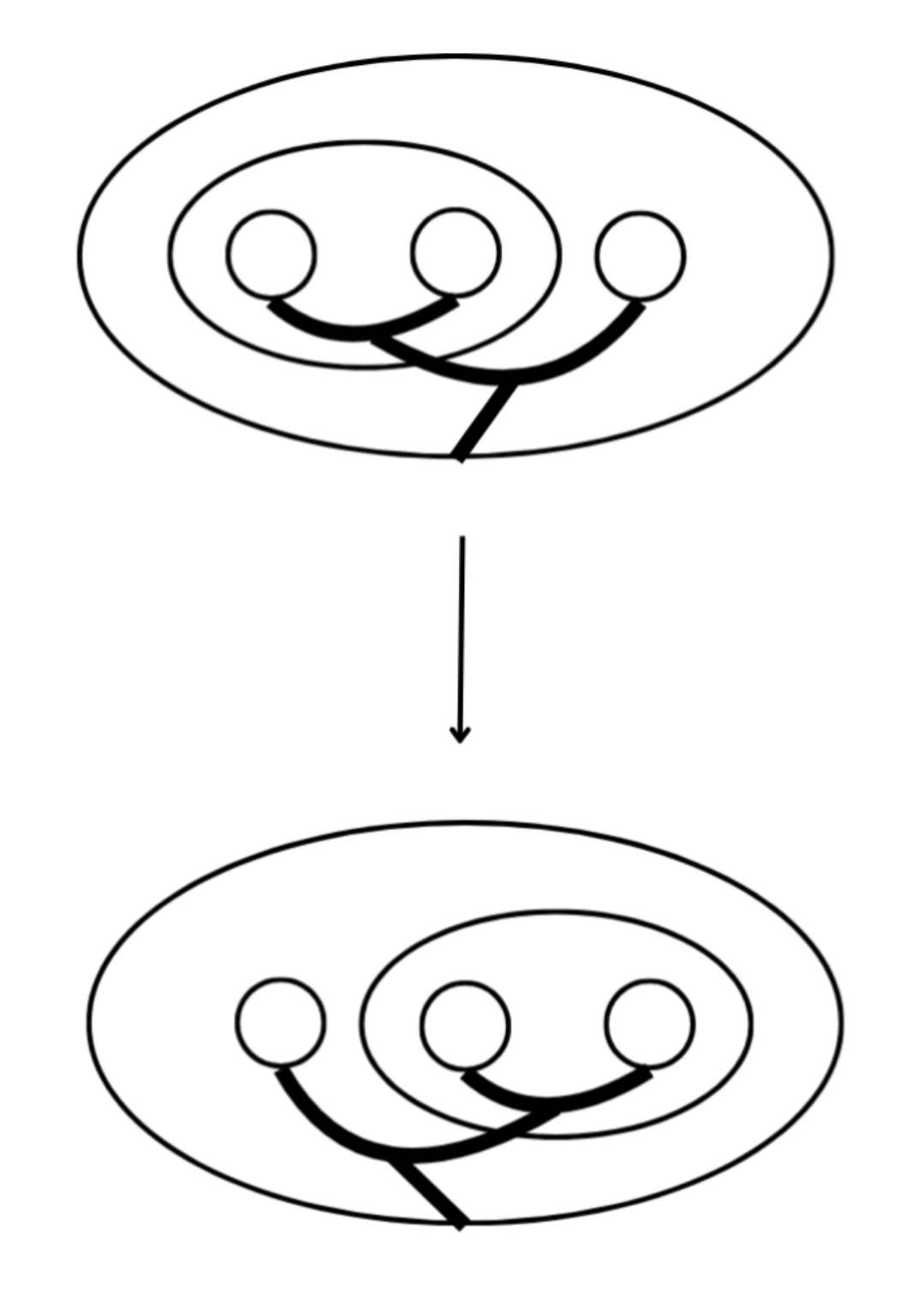}
    \caption{$\overline{A}$-move}
    \label{fig:barA}
 \end{minipage}
\end{figure}
\begin{figure}[p]
 \begin{minipage}[p]{0.45\linewidth}
    \centering
    \includegraphics[width=0.7\linewidth]{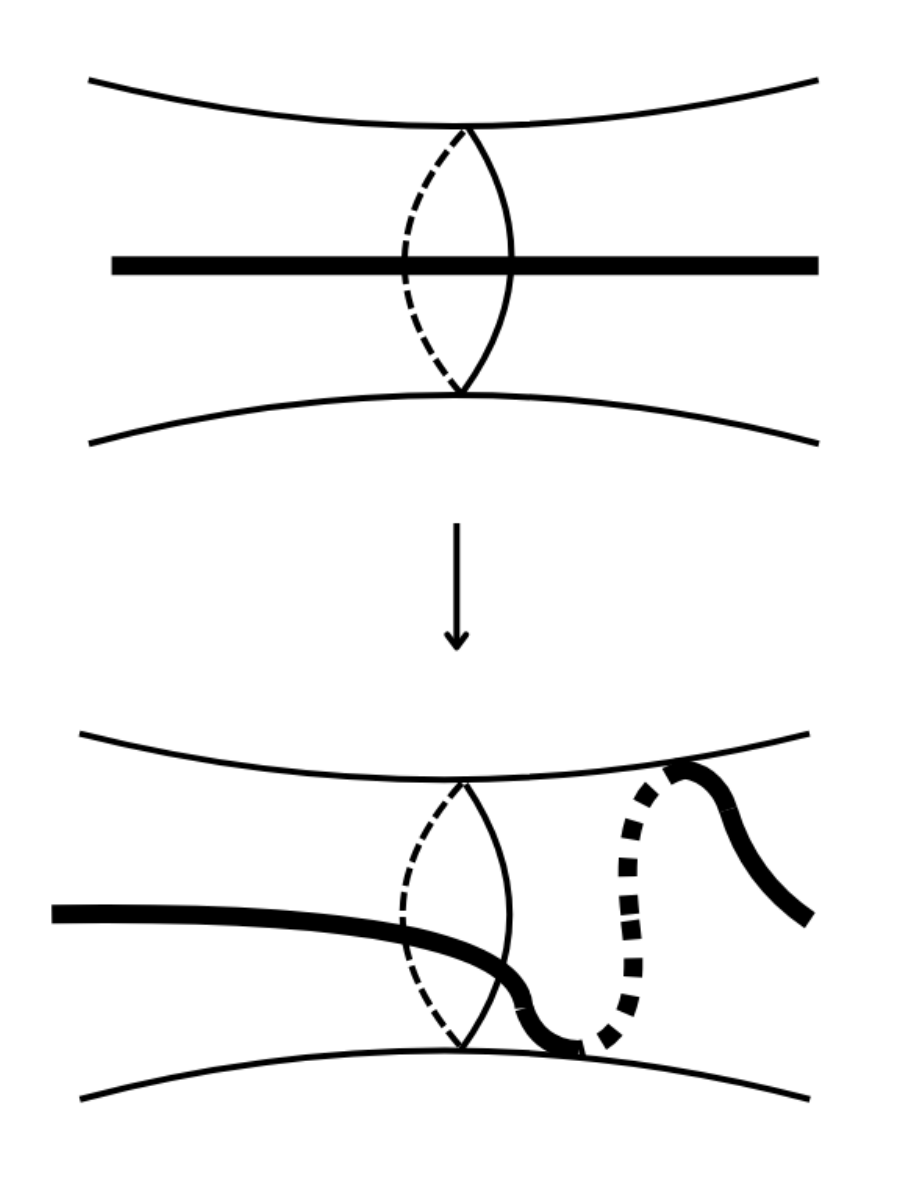}
    \caption{$T$-move}
    \label{fig:T}
 \end{minipage}
 \begin{minipage}[p]{0.45\linewidth}
    \centering
    \includegraphics[width=0.5\linewidth]{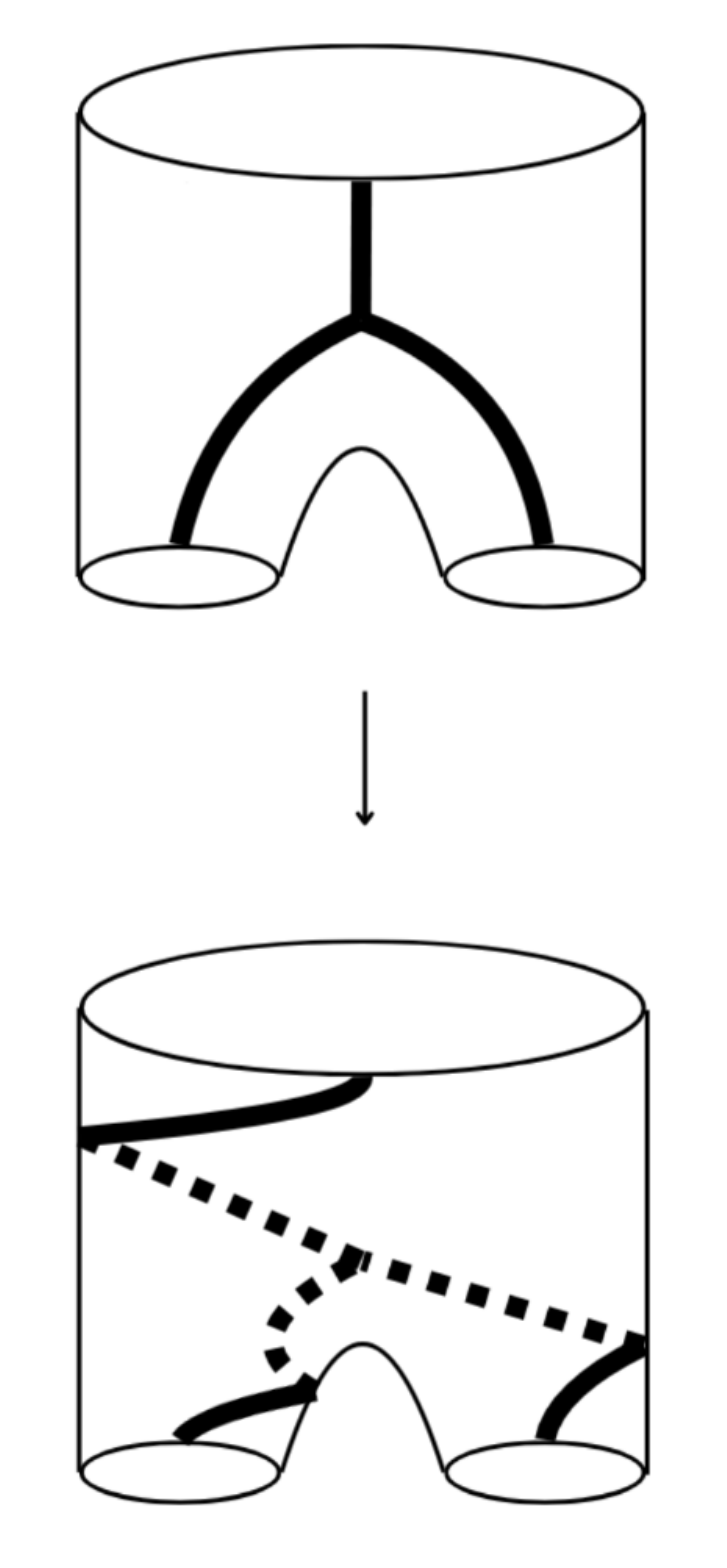}
    \caption{$B'$-move}
    \label{fig:B'}
 \end{minipage}
\end{figure}
\begin{figure}[p]
    \centering
    \includegraphics[width=0.6\linewidth]{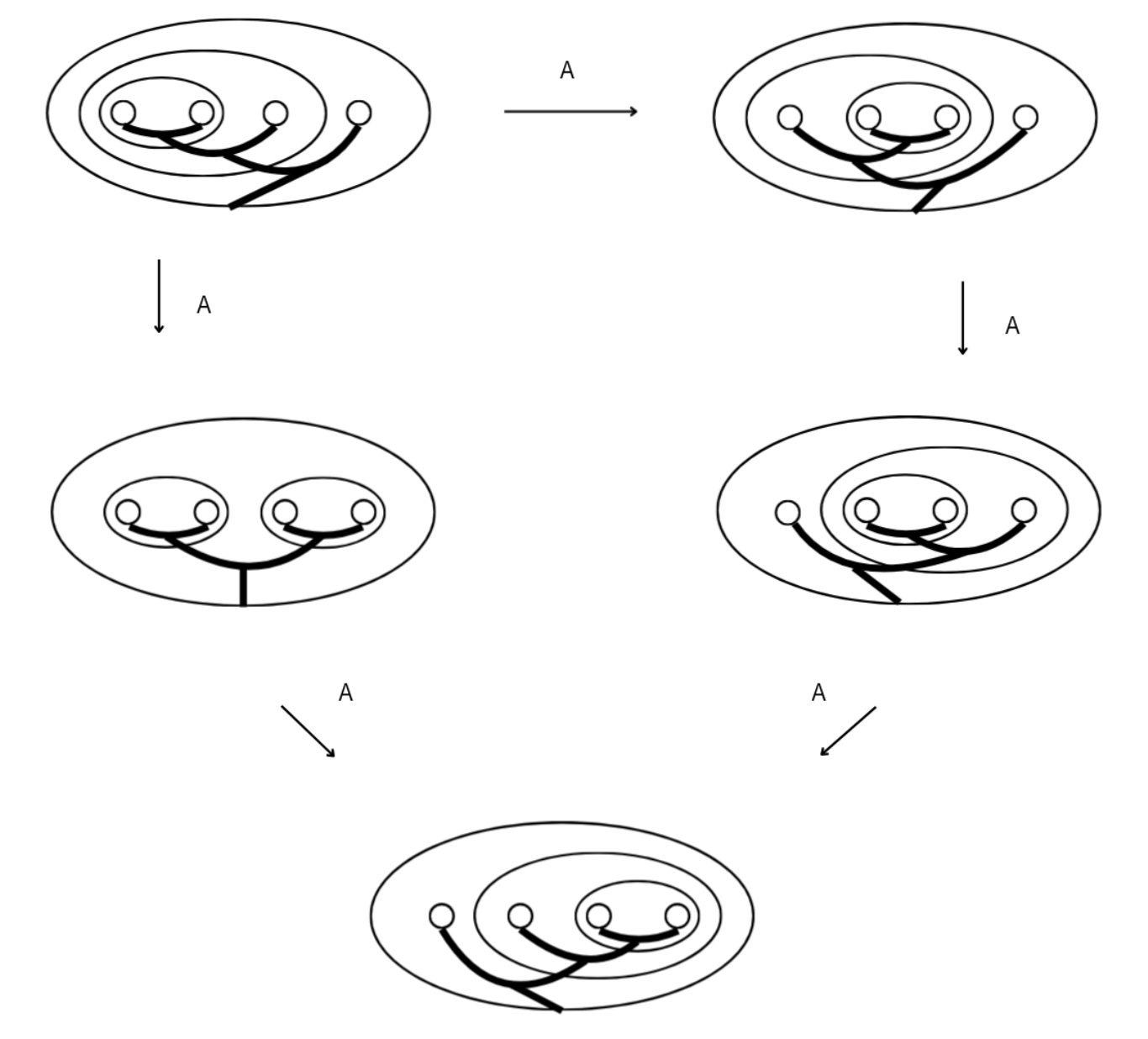}
    \caption{$\overline{5A}$}
    \label{fig:bar5A}
\end{figure}

\begin{figure}[p]
    \centering
    \includegraphics[width=0.55\linewidth]{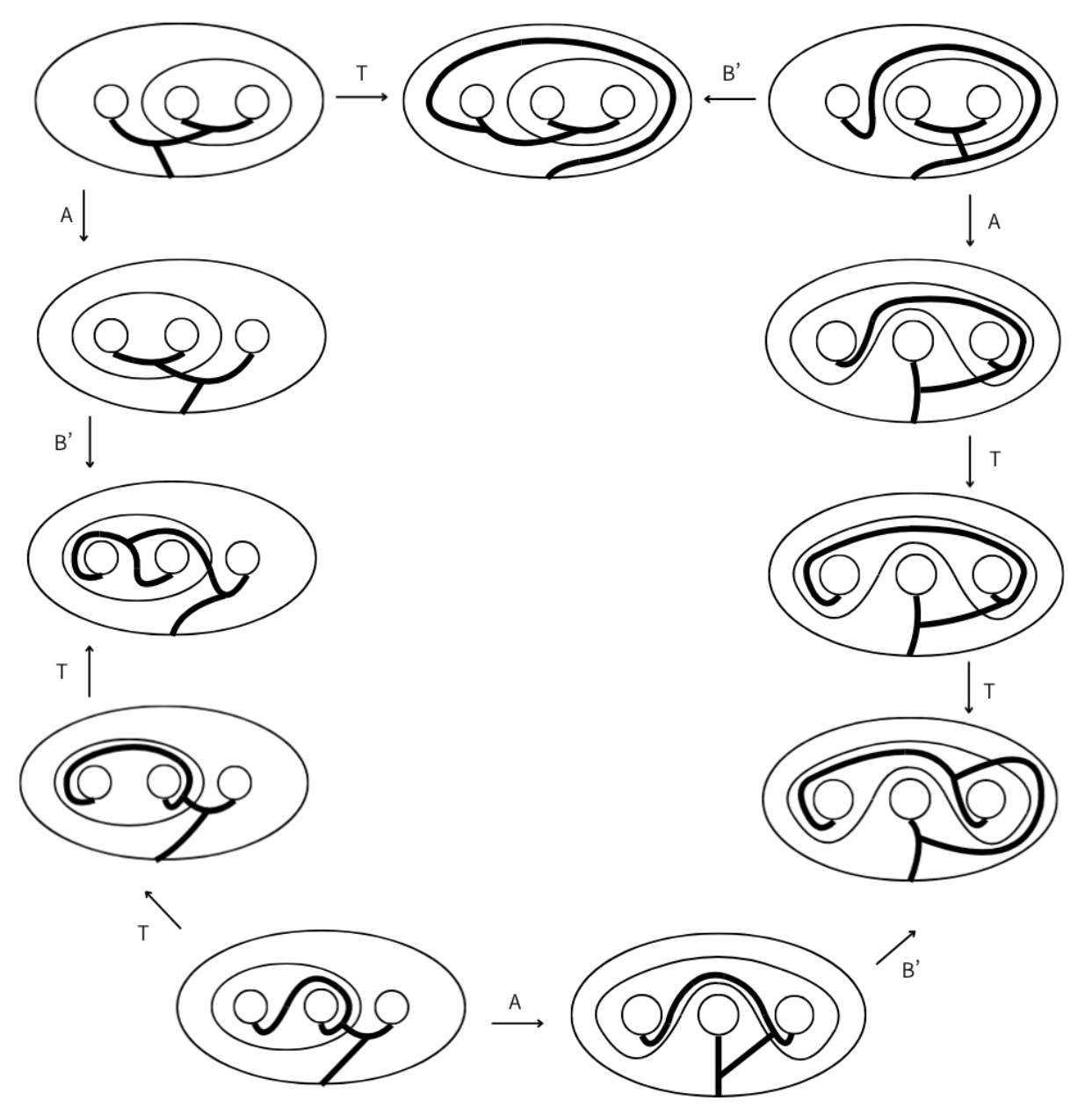}
    \caption{$\overline{3A}$}
    \label{fig:bar3A}
\end{figure}
\begin{figure}[p]
    \centering
    \includegraphics[width=0.4\linewidth]{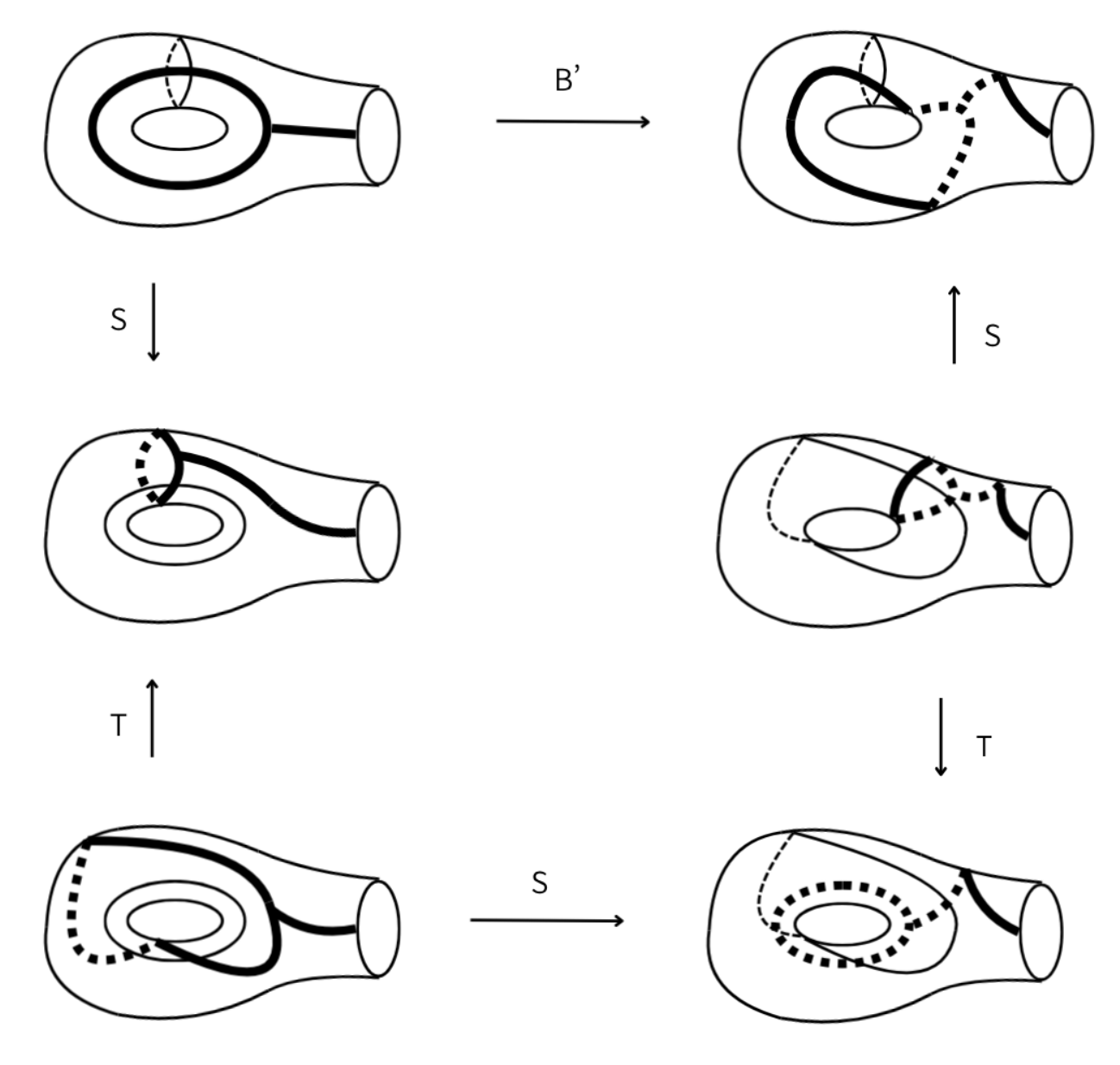}
    \caption{$\overline{3S}$}
    \label{fig:bar3S}
\end{figure}
\begin{figure}[p]
    \centering
    \includegraphics[width=0.6\linewidth]{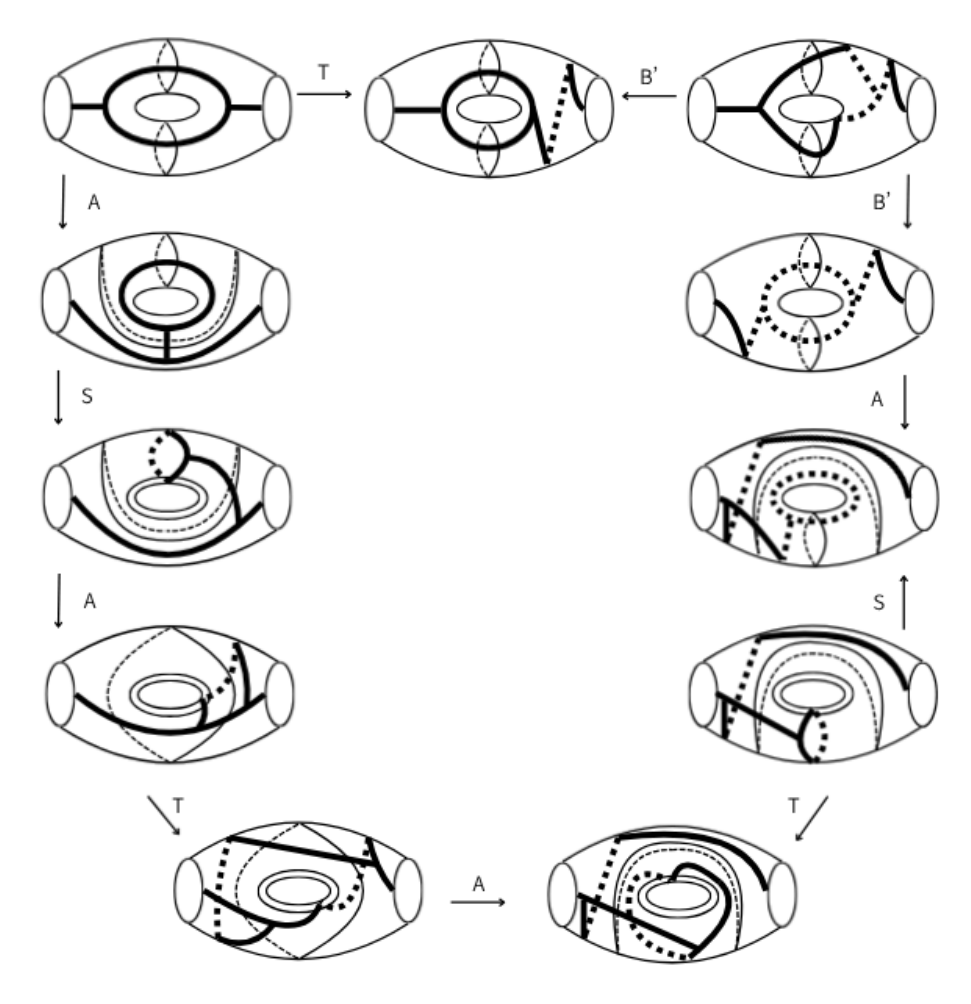}
    \caption{$\overline{6AS}$}
    \label{fig:bar6AS}
\end{figure}
\begin{figure}[p]
    \centering
    \includegraphics[width=0.5\linewidth]{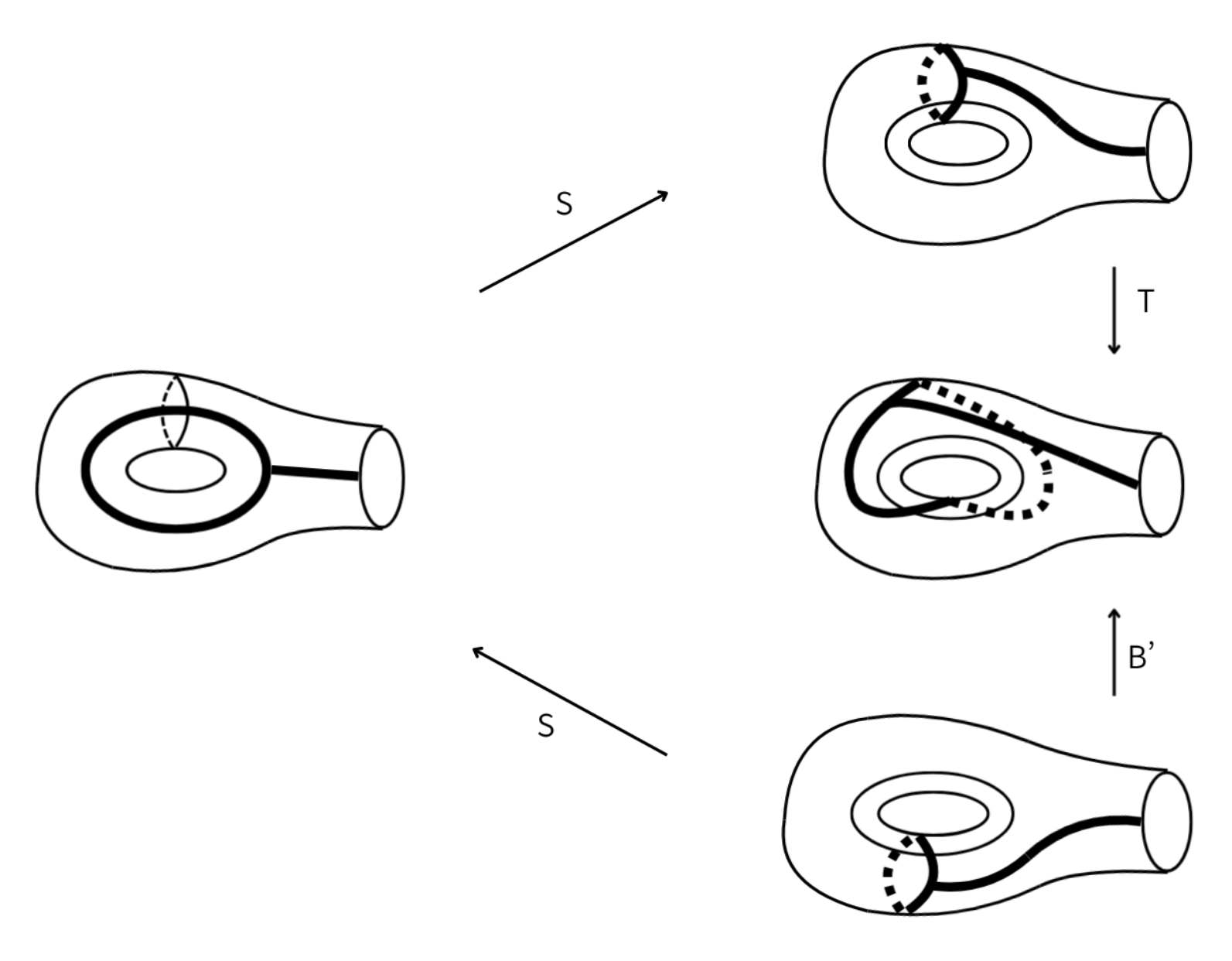}
    \caption{$\overline{2S}$}
    \label{fig:bar2S}
\end{figure}
\begin{figure}[p]
    \centering
    \includegraphics[width=0.5\linewidth]{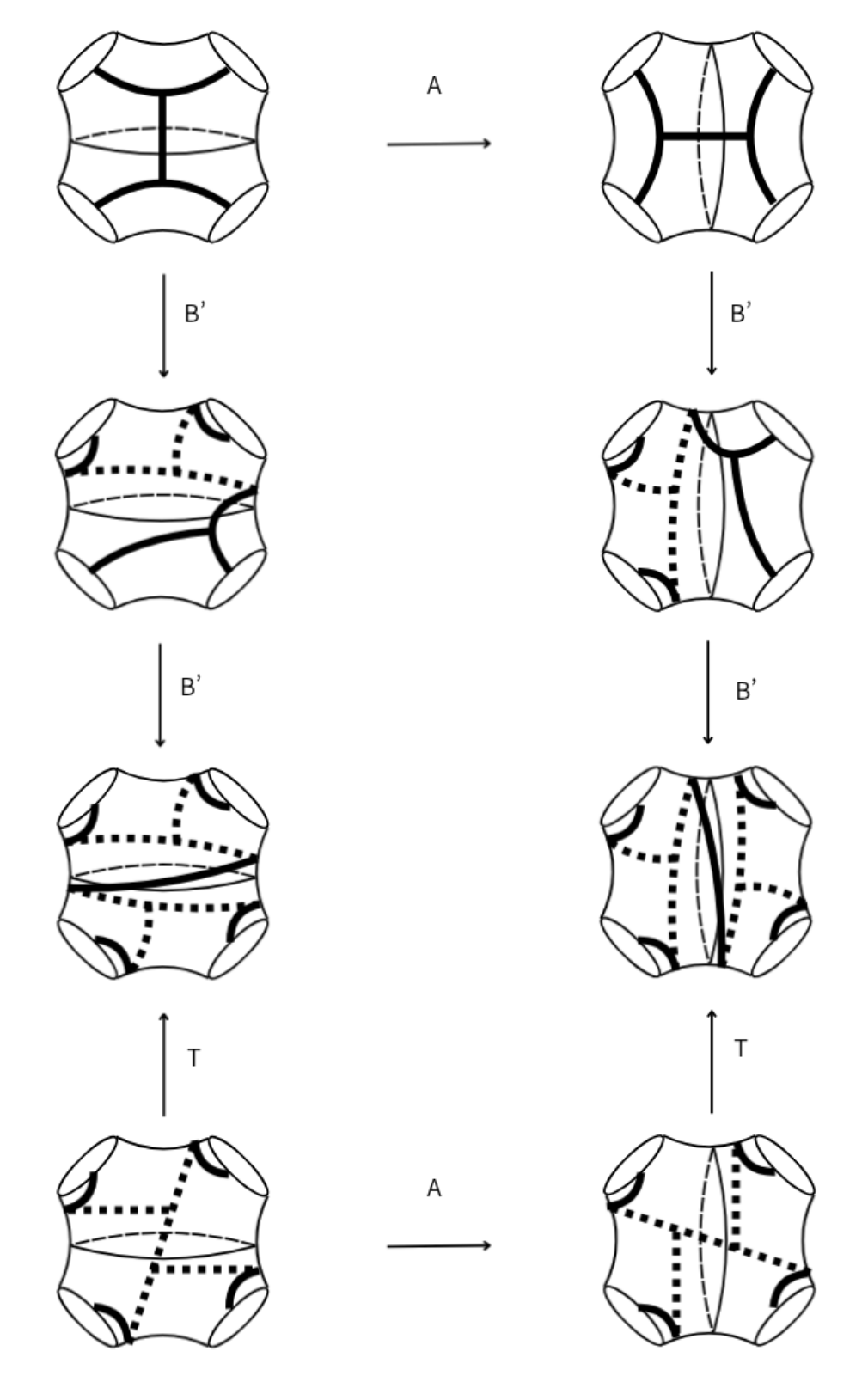}
    \caption{$2A$}
    \label{fig:2A}
\end{figure}
\begin{figure}[p]
    \centering
    \includegraphics[width=0.6\linewidth]{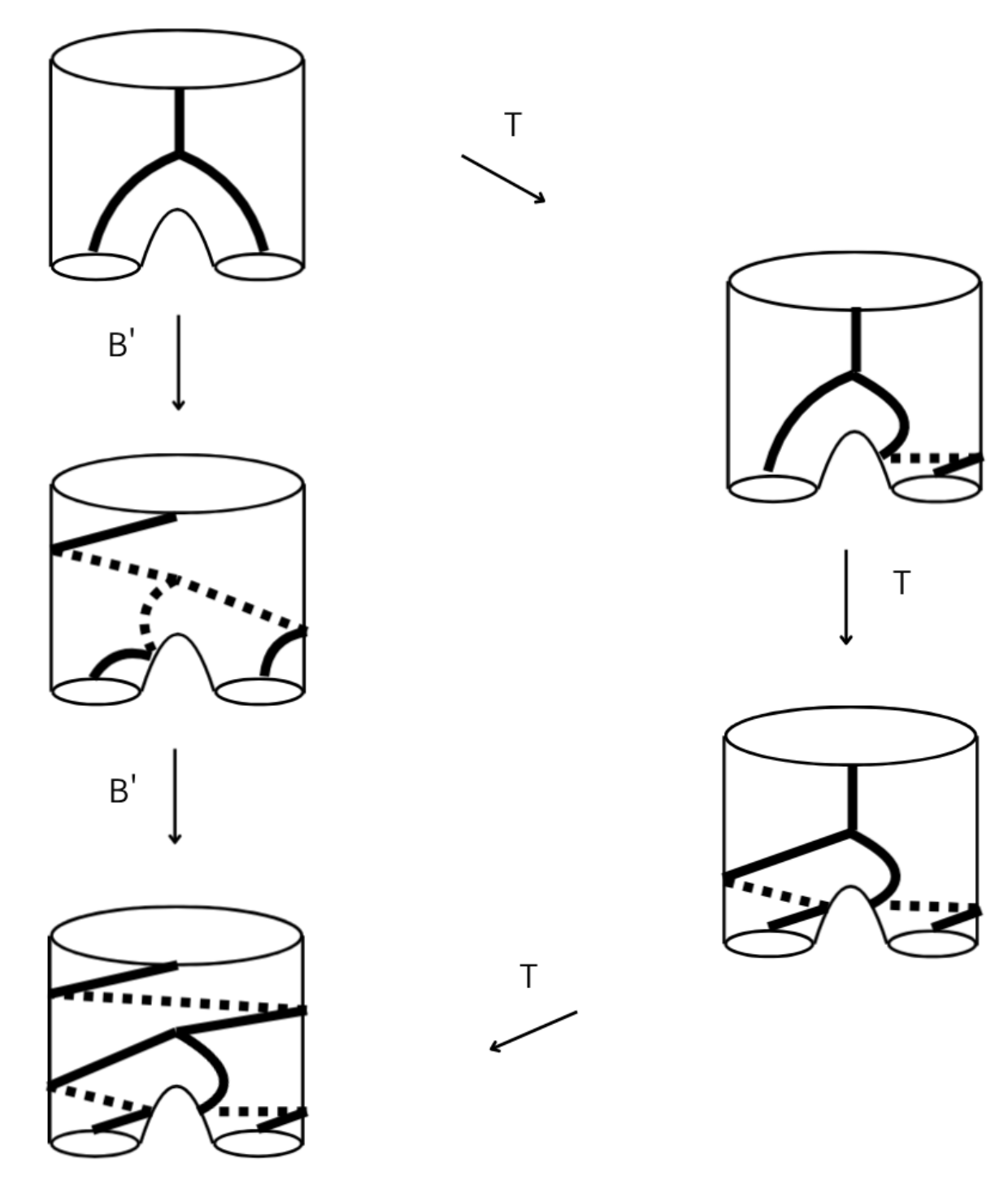}
    \caption{$2B'3T$}
    \label{fig:2B'3T}
\end{figure}

\begin{figure}[p]
    \centering
    \includegraphics[width=0.3\linewidth]{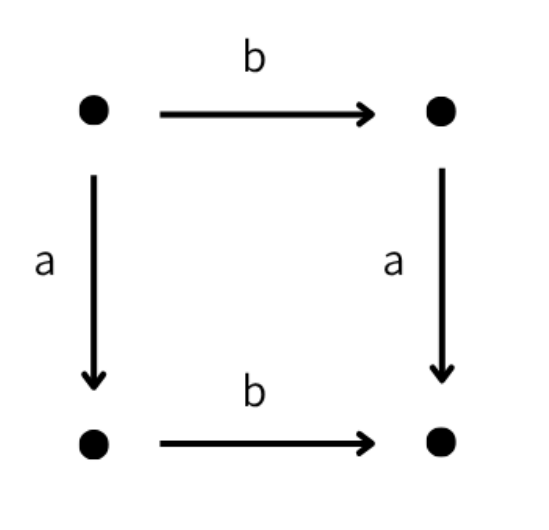}
    \caption{$DC$. Supports of two moves $a$ and $b$ are disjoint.}
    \label{fig:DC}
\end{figure}

\begin{dfn}(a colored quilt decomposition of a surface)[F2]

Let a {\it ``marking" of a pair of pants} be an isotopy class of a structure like that shown in Figure \ref{fig:letterYonpants} on the pair of pants. It is a ``letter Y" embedded on the pair of pants, with the ``letter Y" having three tails, each combined to one of the three boundaries of the pair of pants. 
Let a {\it ``marking" of a surface} $\Sigma_{g,r}$ be an isotopy class of a pants decomposition where each pair of pants has a ``marking" of a pair of pants, and each ``marking" of a pair of pants is connected on each SCC of the pants decomposition.   
We call the ``marking" a {\it colored quilt decomposition} of $\Sigma_{g,r}$ (see Figure \ref{fig:letterYonsurface}).

We can also define it as follows. Let a colored quilt decomposition of a pair of pants be a quilt decomposition of it such that one of two patches is colored (see Figure \ref{fig:coloredquiltonpants}). And, let a colored quilt decomposition of $\Sigma_{g,r}$ be a quilt decomposition of $\Sigma_{g,r}$ where each pair of pants has a colored quilt decomposition of the pair of pants, and each color of the two colored quilt decompositions of the pairs of pants next to each other is connected on each SCC of the pants decomposition (see Figure \ref{fig:coloredquiltonsurface}).

It is easy to see that the two definitions of colored quilt decomposition are equivalent in the sense that a trivalent graph is homotopy equivalent to a colored quilt on $\Sigma_{g,r}$.

\end{dfn}

\begin{dfn}(colored extended Hatcher complex)[F2]

The {\it colored extended Hatcher complex} $\overline{H}_{g,r}$ is a two-dimensional cell complex having these properties:

Each 0-cell is a colored quilt decomposition of the surface $\Sigma_{g,r}$. Each 1-cell is either an $\overline{S}$-move (Figure \ref{fig:barS}), 
an $\overline{A}$-move (Figure \ref{fig:barA}), $T$-move (Figure \ref{fig:T}),or $B'$-move (Figure \ref{fig:B'}). Each 2-cell is either type  $\overline{5A}$ 
 (Figure \ref{fig:bar5A}),
 $\overline{3A}$ (Figure \ref{fig:bar3A}), $\overline{3S}$ (Figure \ref{fig:bar3S}), 
 $\overline{6AS}$ (Figure \ref{fig:bar6AS}), $\overline{2S}$ (Figure \ref{fig:bar2S}), $2A$ (Figure \ref{fig:2A}), $2B'3T$ (Figure \ref{fig:2B'3T}),
 or $DC$ (Figure \ref{fig:DC}).

We will sometimes write them without bars.

\end{dfn}

These figures mean that they are subsurfaces and their complements are invariant under the moves. We call the subsurfaces {\it supports} of the moves. 
In Figure \ref{fig:DC}, supports of two moves $a$ and $b$ are disjoint.

\begin{thm}

\textup{[F2]}

A colored extended Hatcher complex $\overline{H}(\Sigma_{g,r})$ is simply-connected.

\end{thm}

Let $\mathcal{M}_{g,r}$ be a mapping class group for an orientable surface $\Sigma_{g,r}$. 
$\mathcal{M}_{g,r}$ acts on 
$\overline{H}(\Sigma_{g,r})$ cellularly.

\begin{thm}
\textup{[F2]}
      \mbox{}

A torsion-free subgroup of $\mathcal{M}_{g,r}$ acts freely on the colored extended Hatcher complex $\overline{H}(\Sigma_{g,r})$.

\end{thm}

Let $r=0,1$. {\it Torelli group $\mathcal{I}_{g,r}$} is a subgroup of $\mathcal{M}_{g,r}$ such that acts on $H_1(\Sigma_{g,r}, \mathbb{Z})$ trivially
 (See [J]).
 
The Torelli groups are torsion-free ([I]). Therefore, the following corollary holds.

\begin{crl}
\textup{[F2]}
    \mbox{}
    
    $\mathcal{I}_{g,r}$ acts freely on the colored extended Hatcher complex $\overline{H}(\Sigma_{g,r})$.

\end{crl}

In [F2], we proved theorem 5 using the fact that the groups are torsion-free ((1) of theorem 18 of [F2]).
In this paper,
we prove corollary 6 without using the fact that the Torelli groups are torsion-free. The proof is given below.

Corollary 6 is equivalent to the following proposition.

\begin{prop}\label{thm:mainprop}

    For any $x\in\overline{H}(\Sigma_{g,r})$ and any $f\in\mathcal{I}_{g,r}$, the following holds.
    \[
    f(x)=x \Rightarrow f=1
    \]
\end{prop}

\begin{proof}[Proof of Proposition \ref{thm:mainprop}]

We will prove proposition \ref{thm:mainprop} by considering the following cases: (1) when $x$ is a 0-cell of $\overline{H}(\Sigma)$,
 (2) when $x$ is an inner point of 1-cell of $\overline{H}(\Sigma)$, 
 and (3) when $x$ is an inner point of 2-cell of $\overline{H}(\Sigma)$.

The assumption that $r=0,1$ is not essential to the following discussion.

(1) Let $x$ be a 0-cell of $\overline{H}(\Sigma)$. Then, $x$ is a colored quilt decomposition of $\Sigma$. 
Consider the colored quilt as the trivalent graph $\Gamma$ with $r$ tails (see Figure \ref{fig:oldFig.39}).

\begin{figure}[ht]
    \centering
    \includegraphics[width=0.6\linewidth]{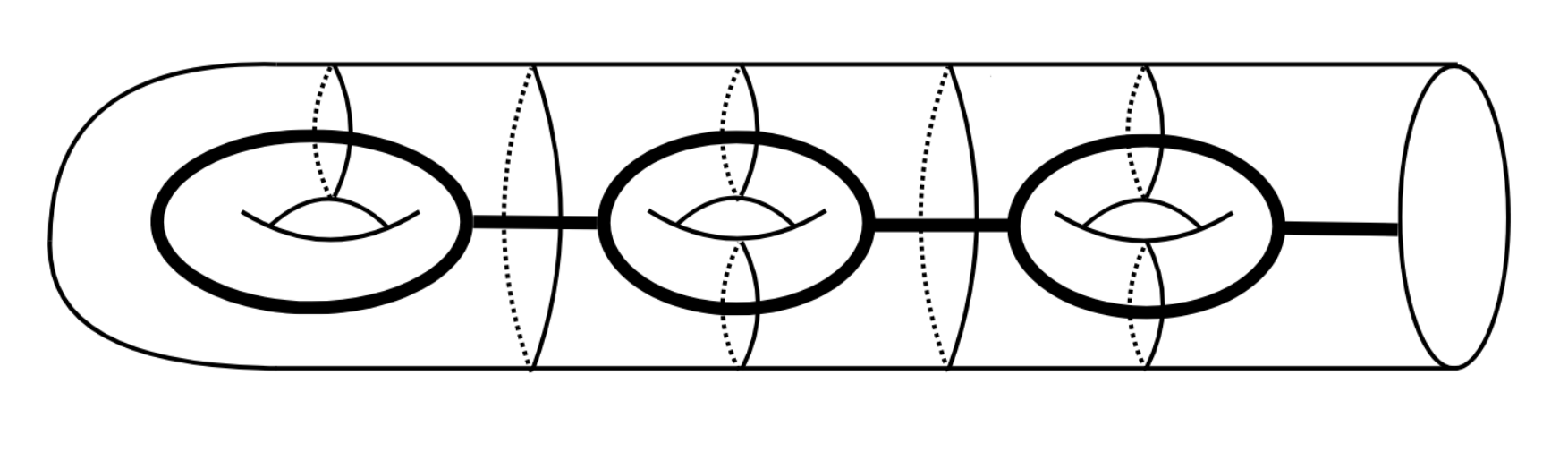}
    \caption{}
    \label{fig:oldFig.39}
\end{figure}

Because $f(x)=x$, one can consider that $f|_{\Gamma}$ is an automorphism of $\Gamma$.

\begin{claim}\label{thm:claim1}
    $ f|_{\Gamma}=1 \Rightarrow f=1$.
\end{claim}

\begin{proof}

The vertices and edges of $\Gamma$ correspond to the pairs of pants and 
 the SCCs of the pants decomposition.
An edge of $\Gamma$ intersects at 1 point to the corresponding SCC of the pants decomposition.
Since  $f|_{\Gamma}=1$, 

\[
f|_{\Gamma\cup\{\text{SCC of the pants decomposition\}}}=1
\]

The complement of $\Gamma\cup\{\text{SCC of the pants decomposition\}}$ on $\Sigma$ is a disjoint union 
 of $2g-2+r$ contractible regions which are each 9-gons.
Hence, $f=1$.

\end{proof}

\begin{figure}[ht]
    \centering
    \includegraphics[width=0.6\linewidth]{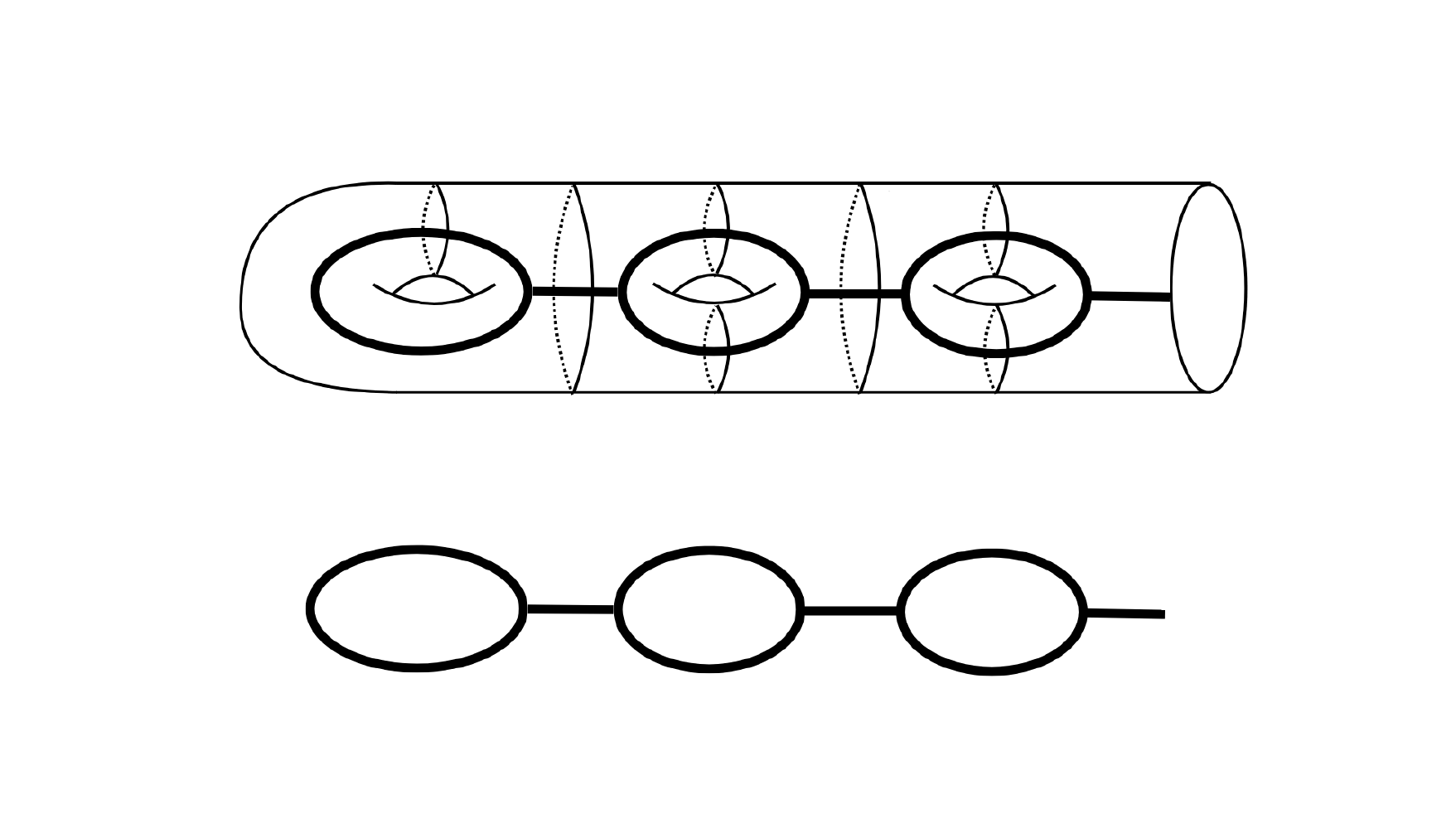}
    \caption{}
    \label{fig:oldFig.40}
\end{figure}

\begin{claim}
    The first Betti number of $\Gamma$ is $g$.
\end{claim}

\begin{proof}

Let $x_0$ be the colored quilt decomposition shown in the upper part of Figure \ref{fig:oldFig.40}.
If $x=x_0$, $\Gamma$ corresponding to $x$ is as shown in the lower part of Figure \ref{fig:oldFig.40}, so the first Betti number of $\Gamma$ is clearly $g$.
(Figure \ref{fig:oldFig.40} shows an example where $g=3$, $r=1$.)

The general $x$ is connected to $x_0$ via only moves of 4 types, $S$, $A$, $T$ or $B'$ in $\overline{H}(\Sigma)$, because $\overline{H}(\Sigma)$ is connected.
So, it is enough to check that the first Betti number of $\Gamma$ is invariant under such moves.

Under $S$ or $T$-move, $\Gamma$ is invariant.

Under $B'$-move, the cyclic order of the three edges extending from the vertex induced by an 
orientation of $\Sigma$ is changed (See Figure \ref{fig:oldFig.41(a)}). The first Betti number is invariant.

\begin{figure}[ht]
    \centering
    \includegraphics[width=0.7\linewidth]{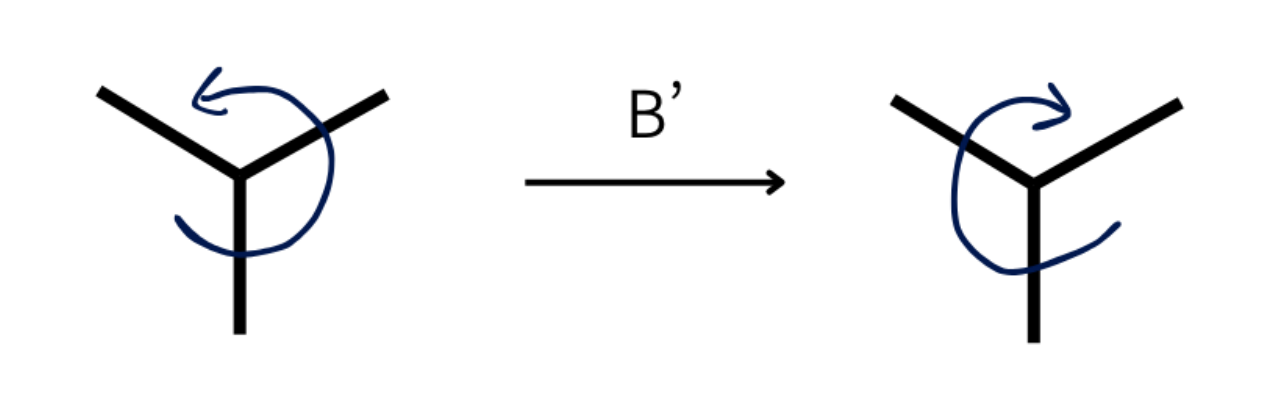}
    \caption{}
    \label{fig:oldFig.41(a)}
\end{figure}

Under $A$-move, $\Gamma$ is changed as shown in Figure \ref{fig:oldFig.41(b)}. The first Betti number is invariant.

\begin{figure}[ht]
   \centering
   \includegraphics[width=0.7\linewidth]{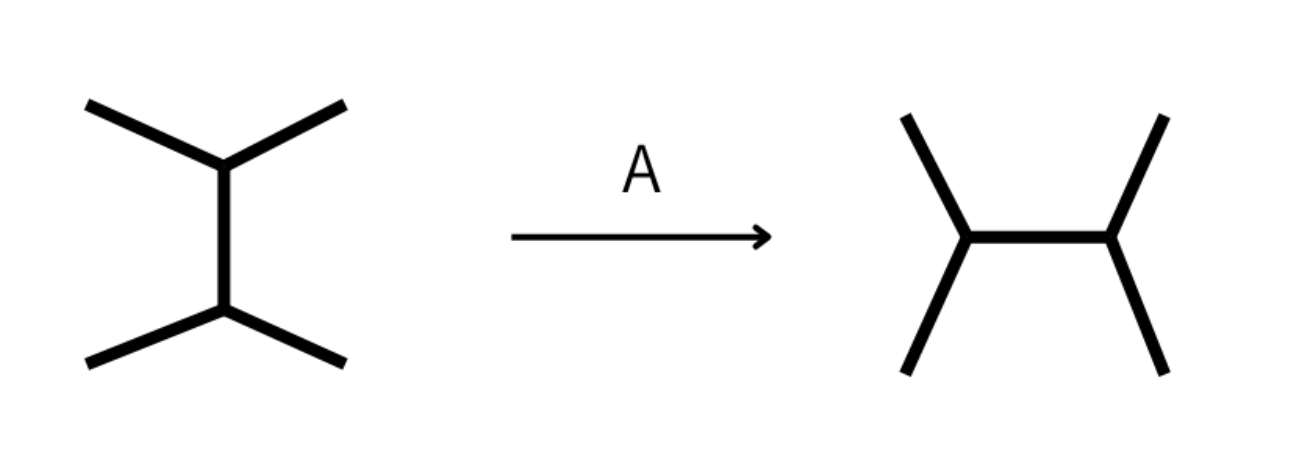}
   \caption{}
   \label{fig:oldFig.41(b)}
\end{figure}
    
\end{proof}

Denote that $\Gamma_x$ is a graph $\Gamma$ corresponding to a colored quilt decomposition $x$. Let $e_x: \Gamma_x\to\Sigma$ be the embedding for $x$.
Let $e_x^{*}:H_1(\Gamma_x)\to H_1(\Sigma)$ be the map of 1-homology induced by $e_x$.
$e_x^{*}$ is homomorphism clearly.

\begin{claim}\label{thm:injectiveclaim}
    $e_x^{*}$ is injective.
\end{claim}
\begin{proof}

Let $x_0$ be the colored quilt decomposition shown in the upper part of Figure \ref{fig:oldFig.42}.
If $x=x_0$, $\Gamma_x$ is as shown in the lower part of Figure \ref{fig:oldFig.42}.
Take a basis of $H_1(\Gamma_x)$ as shown in the lower part of Figure \ref{fig:oldFig.42}.
Then, $e_x^{*}$ is injective clearly.

\begin{figure}[ht]
    \centering
    \includegraphics[width=0.6\linewidth]{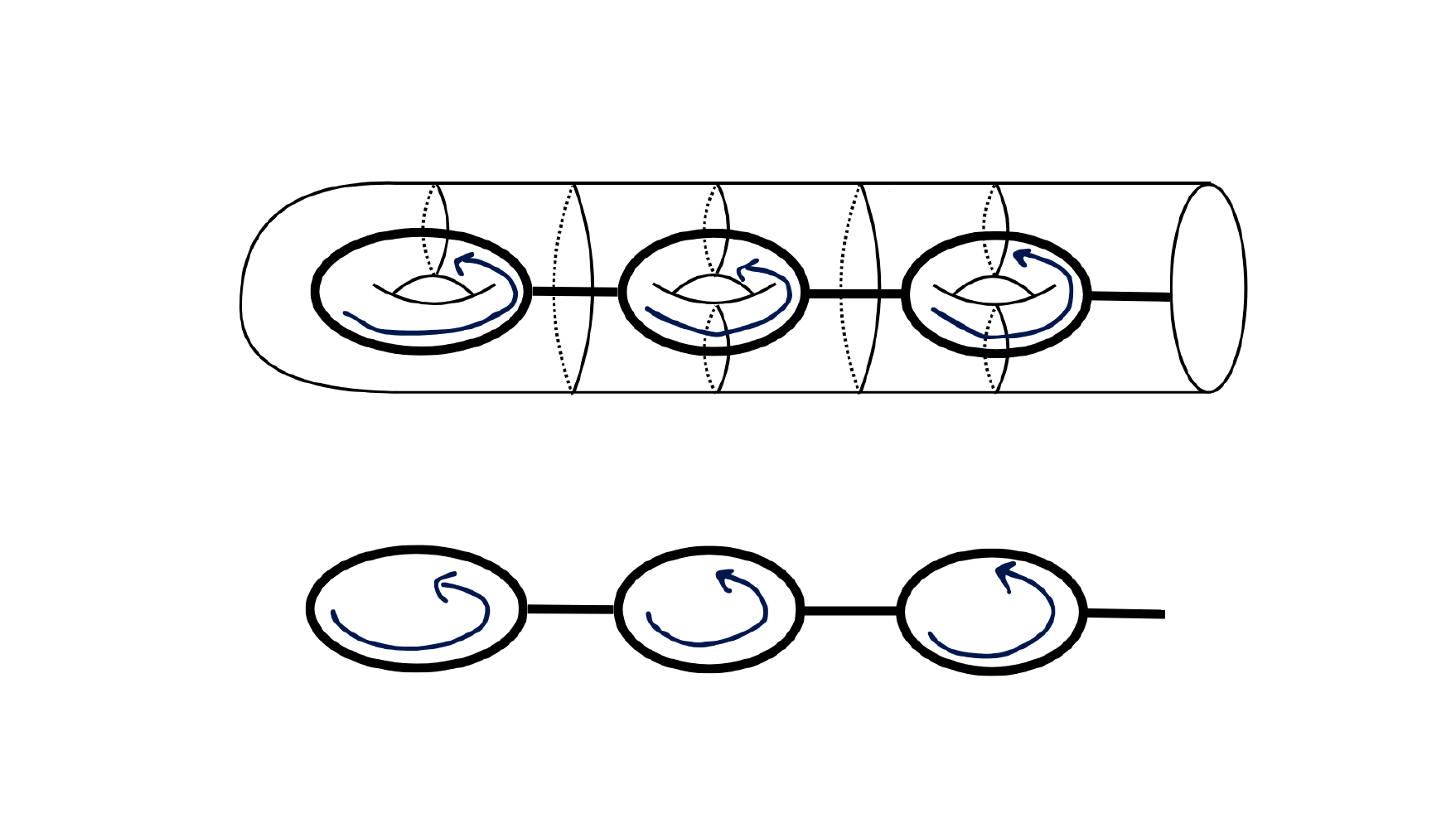}
    \caption{}
    \label{fig:oldFig.42}
\end{figure}


Let $x$ be general 0-cell of $\overline{H}(\Sigma)$. 
Take an edge path $m_n\circ m_{n-1}\circ\dots\circ m_1$ in $\overline{H}(\Sigma)$ such that $x=(m_n\circ m_{n-1}\circ\dots\circ m_1)x_0$ 
 where  $m_i=S,A,T\text{ or }B'$ for $i=1, \dots , n$, i.e. we obtain $x$ by acting $m_n\circ m_{n-1}\circ\dots\circ m_1$ to $x_0$.
 
We want to consider how to take a basis of $H_1(\Gamma_x)$ for $x$.

\begin{figure}[ht]
    \centering
    \includegraphics[width=0.4\linewidth]{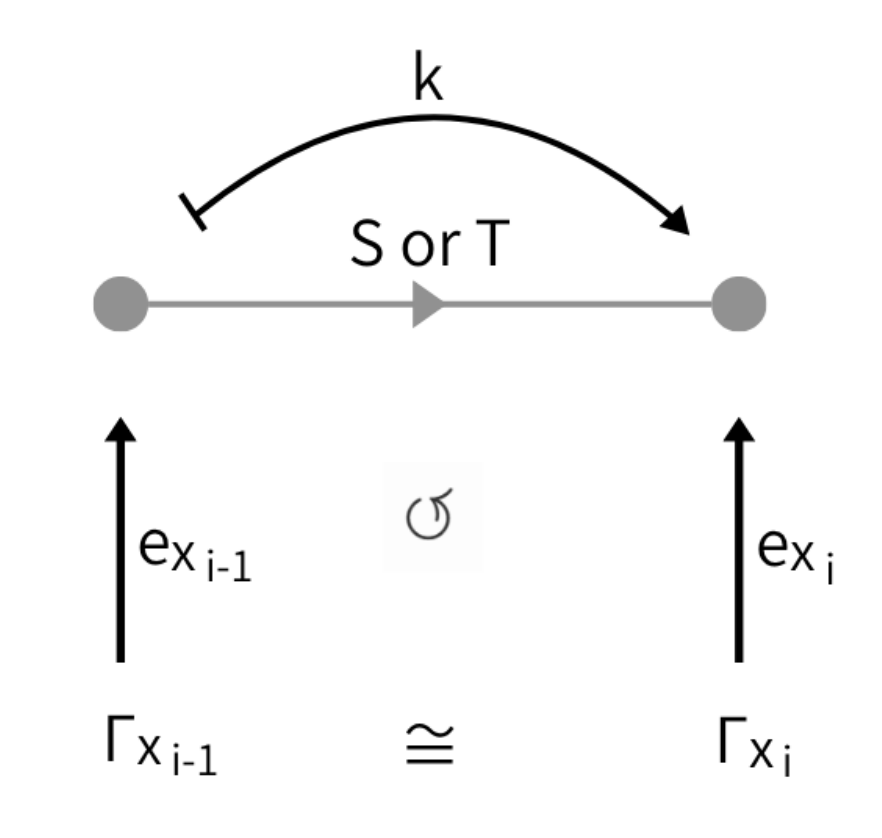}
    \caption{}
    \label{fig:oldFig.43}
\end{figure}

Let $x_{i-1},x_i$ be colored quilt decompositions such that $x_i=m_i x_{i-1}$ 
and suppose that there is a determined way to take the basis of $H_1(\Gamma_{x_{i-1}})$
such that $e_{x_{i-1}}^{*}$ is injective.

\begin{figure}[ht]
    \centering
    \includegraphics[width=0.7\linewidth]{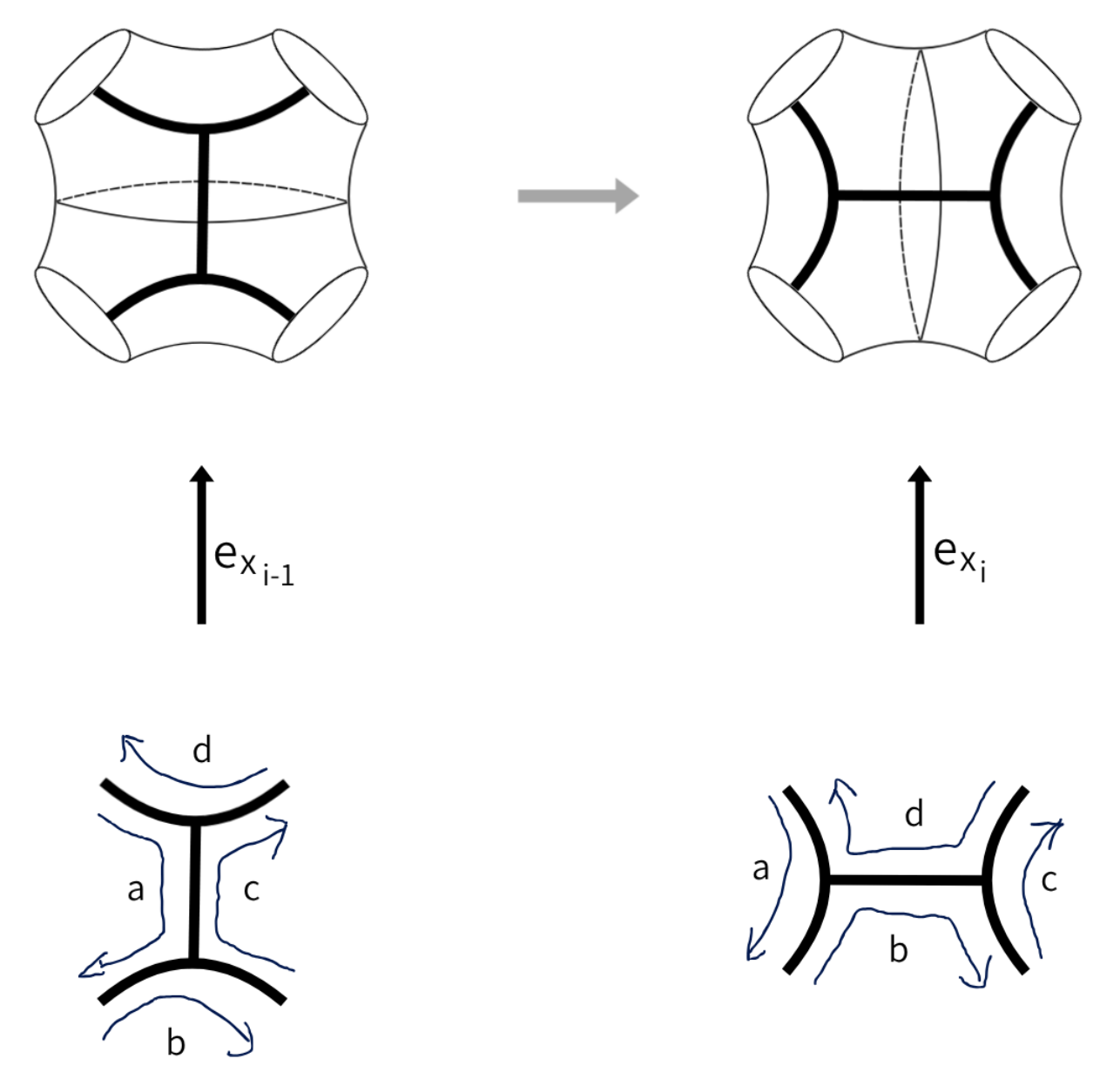}
    \caption{}
    \label{fig:oldFig.44(a)}
\end{figure}
\begin{figure}[ht]
    \centering
    \includegraphics[width=0.7\linewidth]{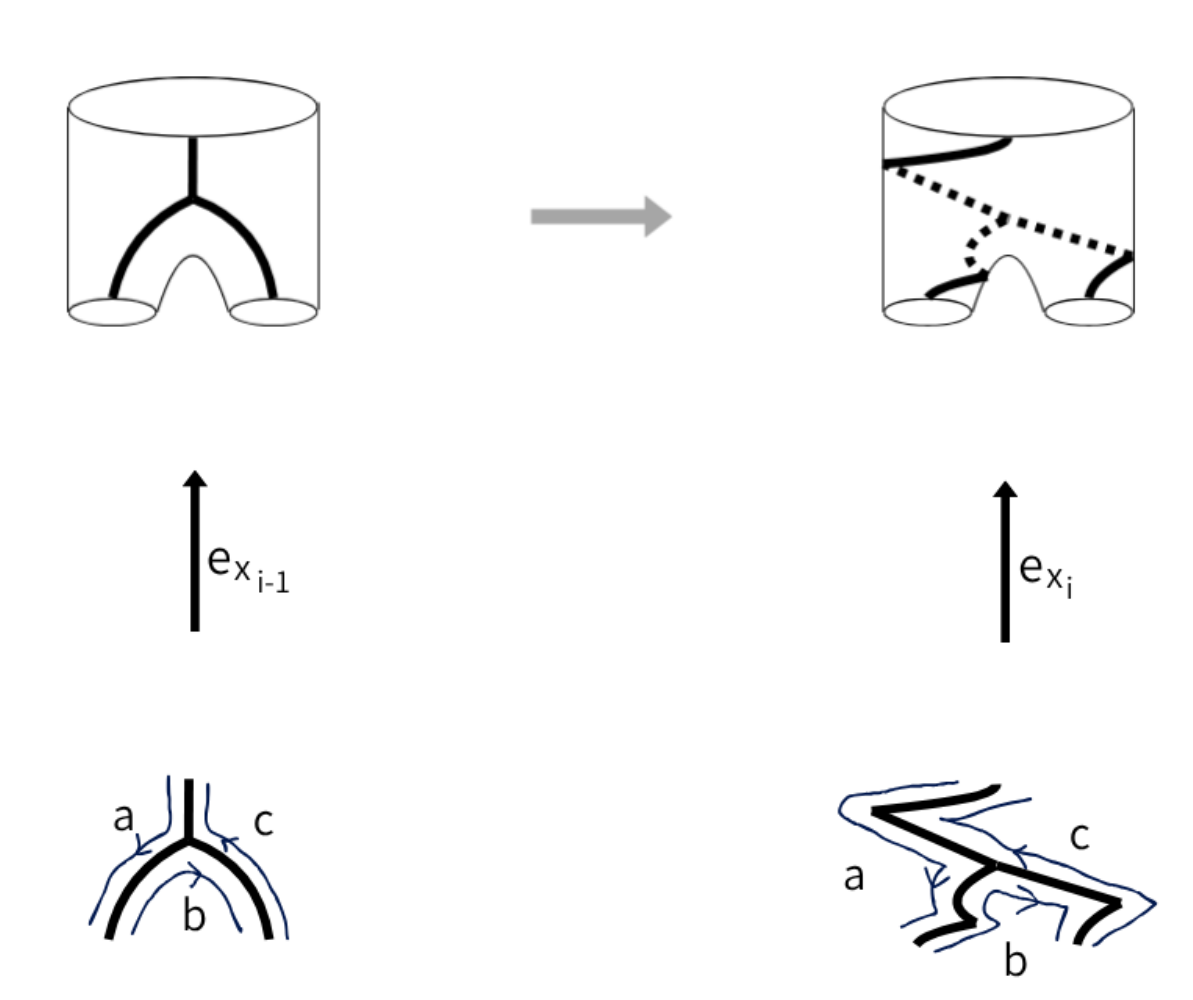}
    \caption{}
    \label{fig:oldFig.44(b)}
\end{figure}

\clearpage

We will determine how to take a basis of $H_1(\Gamma_{x_i})$ as follows.

Suppose that $m_i=S \text{ or } T$. Then, there is $k \in \mathcal{M}_{g,r}$ such that $x_i=k(x_{i-1})$, i.e. if $m_i=S$, then $k$ is the 
``90-degree rotation" in the neighborhood of the support of the move, and if $m_i=T$, then $k$ is 
the Dehn twist in the neighborhood of the support of the move.
(Remark that when we write $m_i x_{i-1}$, we think $x_{i-1}$ to be a 0-cell of $\overline{H}(\Sigma)$, 
while when we write $k(x_{i-1})$, we think $x_{i-1}$ to be a colored quilt decomposition on $\Sigma$.)
Then, we define an isomorphism from $\Gamma_{x_{i-1}}$ to $\Gamma_{x_i}$ as $h:=e_{x_i}^{-1} \circ k \circ e_{x_{i-1}}$
 (See Figure \ref{fig:oldFig.43}). One can see that $h$ is well-defined.

Let $h^{*}:H_1(\Gamma_{x_{i-1}})\to H_1(\Gamma_{x_i})$ be the map of 1-homology induced by $h$.
Take a basis of $H_1(\Gamma_{x_i})$ as $h^{*}(\text{the basis of } H_1(\Gamma_{x_{i-1}}))$.
Then, $\{e_{x_i}^*(\text{the basis of } H_1(\Gamma_{x_i}))\}\subset \{\text{a basis of } H_1(\Sigma)\}$, 
so $e_{x_i}^*$ is injective.

If $m_i=A$, it is enough to take a basis of $H_1(\Gamma_{x_i})$ as Figure \ref{fig:oldFig.44(a)}. In the figure, $a,b,c,d$ are examples of 1-cycle of the graphs.

Similarly, If $m_i=B'$, it is enough to take a basis of $H_1(\Gamma_{x_i})$ as Figure \ref{fig:oldFig.44(b)}. In the figure, $a,b,c$ are examples of 1-cycle of the graphs.

From the above discussion, there exists a basis of $H_1(\Gamma_x)$
such that 

\[\{e_{x}^*(\text{the basis of } H_1(\Gamma_{x}))\}\subset \{\text{a basis of } H_1(\Sigma)\}\]

So, $e_x^*$ is injective.

\end{proof}

$f$ acts on $H_1(\Sigma)$ trivially because $f\in\mathcal{I}_{g,r}$.
From claim \ref{thm:injectiveclaim}, the following claim holds.

\begin{claim}\label{thm:claim2}
$f|_{\Gamma}$ acts on $H_1(\Gamma)$ trivially.
    
\end{claim}


Assume $f|_{\Gamma}\ne 1$. From claim \ref{thm:claim2}, there exists a closed path $l$ of $\Gamma$ that is ``rotated" by $f|_{\Gamma}$ (See Figure \ref{fig:oldFig.45(a)}).

\begin{figure}[ht]
    \centering
    \includegraphics[width=0.2\linewidth]{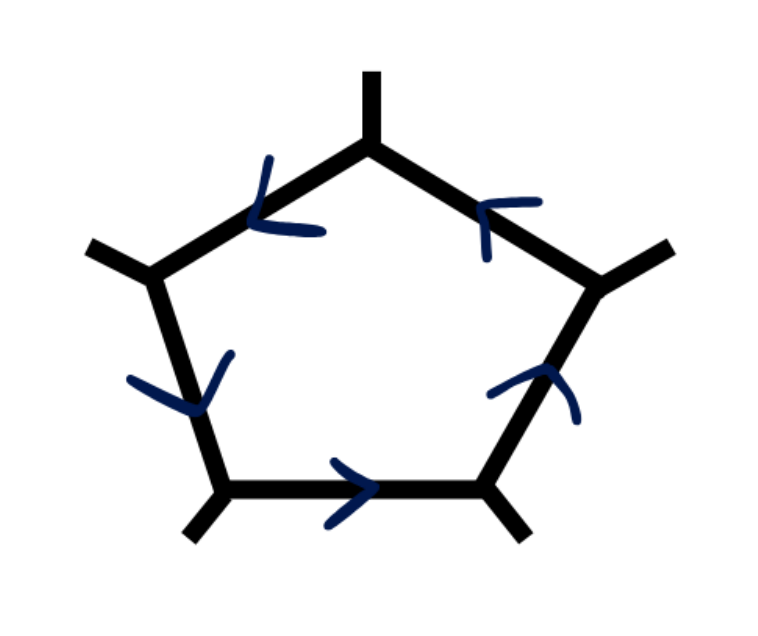}
    \caption{}
    \label{fig:oldFig.45(a)}
\end{figure}

(i) Assume that there exist two edges $a$ and $b$ extending from $l$ which are connected by a path. Then, as shown in Figure \ref{fig:oldFig.45(c)}, 
 there exists another closed path $l'$ which includes $a$ and $b$. As shown in Figure \ref{fig:oldFig.45(d)}, $l'$ must also be rotated by $f|_{\Gamma}$.
 This is a contradiction.

\begin{figure}[ht]
 \begin{minipage}[ht]{0.45\linewidth}
    \centering
    \includegraphics[width=0.75\linewidth]{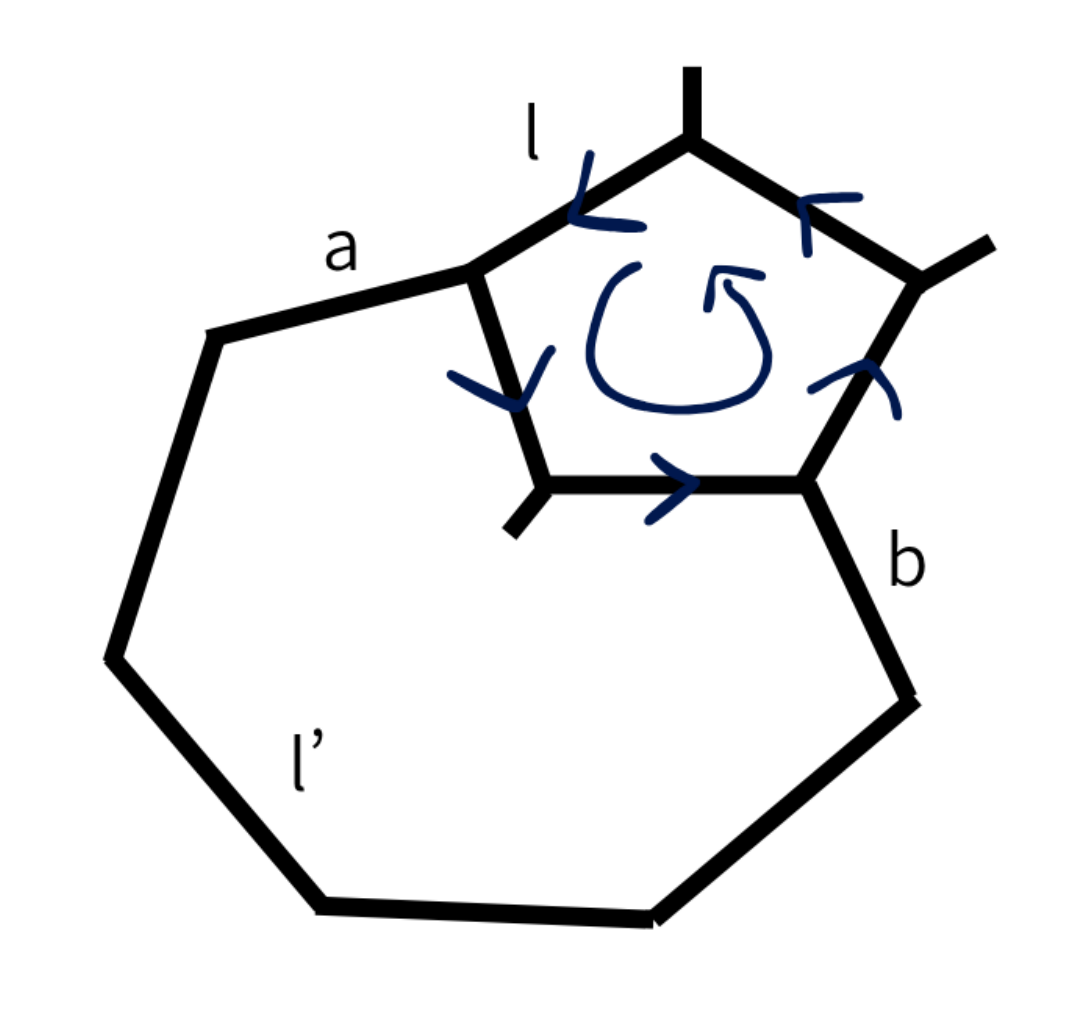}
    \caption{}
    \label{fig:oldFig.45(c)}
 \end{minipage}
 \begin{minipage}[ht]{0.45\linewidth}
    \centering
    \includegraphics[width=0.7\linewidth]{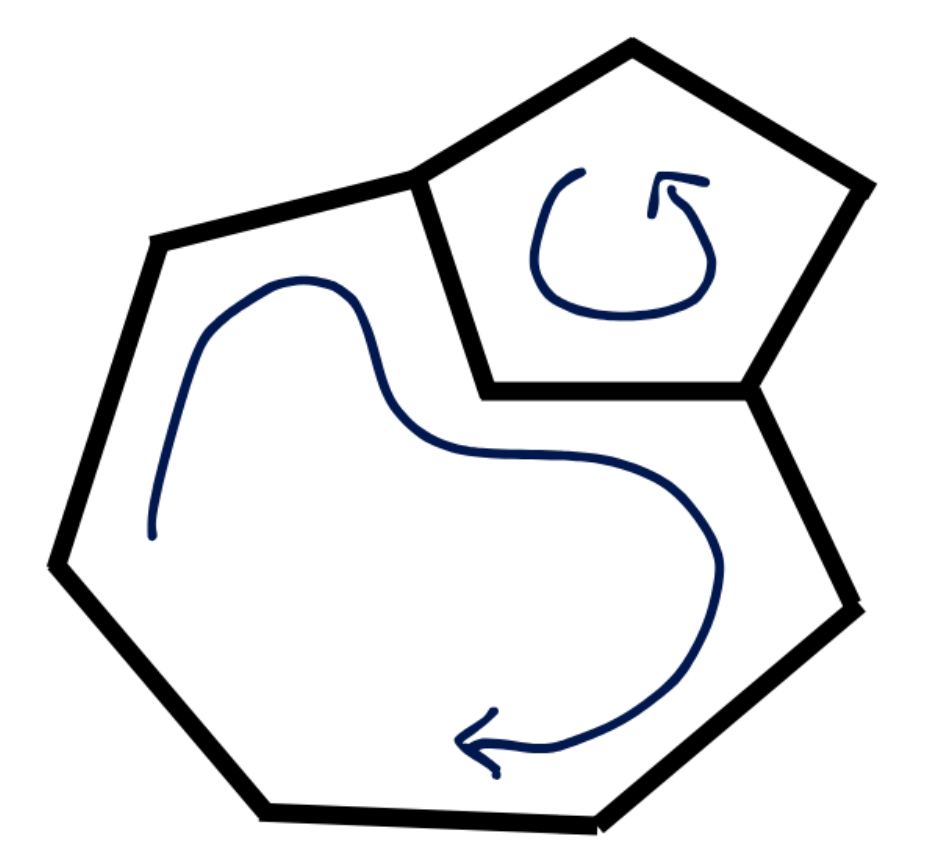}
    \caption{}
    \label{fig:oldFig.45(d)}
 \end{minipage}
\end{figure}

(ii) Assume that no two edges extending from $l$ are connected by a path, and a tail is connected to $l$ by a path. Let $t$ be the tail.
 See Figure \ref{fig:oldFig.45(b)}. $f|_{\Gamma}$ fixes $t$. This is a contradiction.
(Note that even if $r\geq2$, $f|_\Gamma$ fixes each tail.)

\begin{figure}[ht]
    \centering
    \includegraphics[width=0.4\linewidth]{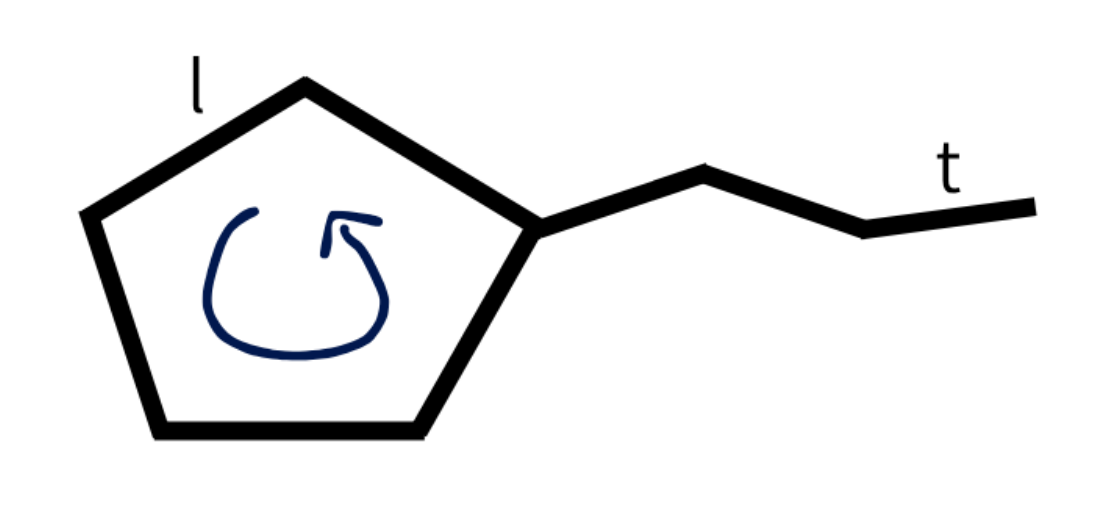}
    \caption{}
    \label{fig:oldFig.45(b)}
\end{figure}

(iii) Assume that no two edges extending from $l$ are connected by a path, and none of them is connected to a tail by a path.
 Then it becomes as shown in Figure \ref{fig:oldFig.45(e)}.
 Because $l$ is rotated by $f|_{\Gamma}$, a closed path $l_1$ is moved to another closed path $l_2$ by $f|_{\Gamma}$ as shown in Figure \ref{fig:oldFig.45(f)}.
 Then, $f|_{\Gamma}$ moves 1-homology of $\Gamma$. This is a contradiction.

\begin{figure}[ht]
 \begin{minipage}[ht]{0.45\linewidth}
    \centering
    \includegraphics[width=0.75\linewidth]{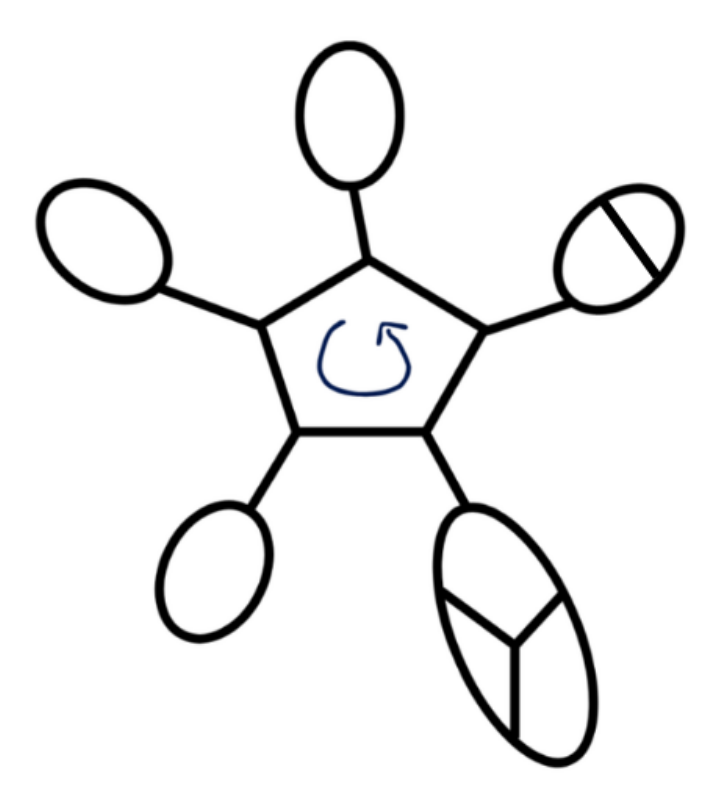}
    \caption{}
    \label{fig:oldFig.45(e)}
 \end{minipage}
 \begin{minipage}[ht]{0.45\linewidth}
    \centering
    \includegraphics[width=0.6\linewidth]{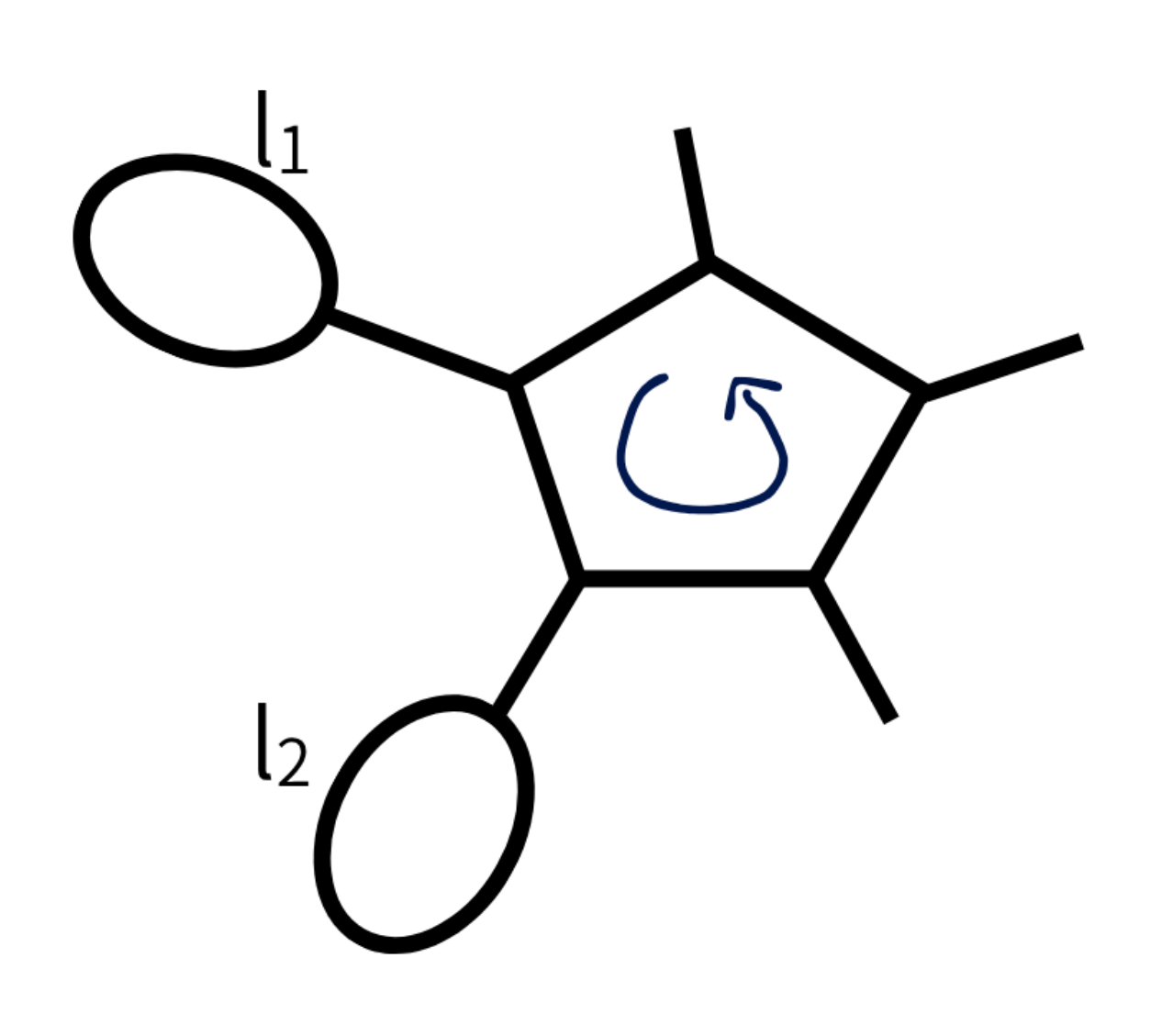}
    \caption{}
    \label{fig:oldFig.45(f)}
 \end{minipage}
\end{figure}
\begin{figure}[ht]
    \centering
    \includegraphics[width=0.4\linewidth]{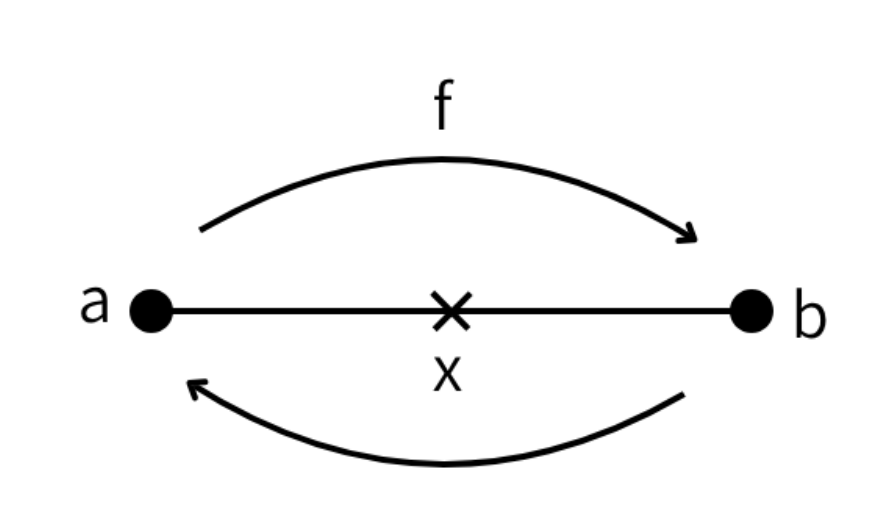}
    \caption{$f$ fixes $x$.}
    \label{fig:f(a)=bandf(b)=a}
\end{figure}
\begin{figure}[ht]
    \centering
    \includegraphics[width=0.5\linewidth]{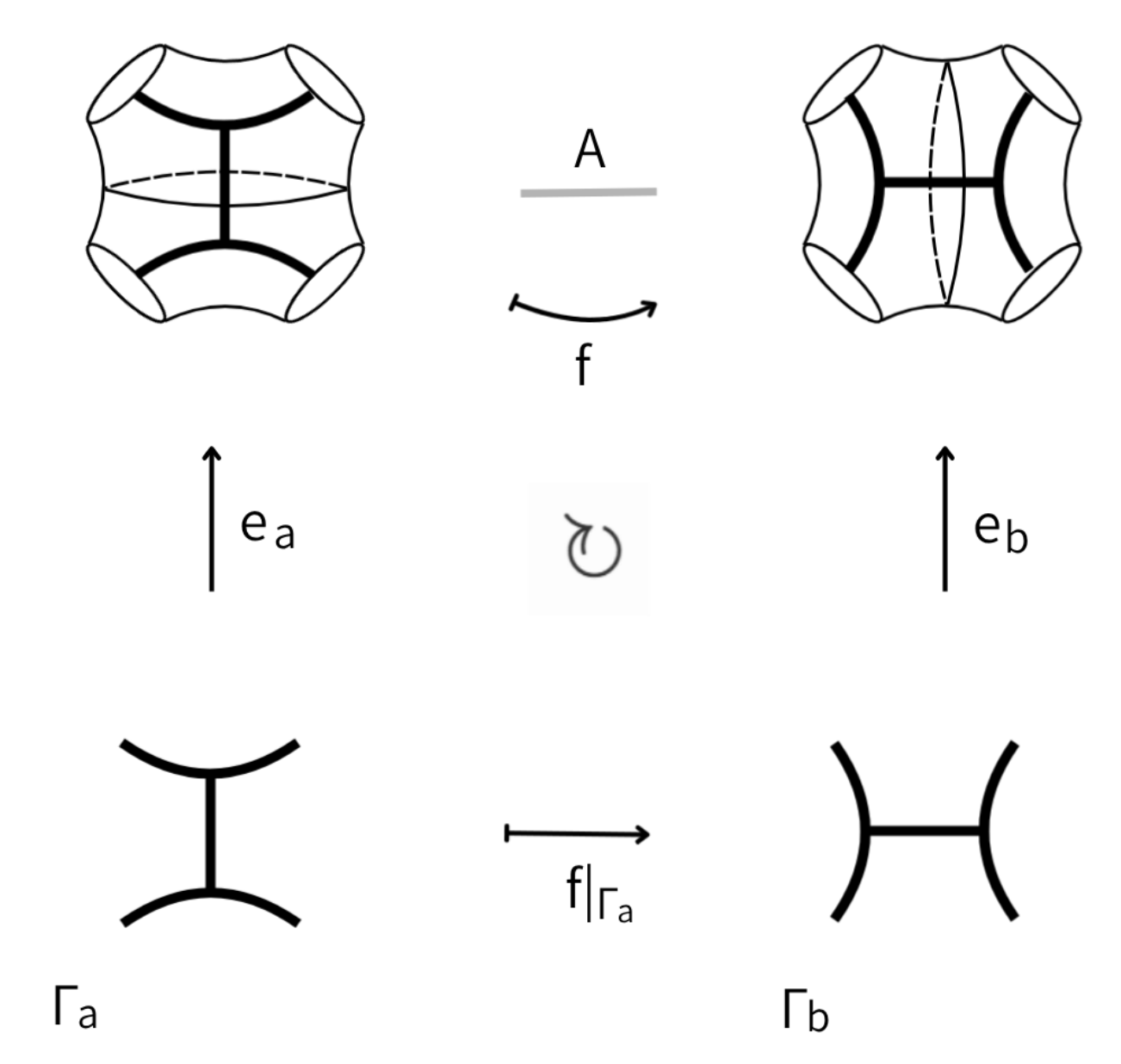}
    \caption{Define $f|_{\Gamma_a}$ by $\circlearrowright$.}
    \label{fig:defoffGammaa}
\end{figure}

 From the above discussion, $f|_{\Gamma}=1$. From claim \ref{thm:claim1}, $f=1$. 
 The proof of (1) of proposition \ref{thm:mainprop} is complete.
 
(2) Let $x$ be an inner point of 1-cell of $\overline{H}(\Sigma)$.
Assume $f \ne 1$. We want to show a contradiction.

From (1) of proposition \ref{thm:mainprop}, every 0-cell of $\overline{H}(\Sigma)$ is moved by $f$.
Consider that $a$ and $b$ are 0-cells at both ends of the 1-cell which includes $x$ as its inner point.
As shown in Figure \ref{fig:f(a)=bandf(b)=a}, $f(a)=b$ and $f(b)=a$.

The 1-cell has four possibilities: $S,A,T$ or $B'$.

Each $S$, $T$ and $B'$ is oriented as 1-cell. For example, $T^{-1}$ is not $T$. See Figure \ref{fig:T}.
$f$ preserves those orientation. In these cases, $f$ as shown in Figure \ref{fig:f(a)=bandf(b)=a} does not exist.

\clearpage

We only need to consider the case where the 1-cell is $A$-move.

\begin{lemma}\label{thm:claim3}
    If $a$ and $b$ are 0-cells at both ends of $A$-move and $f(a)=b$, then $f^k\ne 1$ for $k=2,4,5$.
\end{lemma}

\begin{proof}

$a$ and $b$ are colored quilt decompositions. Let $\Gamma_a$ and $\Gamma_b$ be trivalent graphs with $r$ tails corresponding to $a$ and $b$.
Let $e_a$ and $e_b$ be embeddings from $\Gamma_a$ and $\Gamma_b$ to $\Sigma$. 
Then, we define that $f|_{\Gamma_a}:=e_b^{-1}\circ f\circ e_a$ (see Figure \ref{fig:defoffGammaa}).
Remark that $f|_{\Gamma_a}:\Gamma_a\cong\Gamma_b$ as a trivalent graph.

Now, we consider the following maps of 1-homology induced by the above maps.

\[\xymatrix{
H_1(\Sigma) \ar[rr]^{f^*=\text{id}_{H_1(\Sigma)}} && H_1(\Sigma) \\
& \circlearrowright & \\
H_1(\Gamma_a) \ar[rr]^{(f|_{\Gamma_a})^*} \ar[uu]^{e_a^*} && H_1(\Gamma_b) \ar[uu]^{e_b^*} \\
}\]

$f^*=\text{id}_{H_1(\Sigma)}$ because $f\in\mathcal{I}_{g,r}$. 
$e_a^*$ and $e_b^*$ are injective from claim \ref{thm:injectiveclaim}.
From the definition of $f|_{\Gamma_a}$, the following holds:

\begin{equation}
    (f|_{\Gamma_a})^*=(f^*)|_{H_1(\Gamma_a)}          \tag{3.1}\label{*}
\end{equation}

\begin{figure}[ht]

 \begin{minipage}[ht]{0.3\linewidth}
    \centering
    \includegraphics[width=0.5\linewidth]{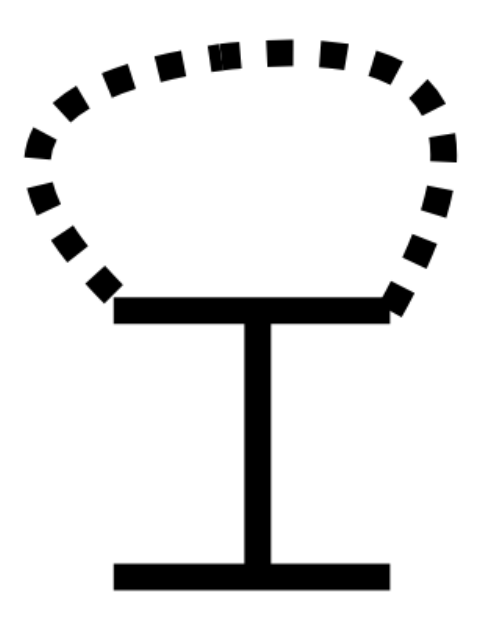}
    \caption{}
    \label{fig:Fig48a}
 \end{minipage}
 \begin{minipage}[ht]{0.3\linewidth}
    \centering
    \includegraphics[width=0.5\linewidth]{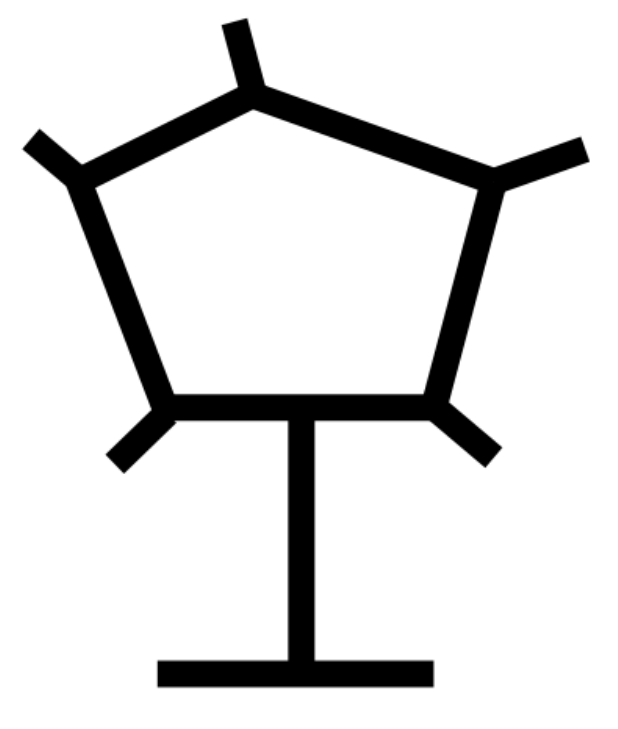}
    \caption{}
    \label{fig:Fig48b}
 \end{minipage}
 \begin{minipage}[ht]{0.3\linewidth}
    \centering
    \includegraphics[width=0.5\linewidth]{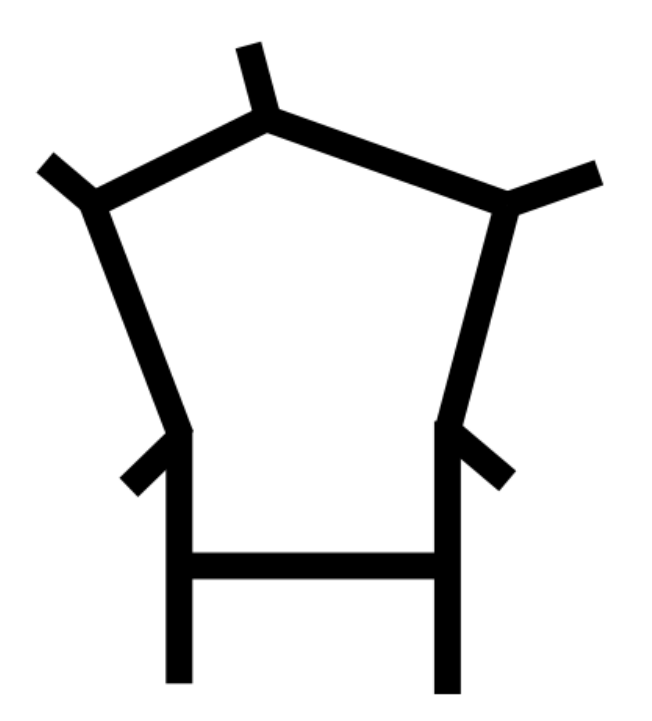}
    \caption{}
    \label{fig:Fig48c}
 \end{minipage}
\end{figure}

(I) Assume that $\Gamma_a$ is connected as shown in Figure \ref{fig:Fig48a}.
In this figure, the ``letter I" represents the support of the move, 
and the dotted line means that the two edges are connected by a path.
That is, there exists a closed path $p$ which is $n$-gon ($n\geq 1$) in $\Gamma_a$. See Figure \ref{fig:Fig48b}.
This figure shows an example where $n=6$.
Then, $\Gamma_b$ has a closed path $p'$ which is $(n+1)$-gon. See Figure \ref{fig:Fig48c}.
From \eqref{*}, $f|_{\Gamma_a}(p)=p'$. Such $f$ does not exist.

(II) Assume that $\Gamma_a$ is connected as shown in Figure \ref{fig:Fig49a}.
That is, there exists a closed path $p$ which is $n$-gon ($n\geq 2$) in $\Gamma_a$. See Figure \ref{fig:Fig49b}.
Then, $\Gamma_b$ has a closed path $p'$ which is $(n-1)$-gon. See Figure \ref{fig:Fig49c}.
From \eqref{*}, $f|_{\Gamma_a}(p)=p'$. Such $f$ does not exist.

\begin{figure}[ht]
 \begin{minipage}[ht]{0.3\linewidth}
    \centering
    \includegraphics[width=0.5\linewidth]{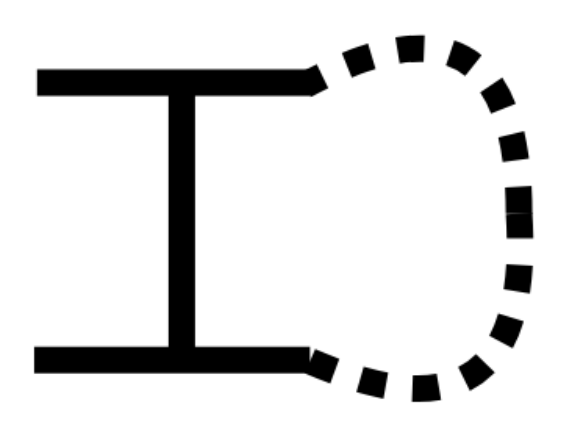}
    \caption{}
    \label{fig:Fig49a}
 \end{minipage}
 \begin{minipage}[ht]{0.3\linewidth}
    \centering
    \includegraphics[width=0.5\linewidth]{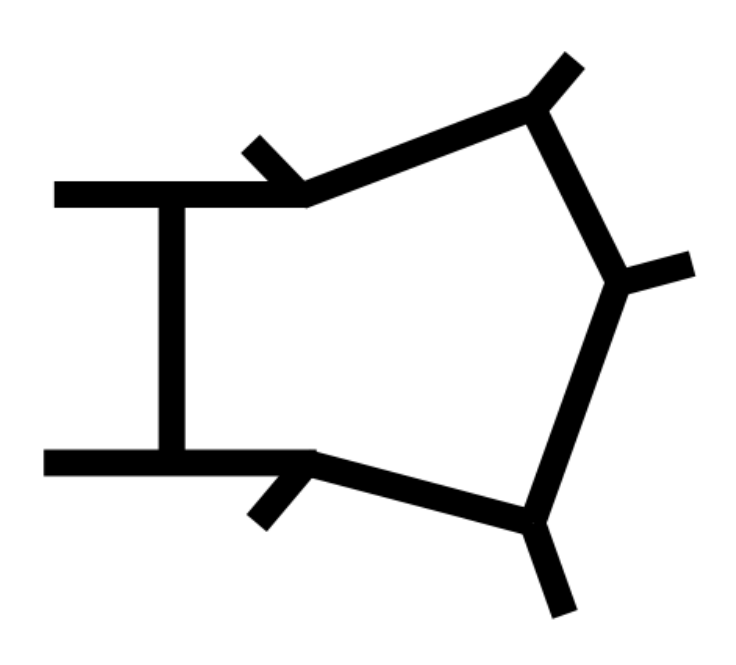}
    \caption{}
    \label{fig:Fig49b}
 \end{minipage}
 \begin{minipage}[ht]{0.3\linewidth}
    \centering
    \includegraphics[width=0.5\linewidth]{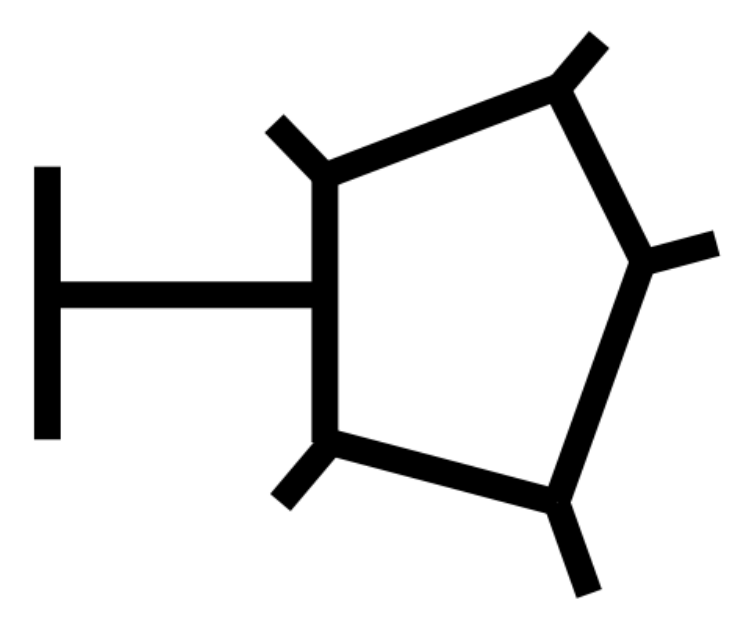}
    \caption{}
    \label{fig:Fig49c}
 \end{minipage}
\end{figure}

(III) Assume that $\Gamma_a$ is connected as shown in Figure \ref{fig:Fig50a}.
That is, there exists a closed path $p$ which is $n$-gon ($n\geq 2$) in $\Gamma_a$. See Figure \ref{fig:Fig50b}.
Then, $\Gamma_b$ has a closed path $p'$ which is $n$-gon. See Figure \ref{fig:Fig50c}.
From \eqref{*}, $f|_{\Gamma_a}(p)=p'$.

\begin{figure}[ht]
 \begin{minipage}[ht]{0.3\linewidth}
    \centering
    \includegraphics[width=0.5\linewidth]{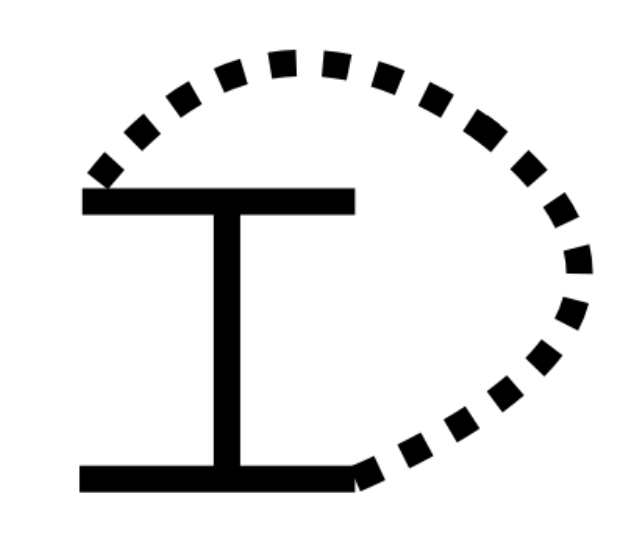}
    \caption{}
    \label{fig:Fig50a}
 \end{minipage}
 \begin{minipage}[ht]{0.3\linewidth}
    \centering
    \includegraphics[width=0.5\linewidth]{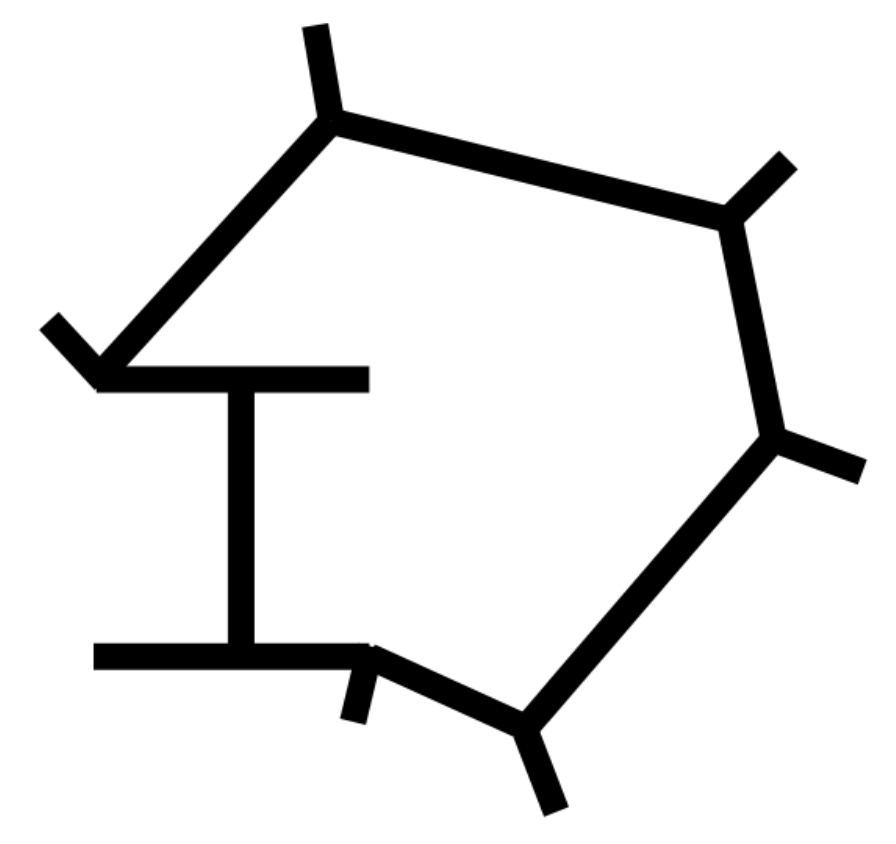}
    \caption{}
    \label{fig:Fig50b}
 \end{minipage}
 \begin{minipage}[ht]{0.3\linewidth}
    \centering
    \includegraphics[width=0.5\linewidth]{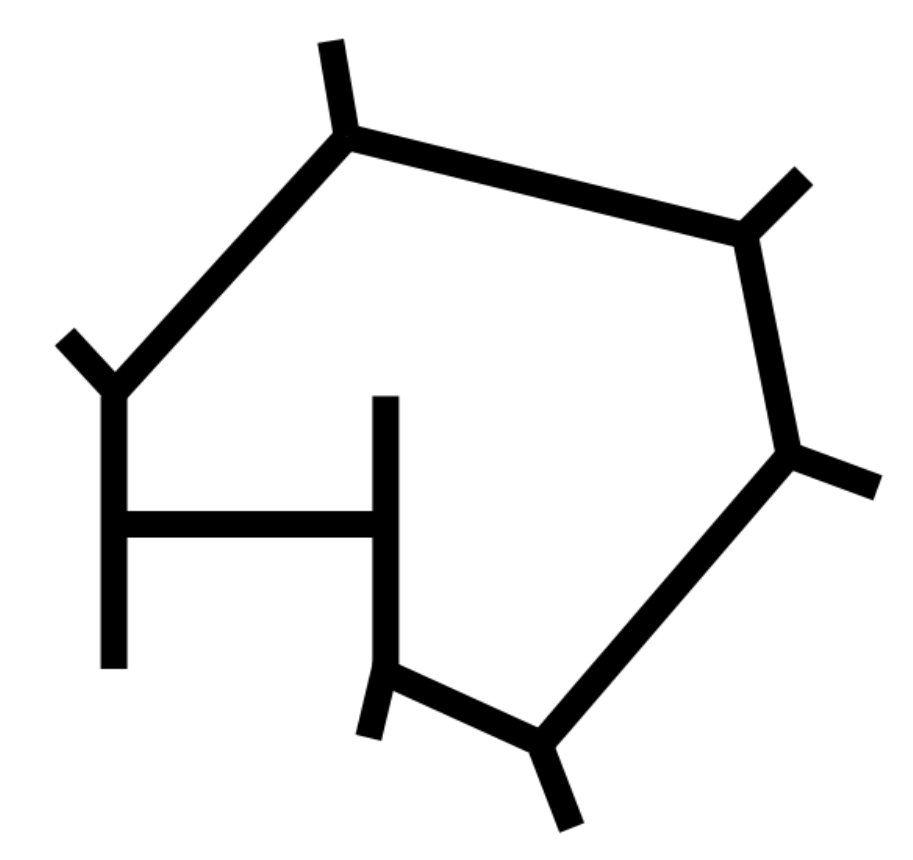}
    \caption{}
    \label{fig:Fig50c}
 \end{minipage}
\end{figure}

(III)-(i) Assume that $\Gamma_a$ is further connected as shown in Figure \ref{fig:Fig51a}.

\begin{figure}[ht]
 \begin{minipage}[ht]{0.3\linewidth}
    \centering
    \includegraphics[width=0.6\linewidth]{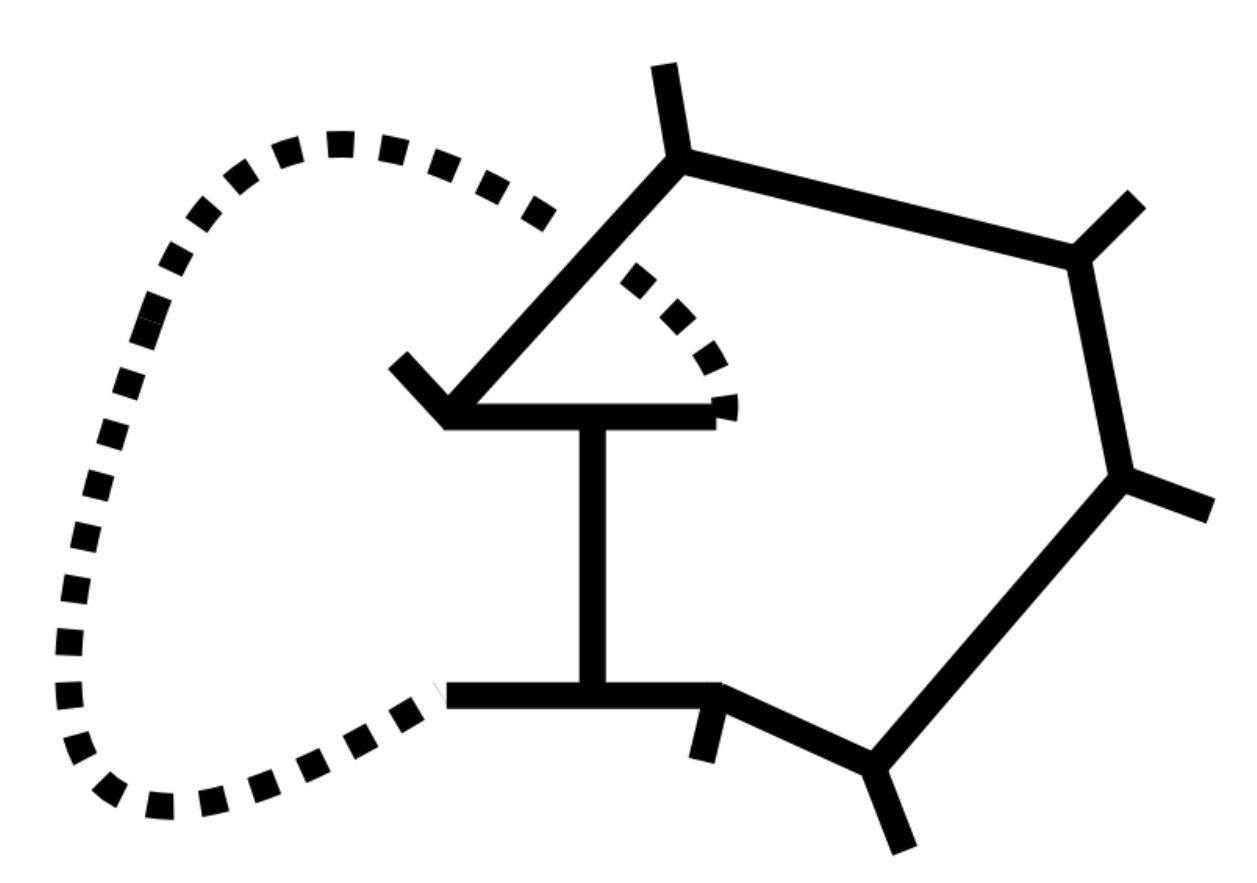}
    \caption{}
    \label{fig:Fig51a}
 \end{minipage}
 \begin{minipage}[ht]{0.3\linewidth}
    \centering
    \includegraphics[width=0.8\linewidth]{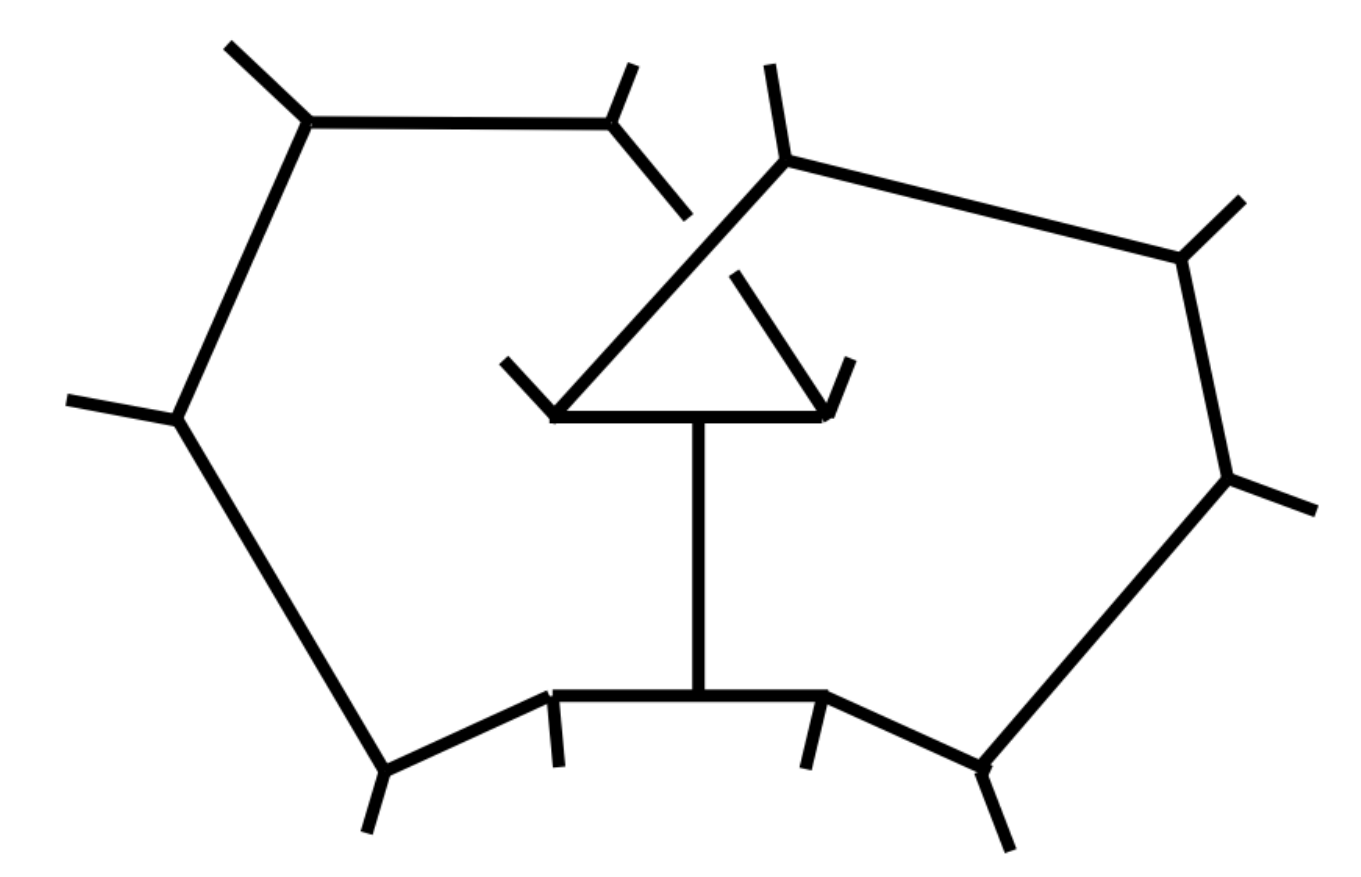}
    \caption{}
    \label{fig:Fig51b}
 \end{minipage}
 \begin{minipage}[ht]{0.3\linewidth}
    \centering
    \includegraphics[width=0.8\linewidth]{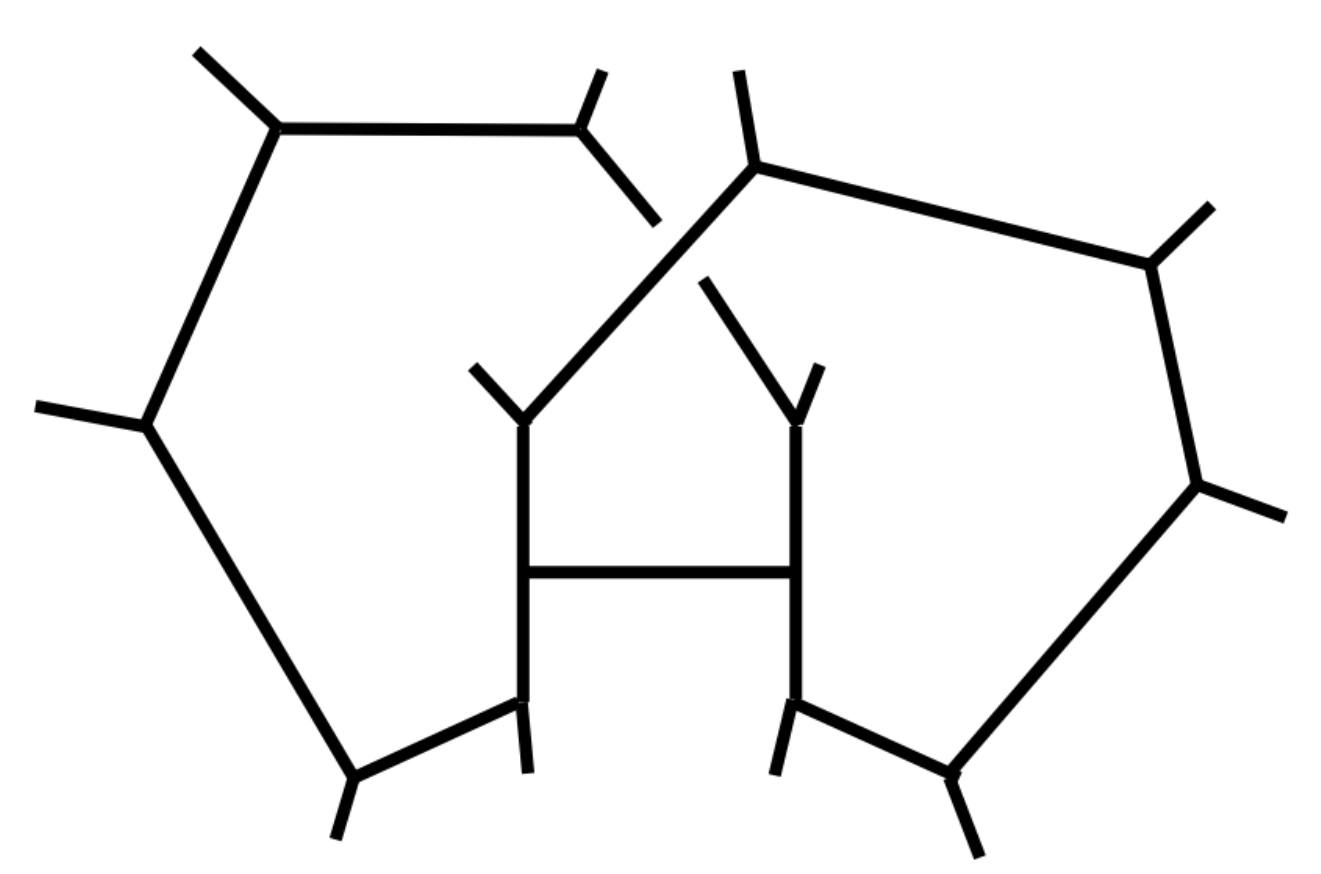}
    \caption{}
    \label{fig:Fig51c}
 \end{minipage}
\end{figure}

That is, there exists a closed path $q$ which is $m$-gon ($m\geq 2$) in $\Gamma_a$. See Figure \ref{fig:Fig51b}.
Then, $\Gamma_b$ has closed paths $q'$ which is $m$-gon. See Figure \ref{fig:Fig51c}. From \eqref{*}, $f|_{\Gamma_a}(q)=q'$. 
Suppose 1-cycle $\widetilde{p}$ and 1-cycle $\widetilde{q}$ on $\Gamma_a$ as shown in Figure \ref{fig:Fig51d}. And,
 suppose 1-cycle $\widetilde{p'}$ and 1-cycle $\widetilde{q'}$ on $\Gamma_b$ as shown in Figure \ref{fig:Fig51e}.
Since $f|_{\Gamma_a}(p)=p'$ and $f|_{\Gamma_a}(q)=q'$, 
$(f|_{\Gamma_a})^*(\widetilde{p})=\widetilde{p'}$ and $(f|_{\Gamma_a})^*(\widetilde{q})=\widetilde{q'}$.
However, as can be seen by comparing Figure \ref{fig:Fig51d} and Figure \ref{fig:Fig51e},
the orientation of their 1-cycles is inverse at the overlapping edge.
Such $f$ does not exist.


(III)-(ii) Assume that $\Gamma_a$ is disconnected as shown in Figure \ref{fig:Fig52a}.
On the figure, the two dotted lines mean they are disconnected.
From (I) and (II), those two edges are not connected to anything else.
So, the two edges have tails or cycles, i.e. the two edges correspond to BSCCs (bounding simple closed curves). 
 See Figure \ref{fig:Fig52b}. On the figure, $B$ means an edge corresponding to BSCC.
$\Gamma_b$ is as shown in Figure \ref{fig:Fig52c}.
Suppose 1-cycle $\widetilde{p}$ on $\Gamma_a$ as shown in Figure \ref{fig:Fig52d},
and suppose 1-cycle $\widetilde{p'}$ on $\Gamma_b$ as shown in Figure \ref{fig:Fig52e}.
Since $f|_{\Gamma_a}(p)=p'$, $(f|_{\Gamma_a})^*(\widetilde{p})=\widetilde{p'}$.

\begin{figure}[H]
 \begin{minipage}[H]{0.45\linewidth}
    \centering
    \includegraphics[width=0.8\linewidth]{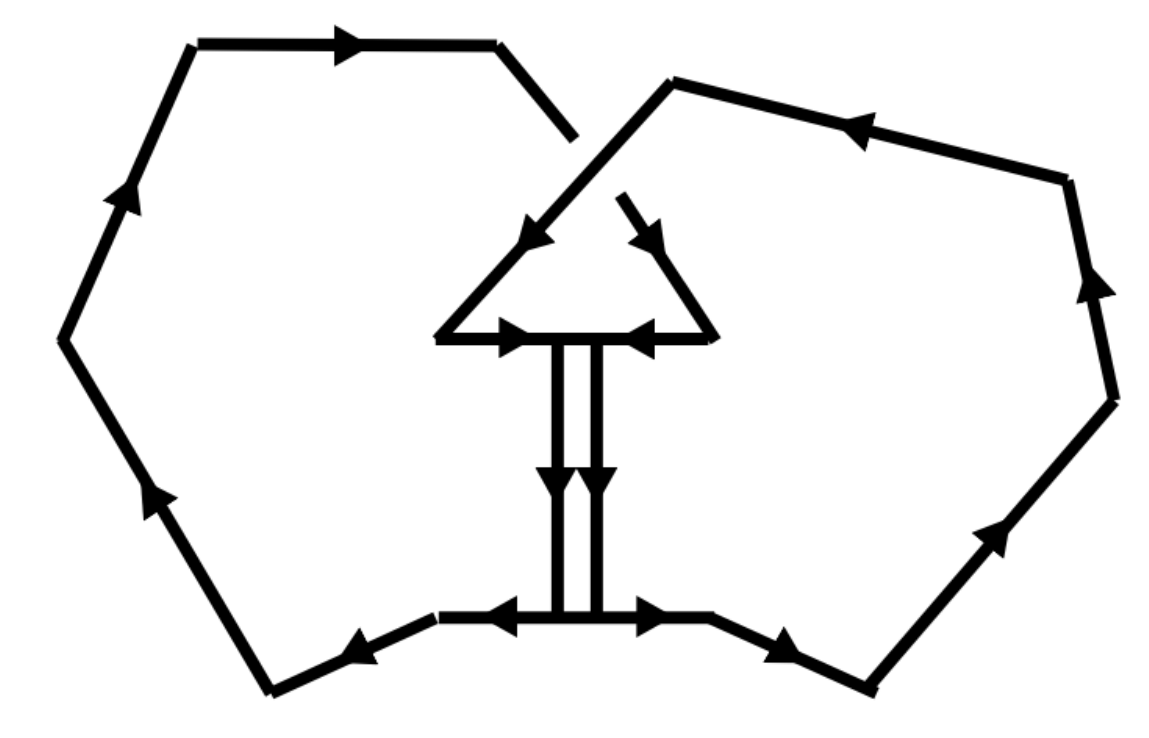}
    \caption{}
    \label{fig:Fig51d}
 \end{minipage}
 \begin{minipage}[H]{0.45\linewidth}
    \centering
    \includegraphics[width=0.8\linewidth]{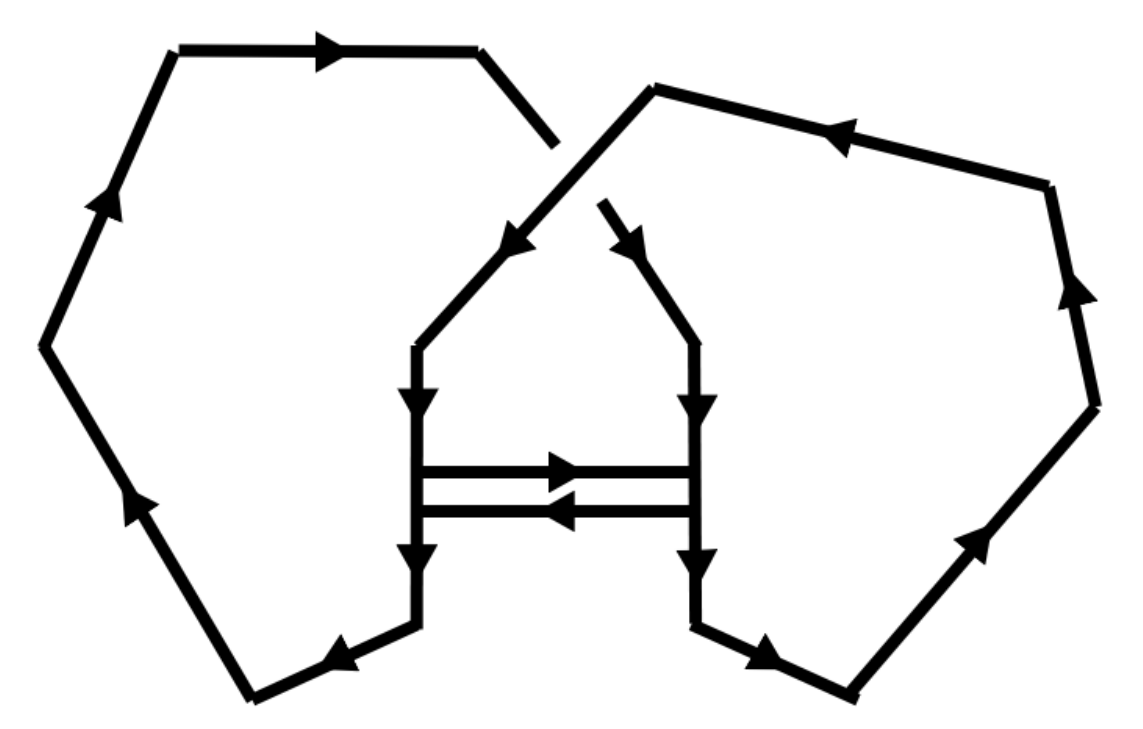}
    \caption{}
    \label{fig:Fig51e}
 \end{minipage}
\end{figure}
\begin{figure}[H]
 \begin{minipage}[H]{0.3\linewidth}
    \centering
    \includegraphics[width=0.9\linewidth]{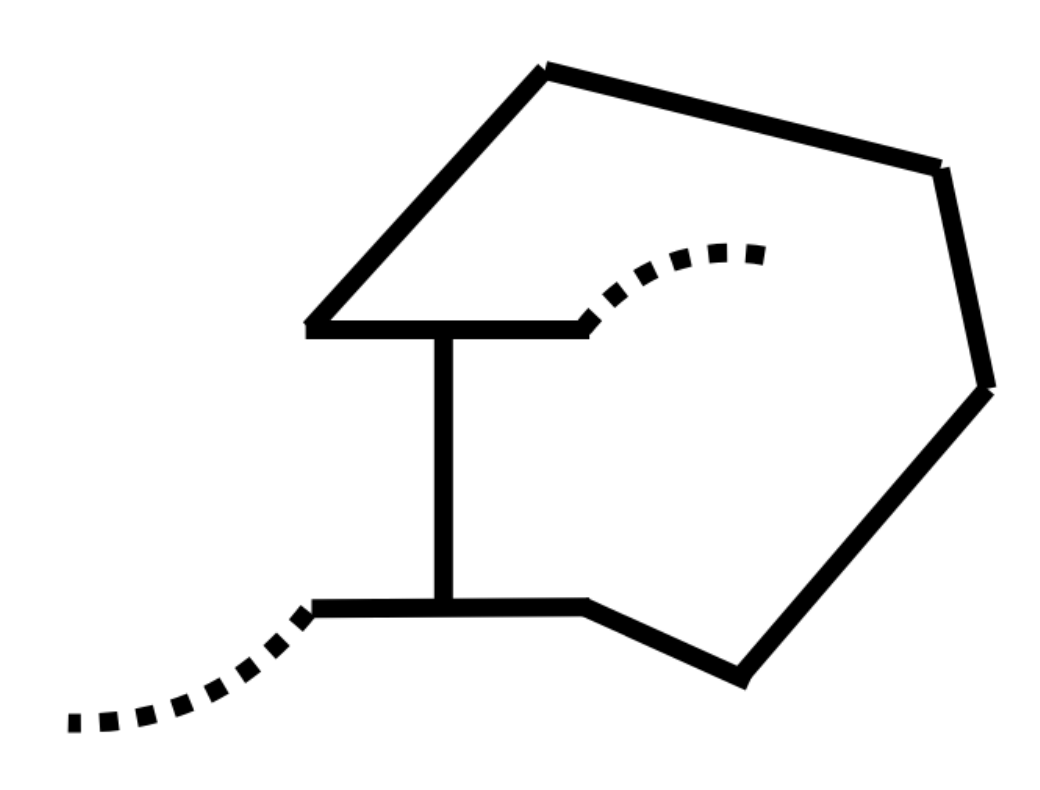}
    \caption{}
    \label{fig:Fig52a}
 \end{minipage}
 \begin{minipage}[H]{0.3\linewidth}
    \centering
    \includegraphics[width=0.8\linewidth]{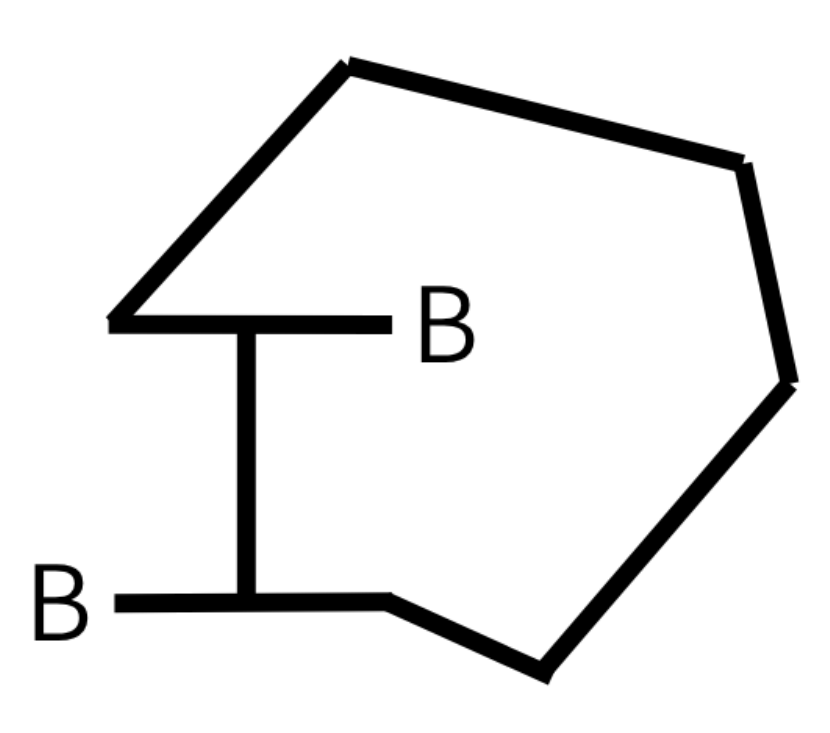}
    \caption{}
    \label{fig:Fig52b}
 \end{minipage}
 \begin{minipage}[H]{0.3\linewidth}
    \centering
    \includegraphics[width=0.8\linewidth]{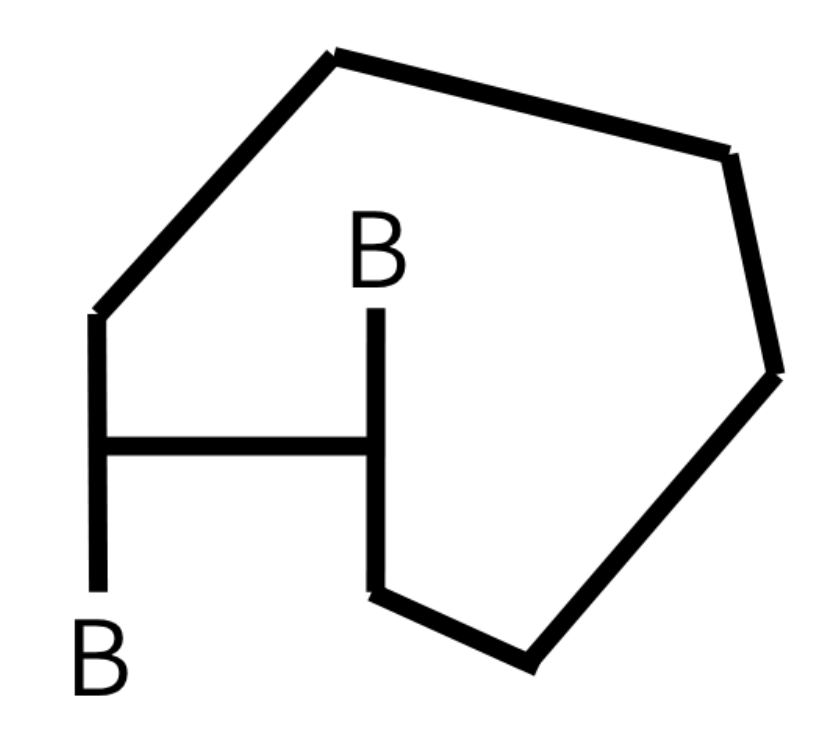}
    \caption{}
    \label{fig:Fig52c}
 \end{minipage}
\end{figure}
\begin{figure}[H]
 \begin{minipage}[H]{0.45\linewidth}
    \centering
    \includegraphics[width=0.5\linewidth]{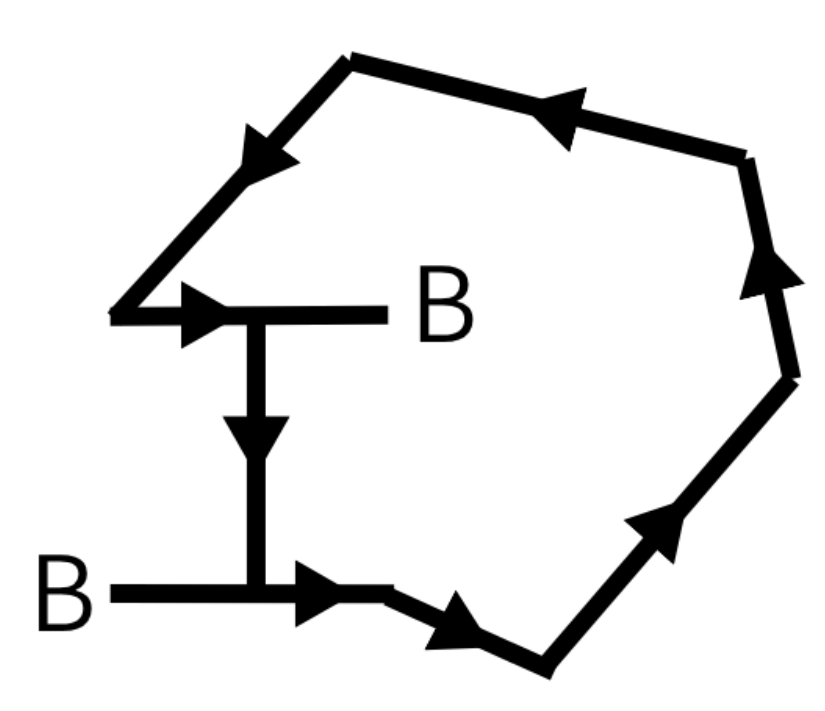}
    \caption{}
    \label{fig:Fig52d}
 \end{minipage}
 \begin{minipage}[H]{0.45\linewidth}
    \centering
    \includegraphics[width=0.5\linewidth]{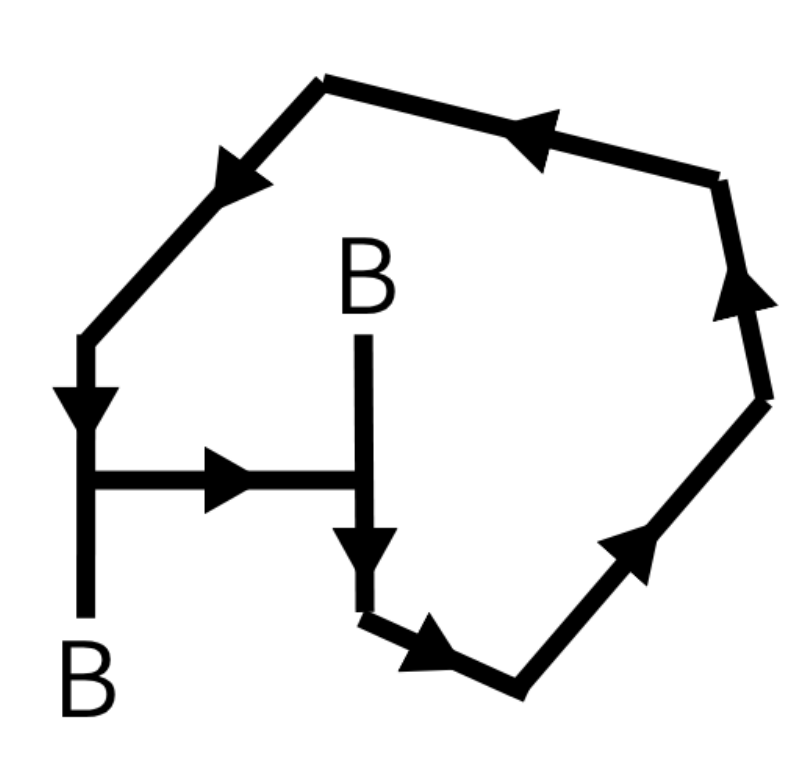}
    \caption{}
    \label{fig:Fig52e}
 \end{minipage}
\end{figure}
\begin{figure}[H]
 \begin{minipage}[H]{0.45\linewidth}
    \centering
    \includegraphics[width=0.5\linewidth]{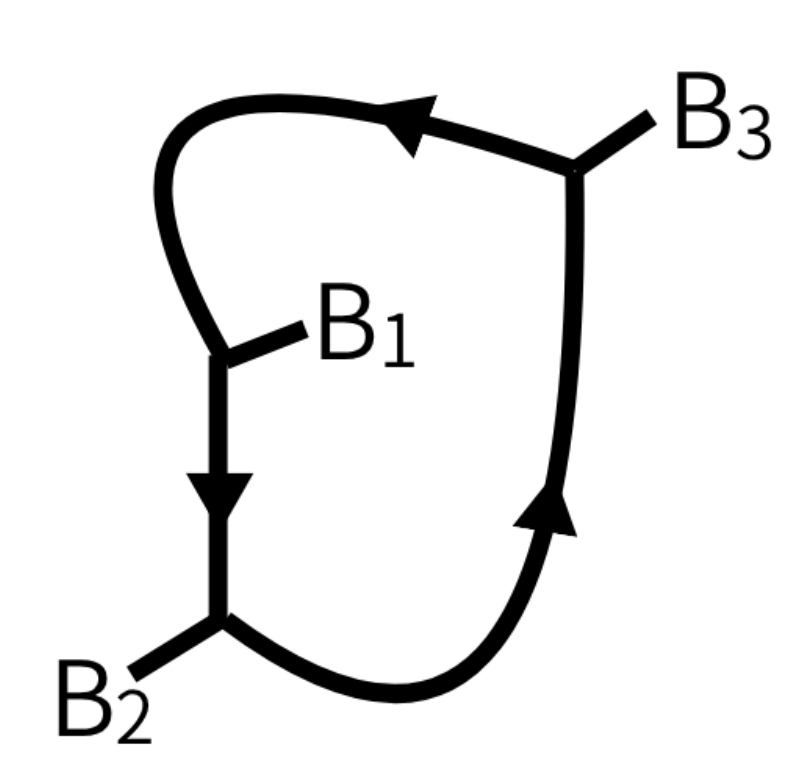}
    \caption{}
    \label{fig:Fig53a}
 \end{minipage}
 \begin{minipage}[H]{0.45\linewidth}
    \centering
    \includegraphics[width=0.5\linewidth]{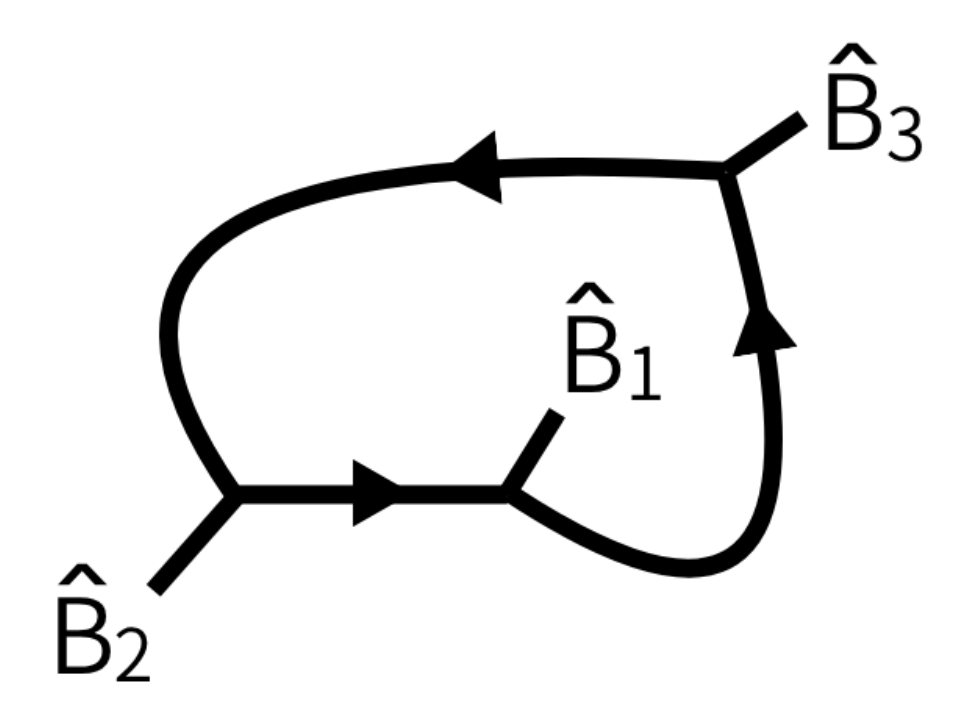}
    \caption{}
    \label{fig:Fig53b}
 \end{minipage}
\end{figure}
\begin{figure}[H]
 \begin{minipage}[H]{0.45\linewidth}
    \centering
    \includegraphics[width=0.65\linewidth]{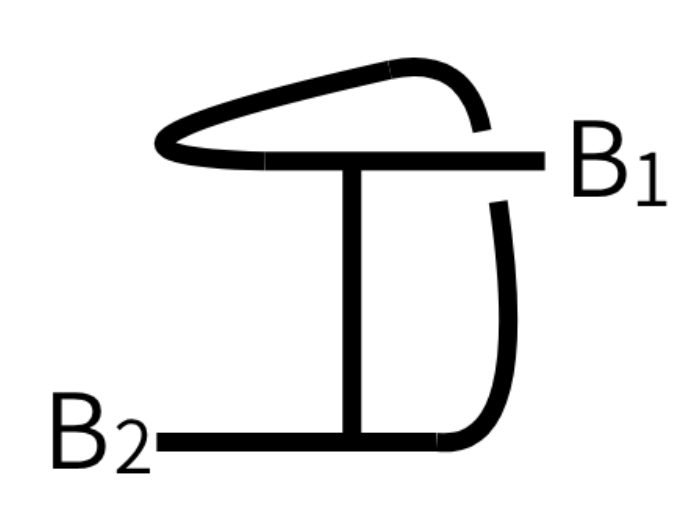}
    \caption{}
    \label{fig:Fig54a}
 \end{minipage}
 \begin{minipage}[H]{0.45\linewidth}
    \centering
    \includegraphics[width=0.5\linewidth]{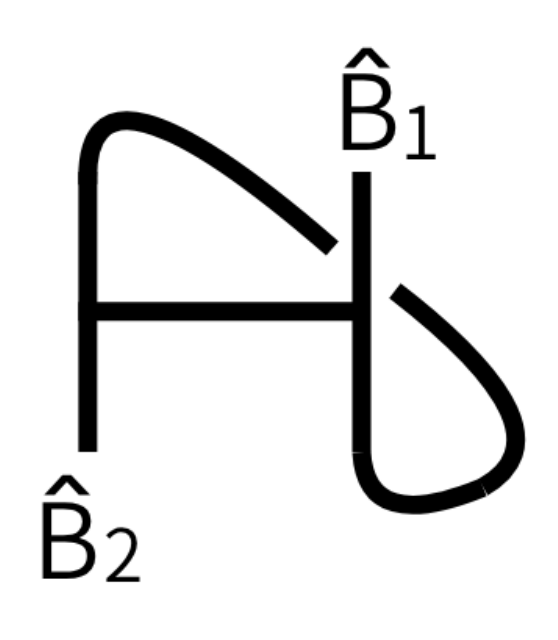}
    \caption{}
    \label{fig:Fig54b}
 \end{minipage}
\end{figure}

(III)-(ii)-(1) Assume $n\geq 3$.
Let $B_1$ and $B_2$ be those two edges corresponding BSCCs on $\Gamma_a$.
There exists at least one edge corresponding BSCC extending from $p$ other than $B_1$ and $B_2$ since $n\geq 3$ and from (I) and (II).
Let $B_3$ be one of the edges.
See Figure \ref{fig:Fig53a}.
Set $\widehat{B}_k:=f|_{{\Gamma}_a}(B_k)$ for $k=1,2,3$. Each $\widehat{B}_k$ is an edge corresponding to BSCC on $\Gamma_b$. 
Then, each $\widehat{B}_k$ is attached to the closed path $p'$
on $\Gamma_b$. See Figure \ref{fig:Fig53b}.
In this case, the cyclic order of $B_1$, $B_2$, $B_3$ induced by the orientation of $\widetilde{p}$ is inverse
to the cyclic order of $\widehat{B}_1$, $\widehat{B}_2$, $\widehat{B}_3$
induced by the orientation of $\widetilde{p'}$.
Since $f|_{{\Gamma}_a}$ is an isomorphism of the graph,
such $f$ does not exist.

(III)-(ii)-(2) Assume $n=2$. Then, $p$ is as shown in Figure \ref{fig:Fig54a}, and $p'$ is as shown in Figure \ref{fig:Fig54b}.
These figures depict the following. The ``letter I" located toward the front in Figure \ref{fig:Fig54a} 
and the ``letter H" located toward the front in Figure \ref{fig:Fig54b} correspond to 
the neighborhood of the support of the move, 
while the curved lines behind them each represent the same single edge in both figures. It is a single edge because $n=2$.
On Figure \ref{fig:Fig54a}, $B_1$ and $B_2$ are edges corresponding to BSCCs in $\Gamma_a$.
On Figure \ref{fig:Fig54b}, $\widehat{B}_1$ and $\widehat{B}_2$ are edges corresponding to BSCCs in $\Gamma_b$, where 
 $f|_{{\Gamma}_a}(B_k)=\widehat{B}_k$ for $k=1,2$.

In this case, the only possibility of our situation is as Figure \ref{fig:Fig54c}.
The upper part of Figure \ref{fig:Fig54c} is $a$, and the lower part of Figure \ref{fig:Fig54c} is $b$.

\begin{figure}[p]
    \centering
    \includegraphics[width=0.6\linewidth]{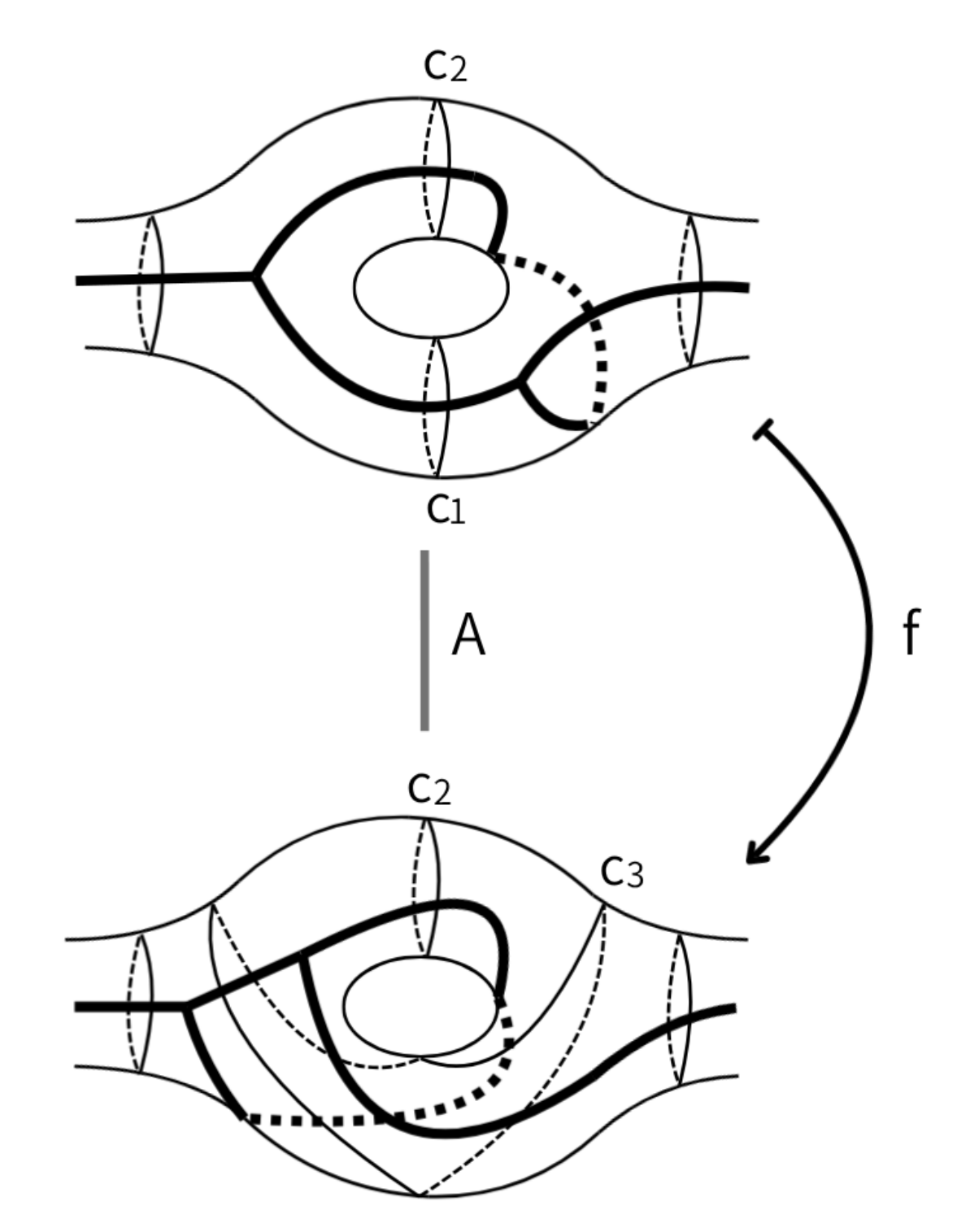}
    \caption{}
    \label{fig:Fig54c}
\end{figure}

\begin{figure}[p]
    \centering
    \includegraphics[width=0.4\linewidth]{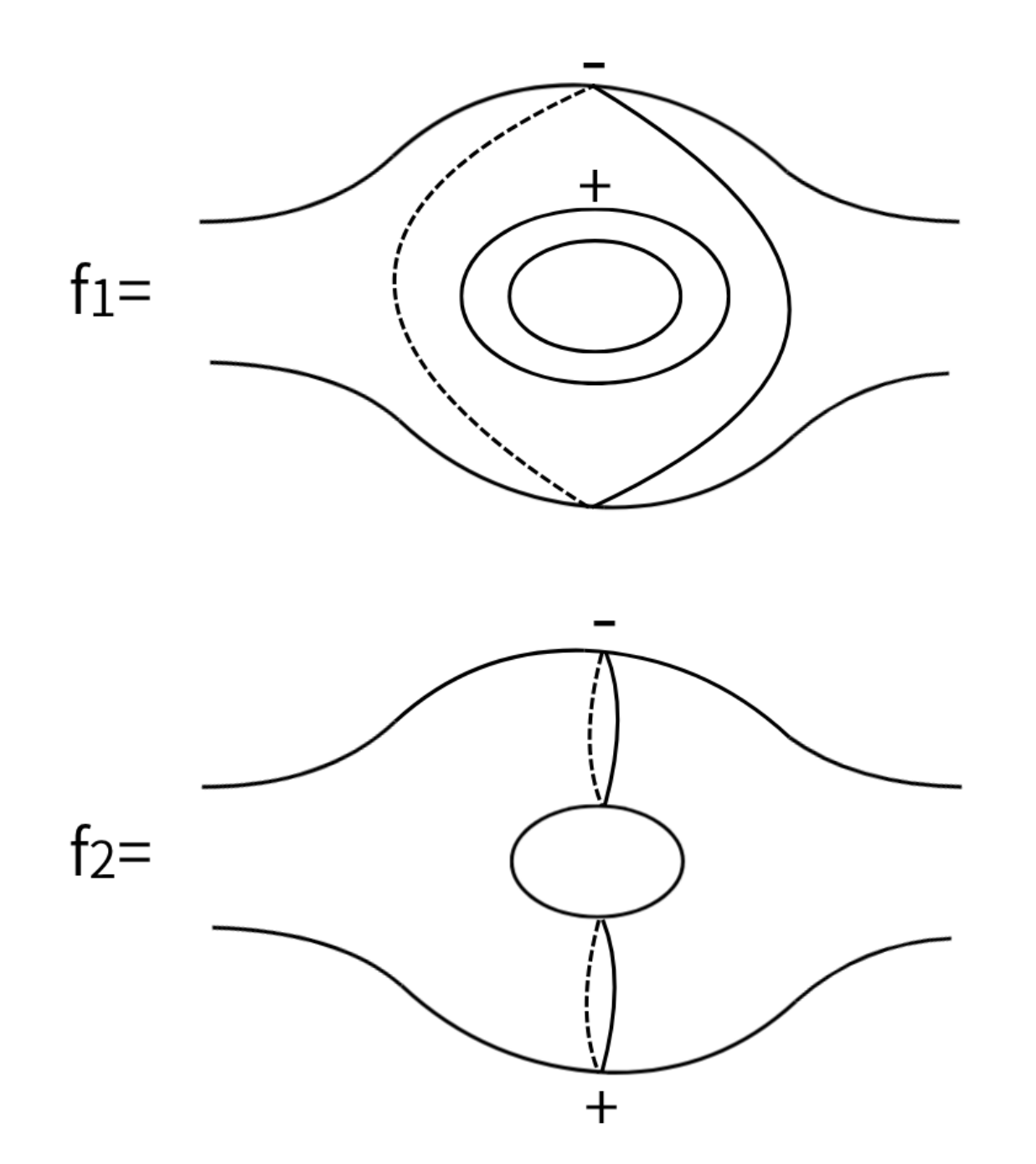}
    \caption{}
    \label{fig:Fig54d}
\end{figure}

Let $f_1,f_2\in \mathcal{M}_{g,r}$ be as shown in Figure \ref{fig:Fig54d}, where ``$+$" means a positive Dehn twist along the SCC and 
 ``$-$" means a negative Dehn twist along the SCC.
Then, one can see that $f=f_2\circ f_1$ and $f(a)=b$. 
 (Note that $f_1,f_2\in\mathcal{I}_{g,r}$ and hence $f\in\mathcal{I}_{g,r}$.)

As shown in Figure \ref{fig:Fig54c}, $f(c_1)=c_2$ and $f(c_2)=c_3$, so, $f^2(c_1)=c_3$.
Therefore, $f^2\ne 1$.
One can also check that $f^4\ne 1$ and $f^5\ne 1$.



(IV) Assume that the four edges of $\Gamma_a$ is all disconnected as shown in Figure \ref{fig:Fig55a}.
Then, all four of those edges correspond to BSCCs by the above discussion. 
Let $B_1,B_2,B_3,B_4$ be those four edges. See Figure \ref{fig:Fig55b}.
From \eqref{*}, $\Gamma_b$ becomes as shown in Figure \ref{fig:Fig55c} by $f|_{\Gamma_a}$,
where $f|_{\Gamma_a}(B_k)=\widehat{B}_k$ for $k=1,2,3,4$.
Such $f$ does not exist. 

\clearpage

\begin{figure}[ht]
 \begin{minipage}[ht]{0.3\linewidth}
    \centering
    \includegraphics[width=0.9\linewidth]{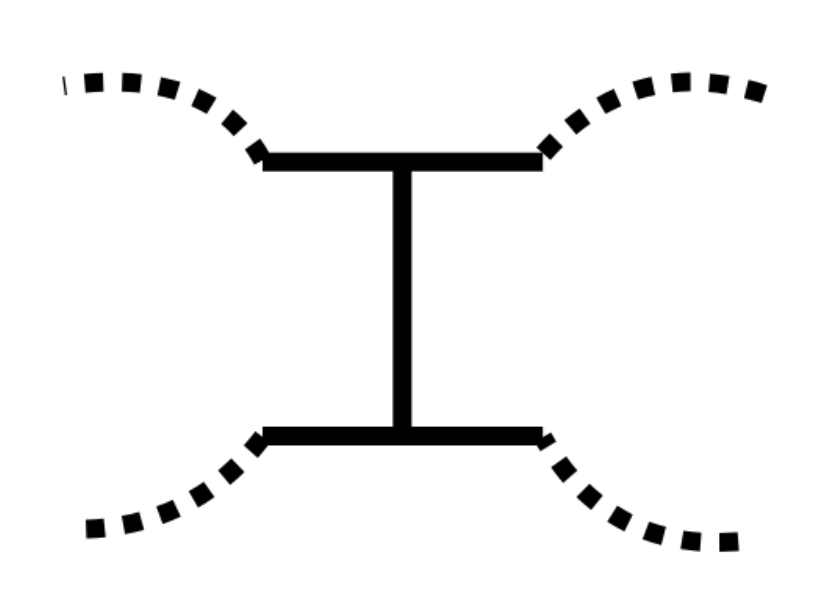}
    \caption{}
    \label{fig:Fig55a}
 \end{minipage}
 \begin{minipage}[ht]{0.3\linewidth}
    \centering
    \includegraphics[width=0.8\linewidth]{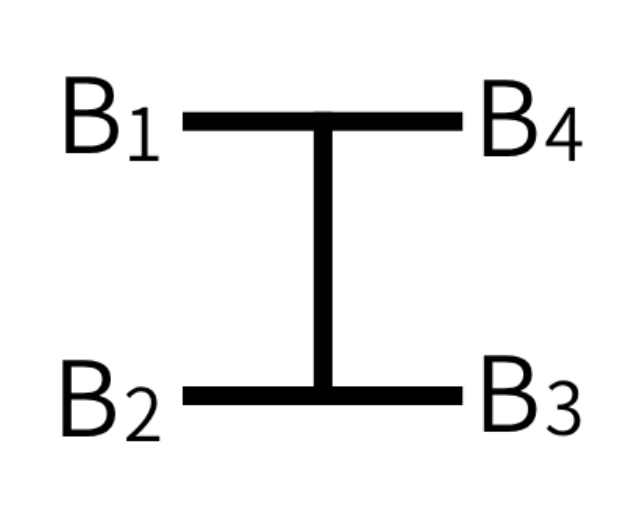}
    \caption{}
    \label{fig:Fig55b}
 \end{minipage}
 \begin{minipage}[ht]{0.3\linewidth}
    \centering
    \includegraphics[width=0.8\linewidth]{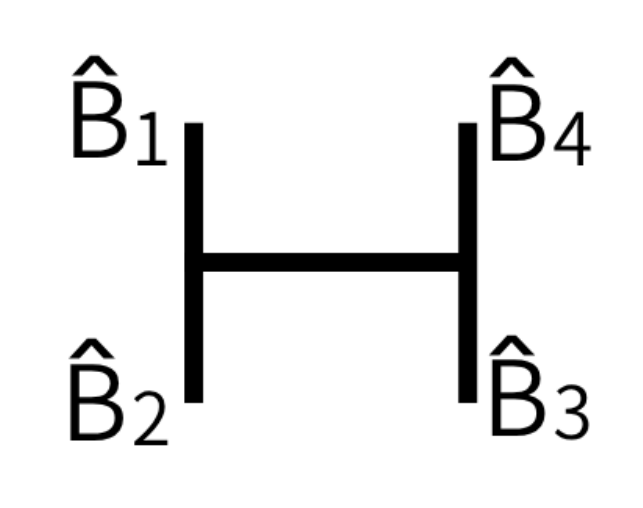}
    \caption{}
    \label{fig:Fig55c}
 \end{minipage}
\end{figure}

This completes the proof of lemma \ref{thm:claim3}.

\end{proof}

Because $f(a)=b$ and $f(b)=a$, $f^2(a)=a$. (See Figure \ref{fig:f(a)=bandf(b)=a} again.) From (1) of proposition \ref{thm:mainprop}, $f^2=1$.
From lemma \ref{thm:claim3}, $f^2\ne 1$. This is a contradiction. 
The proof of (2) of proposition \ref{thm:mainprop} is complete.

(3) Let $x$ be an inner point of 2-cell of $\overline{H}(\Sigma)$.
Assume $f \ne 1$. We want to show a contradiction.

Consider that $a$ is 0-cell of $\overline{H}(\Sigma)$ which is included in the closure of the 2-cell containing $x$.
From (1) and (2) of proposition \ref{thm:mainprop}, every point of 1-skeleton of $\overline{H}(\Sigma)$ is moved by $f$.
Because $f(x)=x$, $f$ moves the 2-cell containing $x$ by ``reflection" or ``rotation" to itself. See Figure \ref{fig:reflectionandrotation}.

\begin{figure}[ht]
    \centering
    \includegraphics[width=0.4\linewidth]{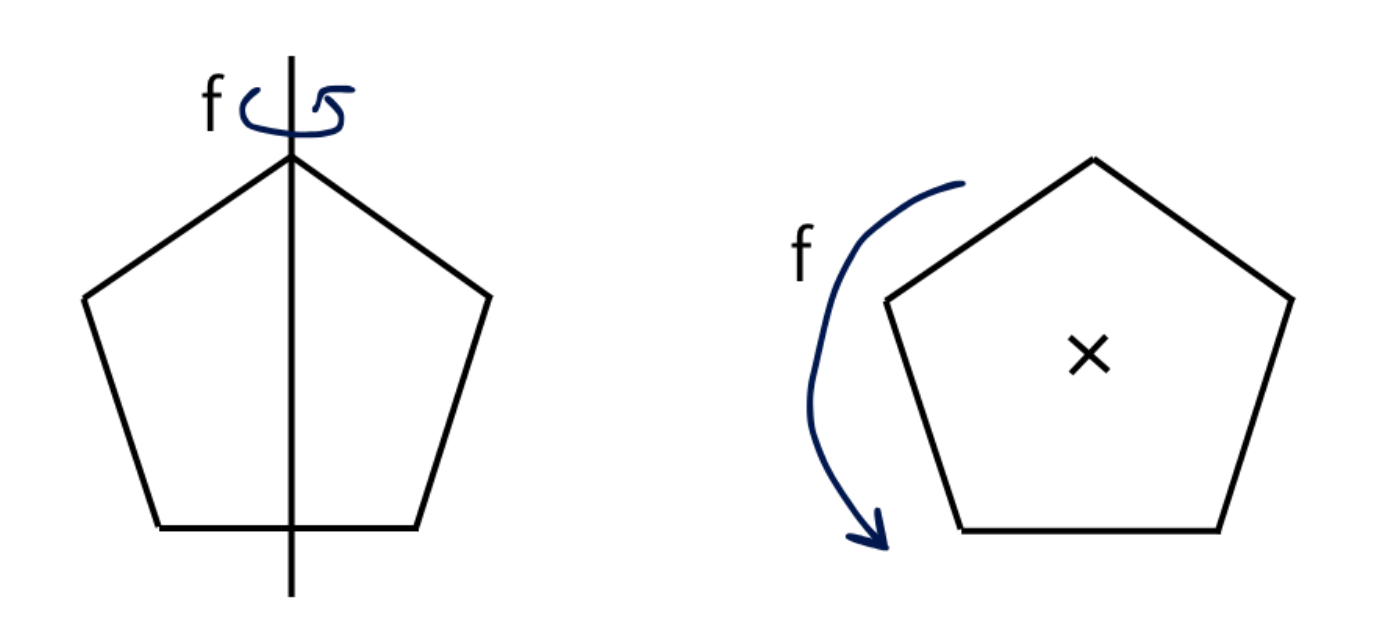}
    \caption{a reflection and a rotation}
    \label{fig:reflectionandrotation}
\end{figure} 

Assume that $f$ moves 2-cell containing $x$ by the reflection.
If the 2-cell is $3A$, $3S$, $6AS$, $2S$ or $2B'3T$, it is oriented as 2-cell. $f$ preserves the orientation. So,
the reflection does not happen. 
Even if the 2-cell is $5A$, $2A$ or $DC$, it is as shown in Figure \ref{fig:reflectionof5Aand2A}, i.e. $f$ fix points on 0-cell or 1-cell.
From (1) and (2) of proposition \ref{thm:mainprop}, this is a contradiction.

\begin{figure}[ht]
    \centering
    \includegraphics[width=0.4\linewidth]{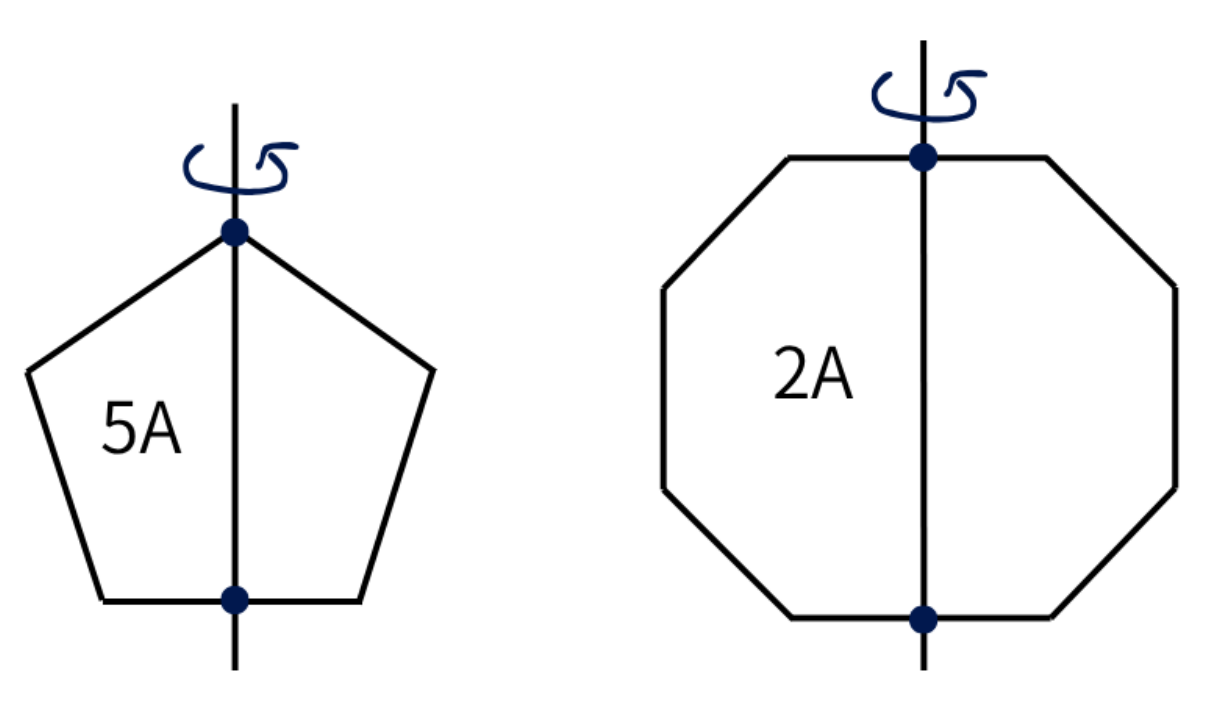}
    \caption{}
    \label{fig:reflectionof5Aand2A}
\end{figure}

We only need to consider the case where $f$ moves the 2-cell containing $x$ by the rotation.

Assume that the 2-cell is $3A$, $3S$, $6AS$, $2S$, $2B'3T$ or $2A$. The 2-cell can not rotate. 
For example, when it is $3S$, see Figure \ref{fig:3Srotation}.

\begin{figure}[ht]
    \centering
    \includegraphics[width=0.2\linewidth]{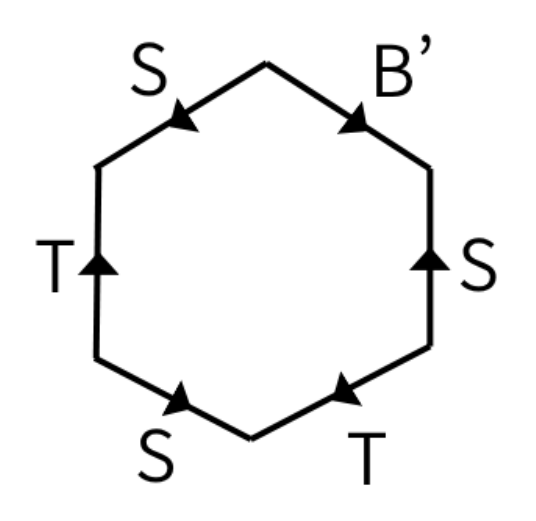}
    \caption{}
    \label{fig:3Srotation}
\end{figure}

Assume that the 2-cell is $5A$. The 2-cell is a pentagon as shown in Figure \ref{fig:5Arotation}.
If it is $1/5$ rotation by $f$, $f^5(a)=a$. From (1) of proposition \ref{thm:mainprop}, $f^5=1$.
And, $a$ and $f(a)$ are 0-cells at both ends of $A$-move. From lemma \ref{thm:claim3}, $f^5\ne 1$. This is a contradiction.

If it is $2/5$ rotation by $f$, $f^2$ is $4/5=-1/5$ rotation.
So, $(f^2)^5(a)=a$. From (1) of proposition \ref{thm:mainprop}, $(f^2)^5=1$.
And, $a$ and $f^2(a)$ are 0-cells at both ends of $A$-move. From lemma \ref{thm:claim3}, $(f^2)^5 \ne 1$.
This is a contradiction.

Assume that the 2-cell is $DC$. There are two moves $\alpha$ and $\beta$ such that the support of $\alpha$ and the support of
 $\beta$ are disjoint as shown in Figure \ref{fig:DCsquare}.
If it is $1/4$ rotation by $f$, $\alpha=\beta$ and $\beta=\alpha^{-1}$ where ``$=$" means equal as a type of moves.
So, $\alpha=\alpha^{-1}$. The only such $\alpha$ is $A$-move. Of course, $\beta$ is also $A$-move.
Then, $f^4(a)=a$ as shown in Figure \ref{fig:DCsquareA}.
From (1) of proposition \ref{thm:mainprop}, $f^4=1$. 
And, $a$ and $f(a)$ are 0-cells at both ends of $A$-move. From lemma \ref{thm:claim3}, $f^4 \ne 1$. This is a contradiction.

\begin{figure}[ht]
 \begin{minipage}[ht]{0.45\linewidth}
    \centering
    \includegraphics[width=0.4\linewidth]{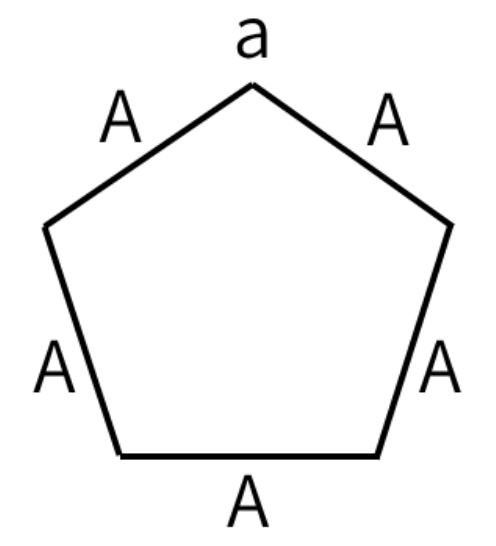}
    \caption{}
    \label{fig:5Arotation}
 \end{minipage}
 \begin{minipage}[ht]{0.45\linewidth}
    \centering
    \includegraphics[width=0.4\linewidth]{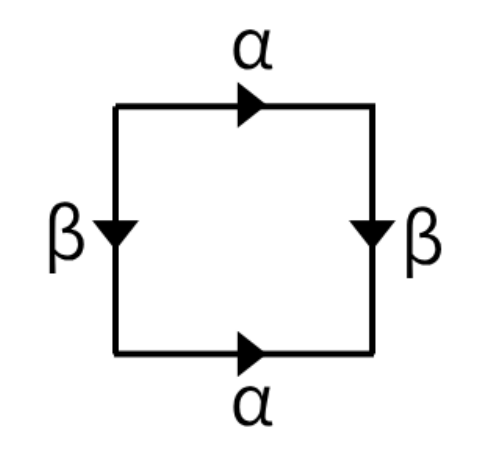}
    \caption{supp$(\alpha)$ $\cap$ supp$(\beta)=\emptyset$}
    \label{fig:DCsquare}
 \end{minipage}
\end{figure}
\begin{figure}[ht]
 \begin{minipage}[ht]{0.45\linewidth}
    \centering
    \includegraphics[width=0.4\linewidth]{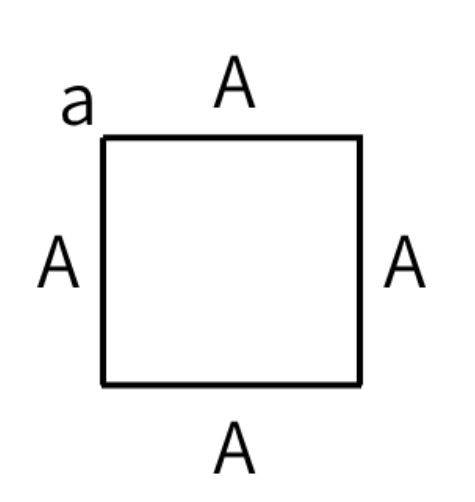}
    \caption{}
    \label{fig:DCsquareA}
 \end{minipage}
 \begin{minipage}[ht]{0.45\linewidth}
    \centering
    \includegraphics[width=0.4\linewidth]{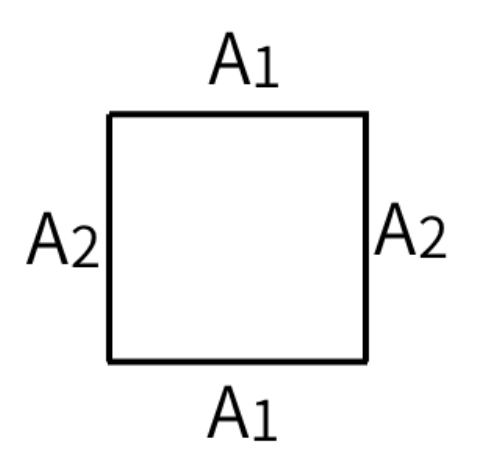}
    \caption{}
    \label{fig:DCsquareA1A2}
 \end{minipage}
\end{figure}
\begin{figure}[ht]
    \centering
    \includegraphics[width=0.8\linewidth]{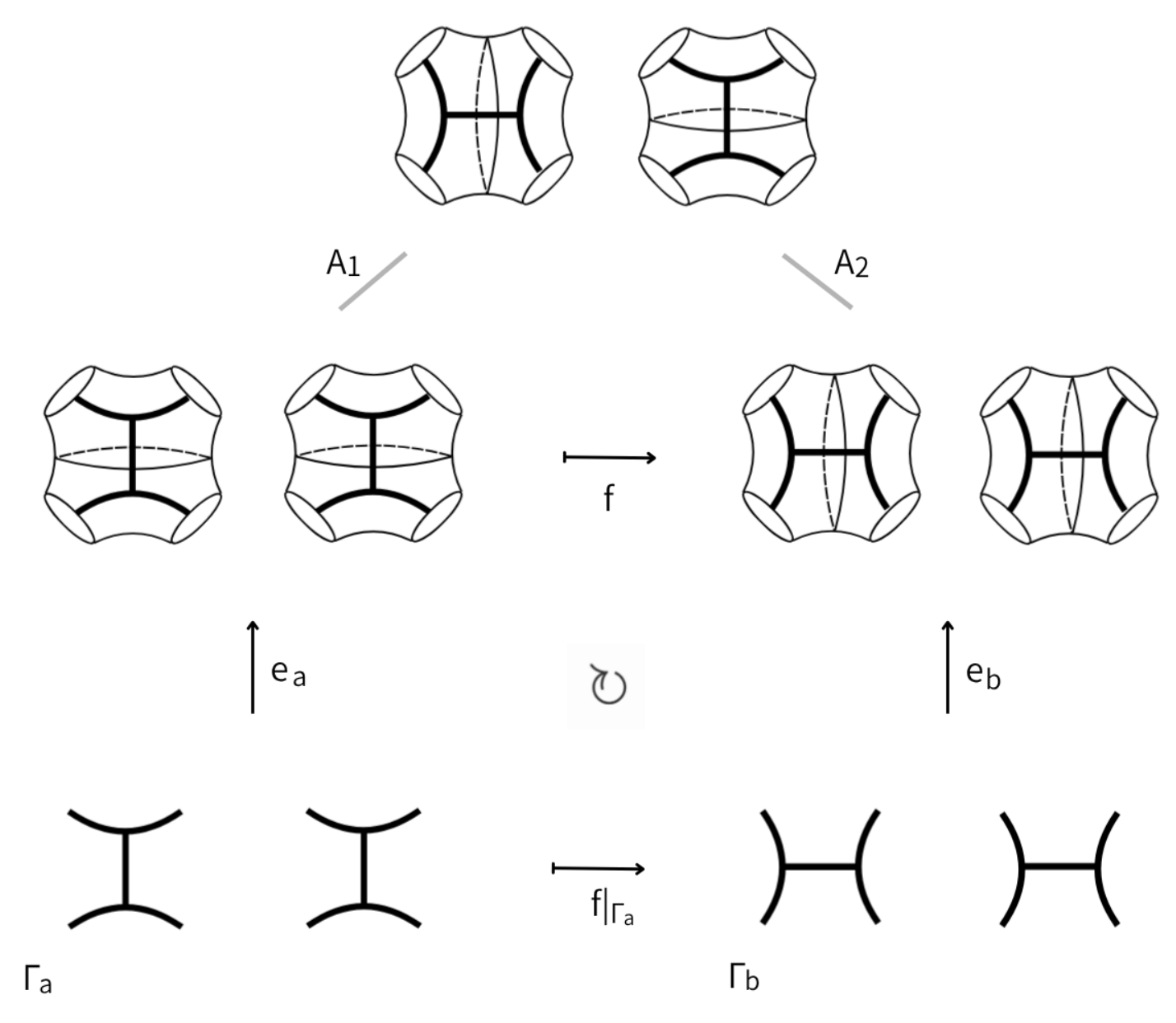}
    \caption{}
    \label{fig:Fig62}
\end{figure}
\clearpage

If it is $2/4$ rotation by $f$, $\alpha=\alpha^{-1}$ and $\beta=\beta^{-1}$. Therefore, both $\alpha$ and $\beta$ are $A$-moves.
Let $A_1$ and $A_2$ be these two $A$-moves. See Figure \ref{fig:DCsquareA1A2}. The support of $A_1$ and the support of $A_2$ are disjoint.

\begin{lemma}\label{thm:claim4}
    $f^2\ne 1$ for above $f$. 
\end{lemma}
\begin{proof}

We will show it by the same manner as in the proof of lemma \ref{thm:claim3}.


Let $f(a)=b$. $a$ and $b$ are colored quilt decompositions.
Let $\Gamma_a$ and $\Gamma_b$ be trivalent graphs with $r$ tails corresponding to $a$ and $b$.
Let $e_a$ and $e_b$ be embeddings from $\Gamma_a$ and $\Gamma_b$ to $\Sigma$. 
Then, we define that $f|_{\Gamma_a}:=e_b^{-1}\circ f\circ e_a$ (see Figure \ref{fig:Fig62}).
Remark that $f|_{\Gamma_a}:\Gamma_a\cong\Gamma_b$ as a trivalent graph.

From the same discussion as the proof of lemma \ref{thm:claim3}, the following holds:

\begin{equation}
    (f|_{\Gamma_a})^*=(f^*)|_{H_1(\Gamma_a)}          \tag{3.2}\label{2}
\end{equation}

(I) Assume that $\Gamma_a$ is connected as shown in Figure \ref{fig:Fig63a}.
That is, there exists a closed path $p$ which is $n$-gon ($n\geq 1$) in $\Gamma_a$.
Then, $\Gamma_b$ has a closed path $p'$ which is $(n+1)$-gon. See Figure \ref{fig:Fig63b}.
From \eqref{2}, $f|_{\Gamma_a}(p)=p'$. Such $f$ does not exist.

\begin{figure}[ht]
 \begin{minipage}[ht]{0.45\linewidth}
    \centering
    \includegraphics[width=0.5\linewidth]{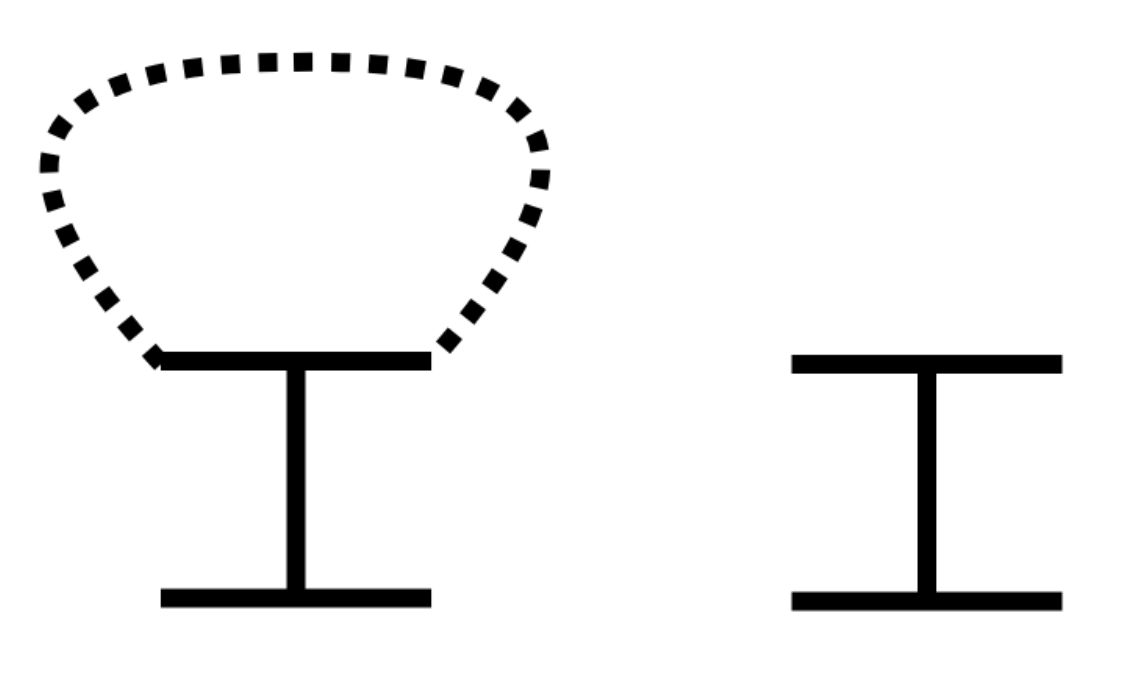}
    \caption{}
    \label{fig:Fig63a}
 \end{minipage}
 \begin{minipage}[ht]{0.45\linewidth}
    \centering
    \includegraphics[width=0.5\linewidth]{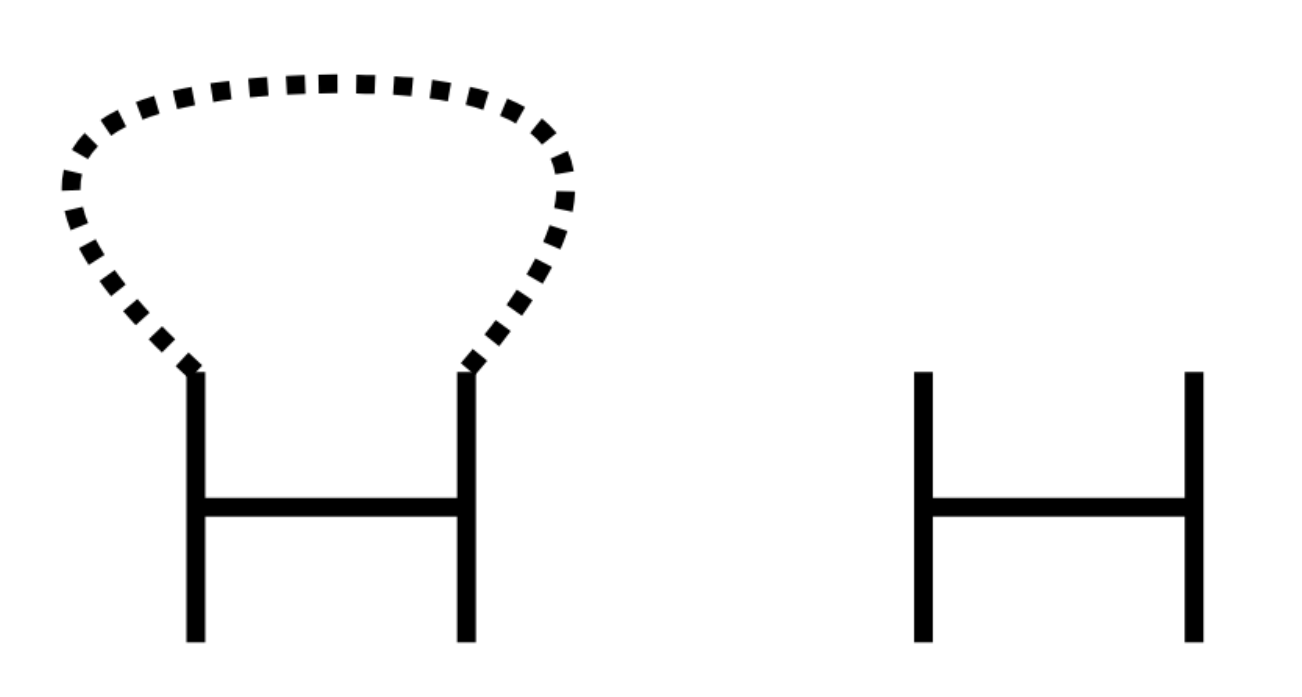}
    \caption{}
    \label{fig:Fig63b}
 \end{minipage}
\end{figure}

(II) Assume that $\Gamma_a$ is connected as shown in Figure \ref{fig:Fig64}.
From the same discussion as (II) of the proof of lemma \ref{thm:claim3}, such $f$ does not exist.

\begin{figure}[ht]
    \centering
    \includegraphics[width=0.25\linewidth]{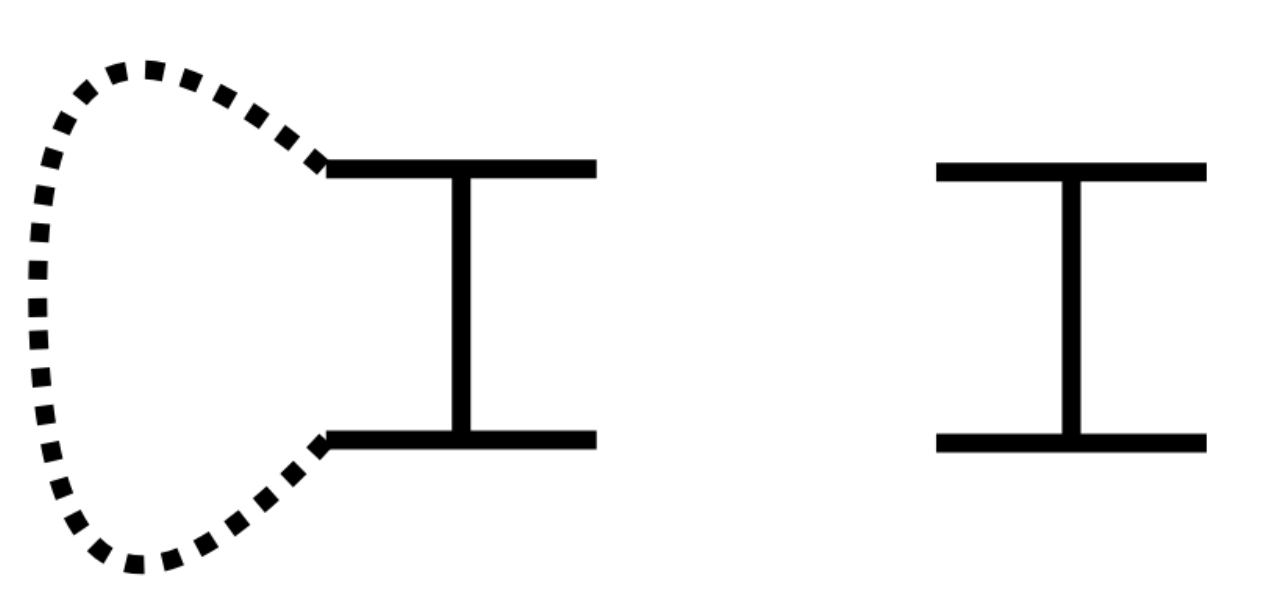}
    \caption{}
    \label{fig:Fig64}
\end{figure}

(III) Assume that $\Gamma_a$ is connected as shown in Figure \ref{fig:Fig65}.

\begin{figure}[ht]
    \centering
    \includegraphics[width=0.22\linewidth]{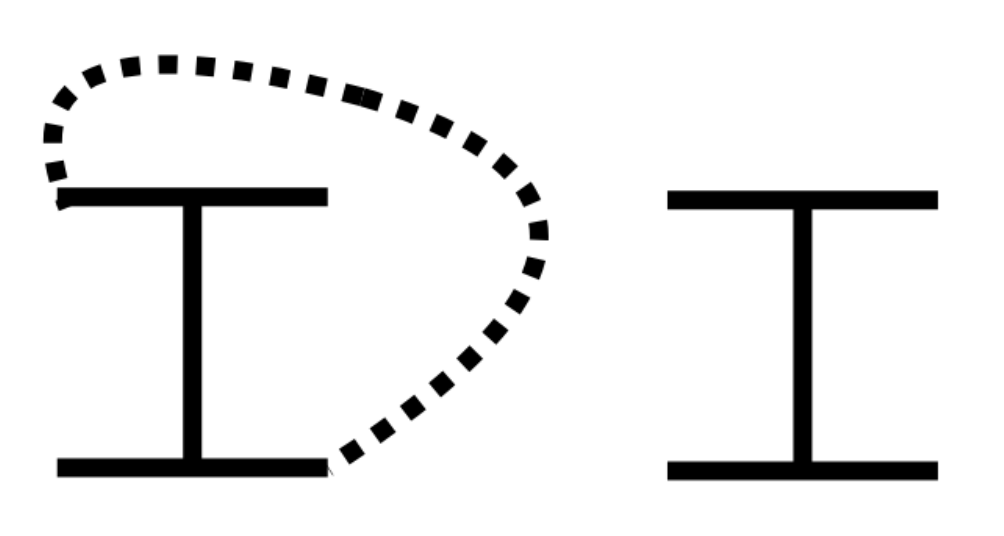}
    \caption{}
    \label{fig:Fig65}
\end{figure}

From the same discussion as (III) of the proof of lemma \ref{thm:claim3}, there are the following possibilities.

(III)-(i) Assume that $\Gamma_a$ is further connected as shown in Figure \ref{fig:Fig66a}.
This includes cases like the one shown in Figure \ref{fig:Fig66b}, for example.
From the same discussion as (III)-(i) of the proof of lemma \ref{thm:claim3}, such $f$ does not exist.

\begin{figure}[ht]
 \begin{minipage}[ht]{0.45\linewidth}
    \centering
    \includegraphics[width=0.5\linewidth]{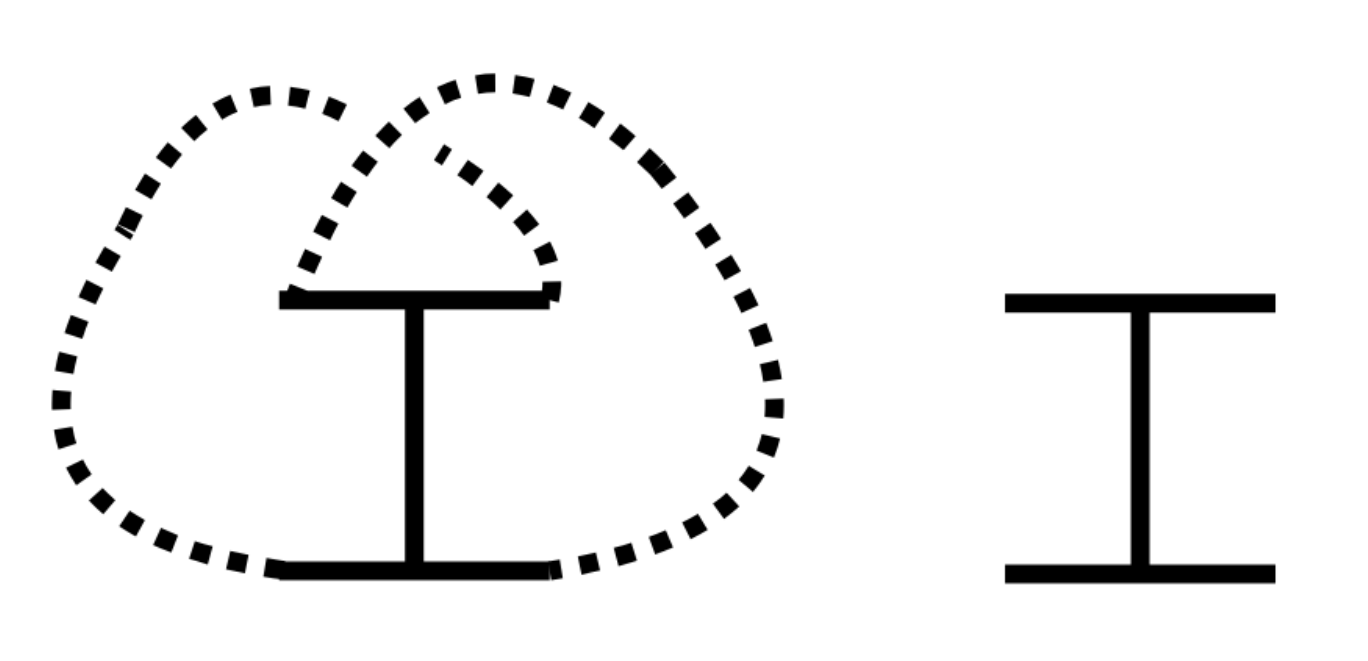}
    \caption{}
    \label{fig:Fig66a}
 \end{minipage}
 \begin{minipage}[ht]{0.45\linewidth}
    \centering
    \includegraphics[width=0.4\linewidth]{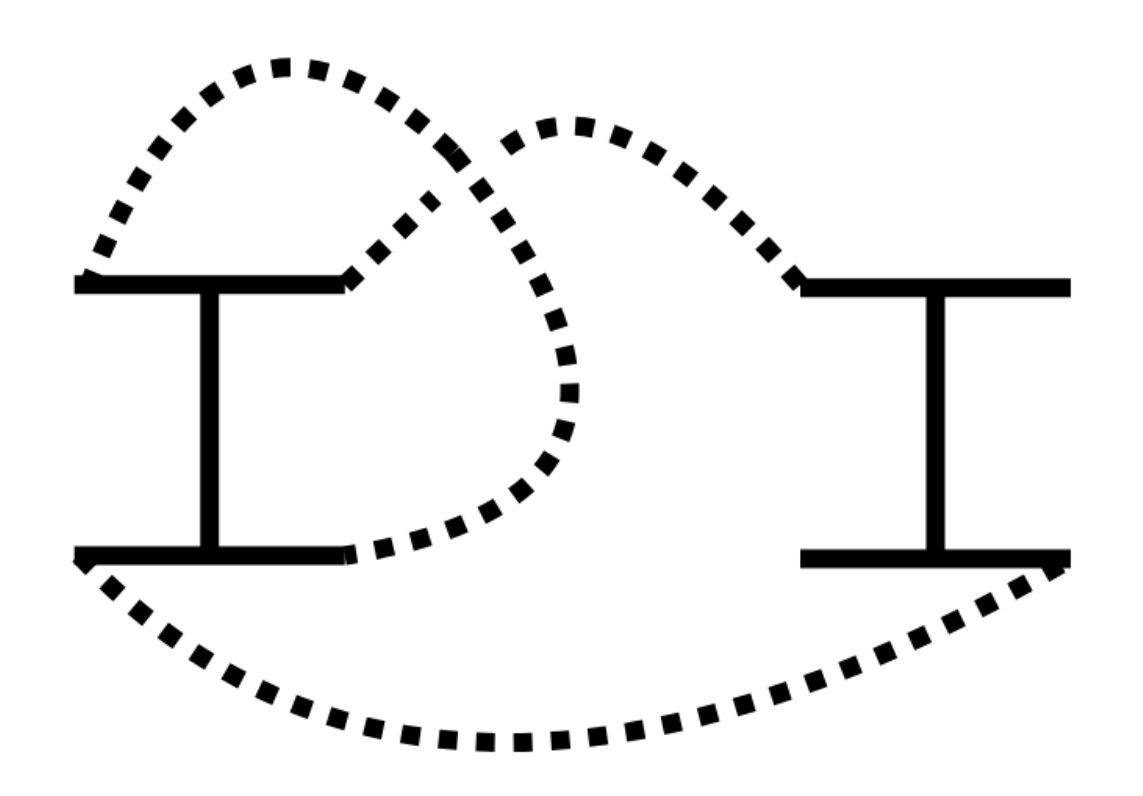}
    \caption{}
    \label{fig:Fig66b}
 \end{minipage}
\end{figure}

(III)-(ii) Assume that $\Gamma_a$ is disconnected as shown in Figure \ref{fig:Fig67a}.
From (I) and (II), those two edges are not connected to anything else.
So, the two edges correspond to BSCCs. See Figure \ref{fig:Fig67a} again. $B$ means BSCC.
From the same discussion as (III)-(ii) of the proof of lemma \ref{thm:claim3}, $\Gamma_a$ is as shown in Figure \ref{fig:Fig67b}.

\begin{figure}[ht]
 \begin{minipage}[ht]{0.45\linewidth}
    \centering
    \includegraphics[width=0.6\linewidth]{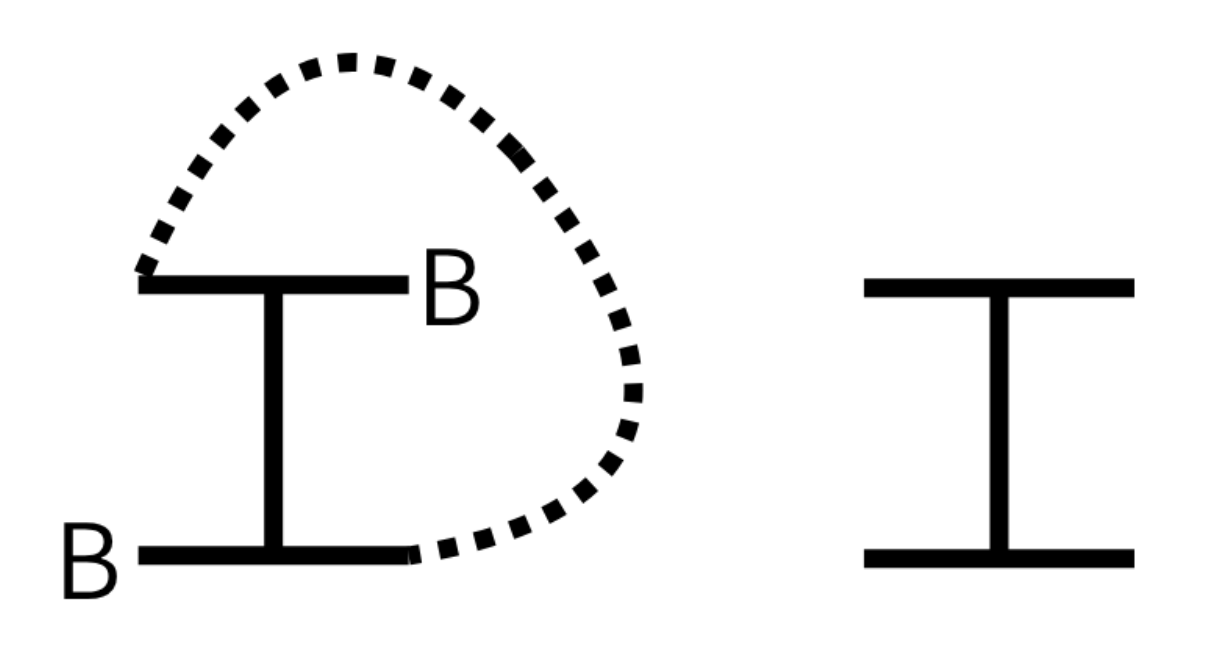}
    \caption{}
    \label{fig:Fig67a}
 \end{minipage}
 \begin{minipage}[ht]{0.45\linewidth}
    \centering
    \includegraphics[width=0.6\linewidth]{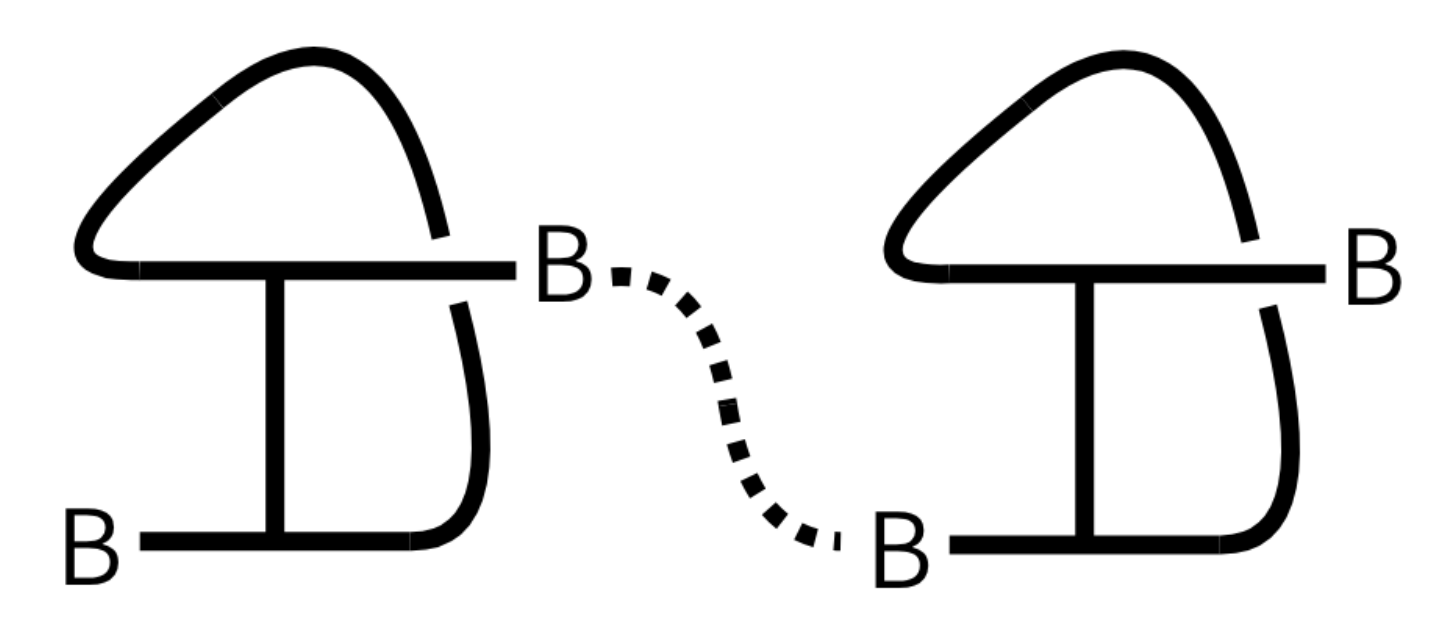}
    \caption{}
    \label{fig:Fig67b}
 \end{minipage}
\end{figure}

The colored quilt decomposition $a$ on $\Sigma$ is as shown in Figure \ref{fig:Fig67c}.
Let $f_1,f_2\in \mathcal{M}_{g,r}$ be as shown in Figure \ref{fig:Fig67d}, where ``$+$" means a positive Dehn twist along the SCC and 
 ``$-$" means a negative Dehn twist along the SCC.
 Then, one can see that $f=f_2\circ f_1$ and $f(a)=b$. 
 (Note that $f_1,f_2\in\mathcal{I}_{g,r}$ and hence $f\in\mathcal{I}_{g,r}$.)
It is not hard to check that $f^2\ne 1$.

\begin{figure}[ht]
    \centering
    \includegraphics[width=0.7\linewidth]{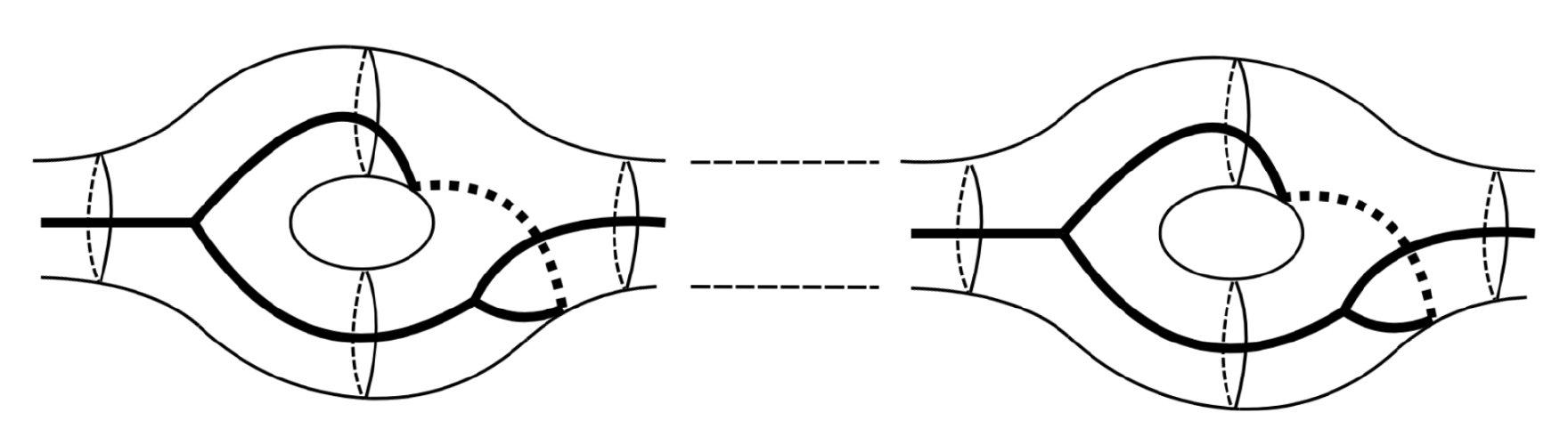}
    \caption{}
    \label{fig:Fig67c}
\end{figure}

\begin{figure}[ht]
    \centering
    \includegraphics[width=0.7\linewidth]{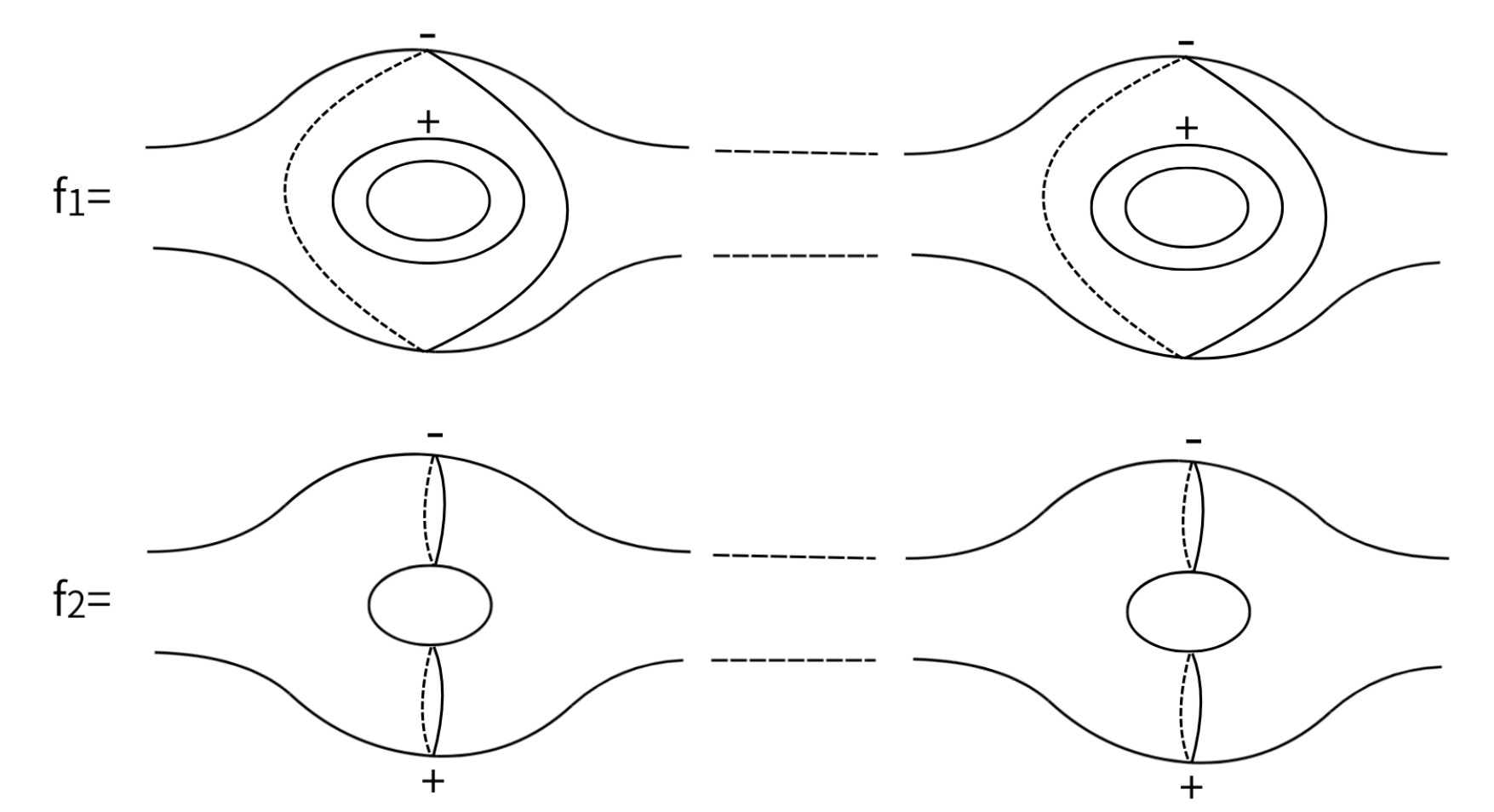}
    \caption{}
    \label{fig:Fig67d}
\end{figure}


(IV) Assume that the four edges of $\Gamma_a$ is all disconnected as shown in Figure \ref{fig:Fig68}.
Then, all four of those edges correspond to BSCCs. 

\begin{figure}[ht]
    \centering
    \includegraphics[width=0.3\linewidth]{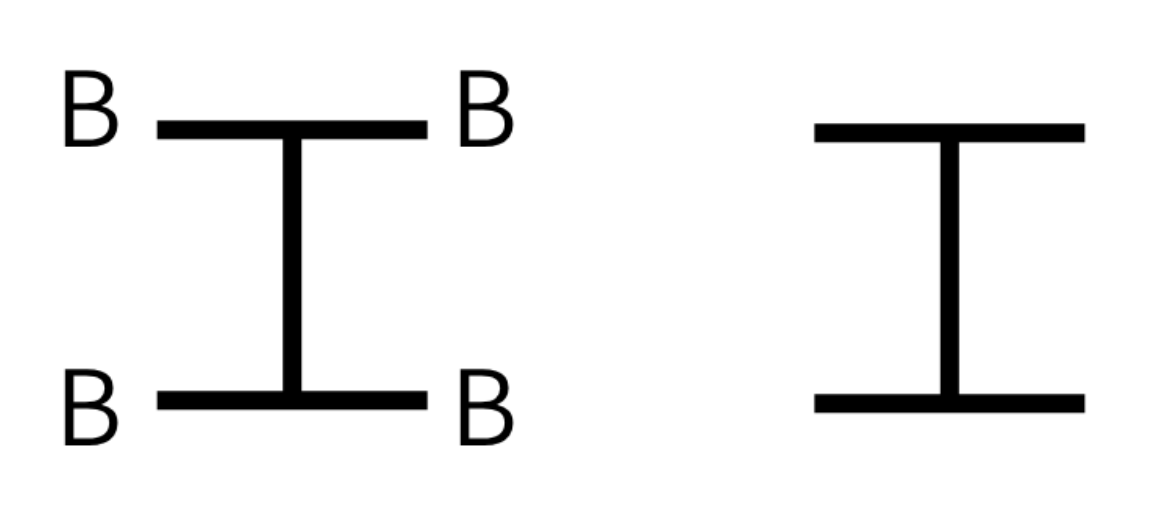}
    \caption{}
    \label{fig:Fig68}
\end{figure}

From the same discussion as (IV) of the proof of lemma \ref{thm:claim3}, such $f$ does not exist. 

(V) Assume that $\Gamma_a$ is connected as shown in Figure \ref{fig:Fig69a}.
That is, there exists a closed path $p$ which is $n$-gon ($n\geq 4$) in $\Gamma_a$.
Then, $\Gamma_b$ has a closed path $p'$ which is $(n-2)$-gon. See Figure \ref{fig:Fig69b}.
From \eqref{2}, $f|_{\Gamma_a}(p)=p'$. Such $f$ does not exist.


\begin{figure}[ht]
 \begin{minipage}[ht]{0.45\linewidth}
    \centering
    \includegraphics[width=0.6\linewidth]{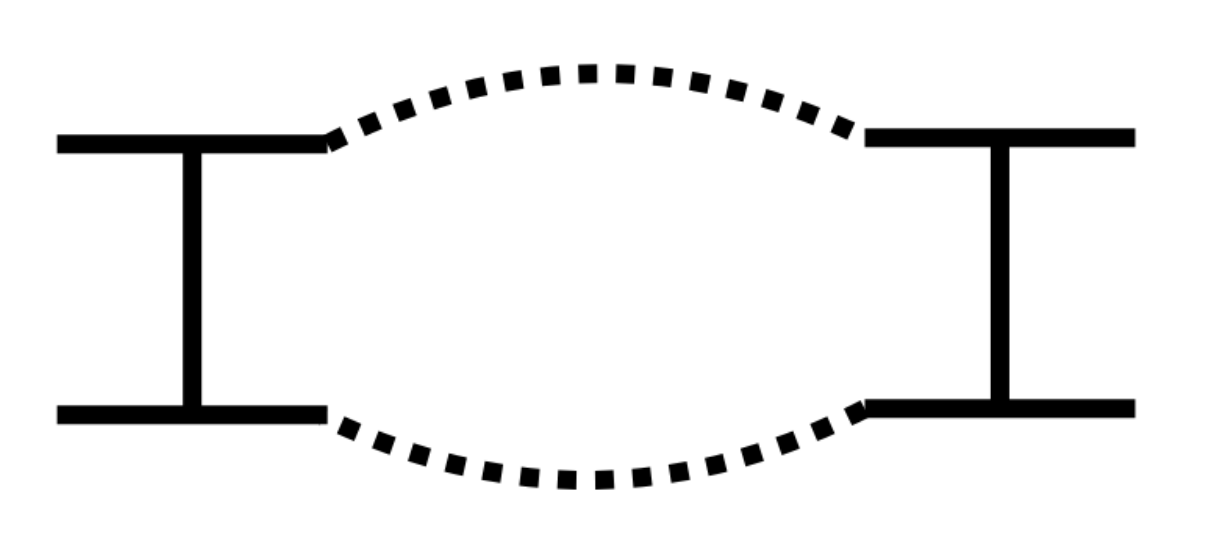}
    \caption{}
    \label{fig:Fig69a}
 \end{minipage}
 \begin{minipage}[ht]{0.45\linewidth}
    \centering
    \includegraphics[width=0.6\linewidth]{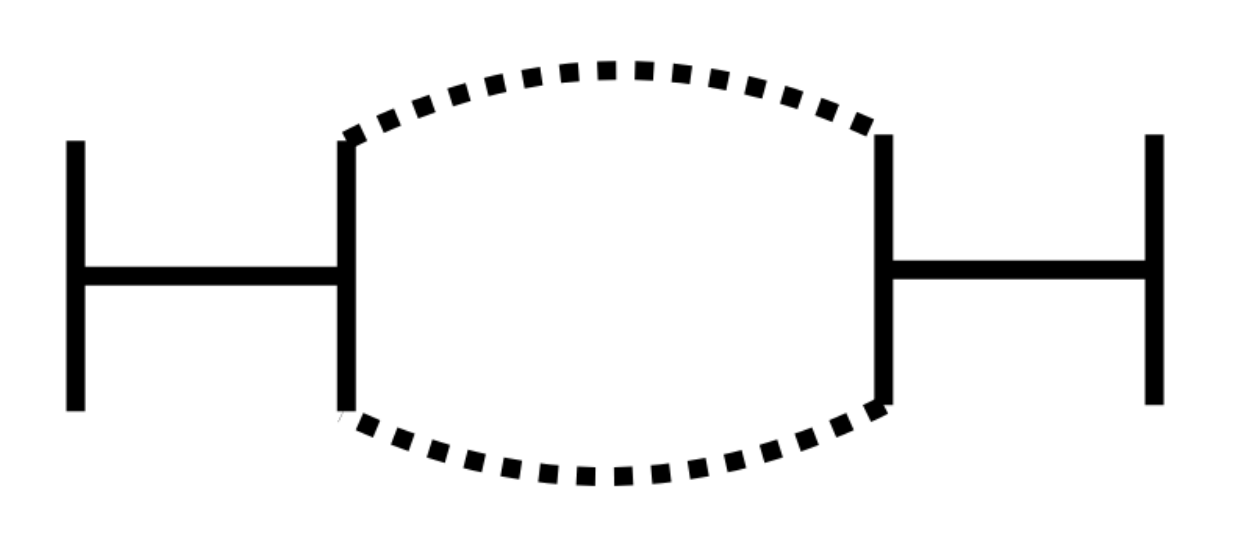}
    \caption{}
    \label{fig:Fig69b}
 \end{minipage}
\end{figure}

(VI) Assume that $\Gamma_a$ is connected as shown in Figure \ref{fig:Fig70}.
Such $f$ does not exist.

\begin{figure}[ht]
    \centering
    \includegraphics[width=0.3\linewidth]{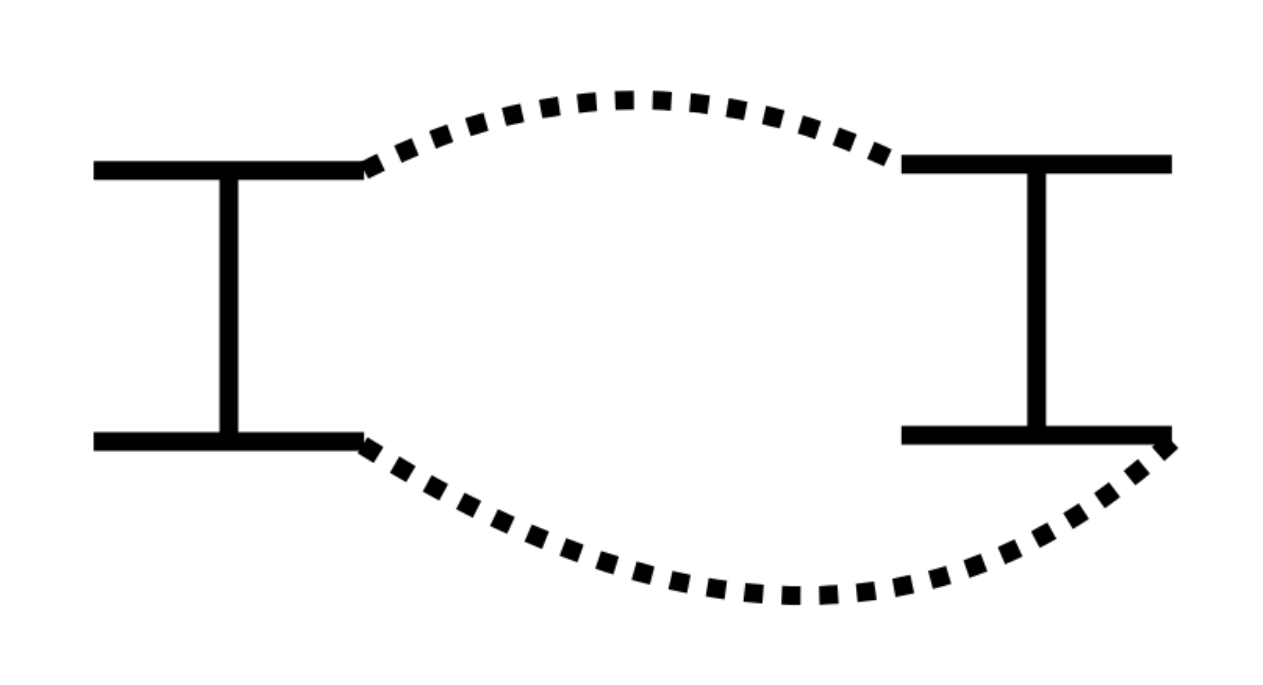}
    \caption{}
    \label{fig:Fig70}
\end{figure}

(VII) Assume that $\Gamma_a$ is connected as shown in Figure \ref{fig:Fig71}.
There exists a closed path $p$ as shown in the figure.

\begin{figure}[ht]
    \centering
    \includegraphics[width=0.3\linewidth]{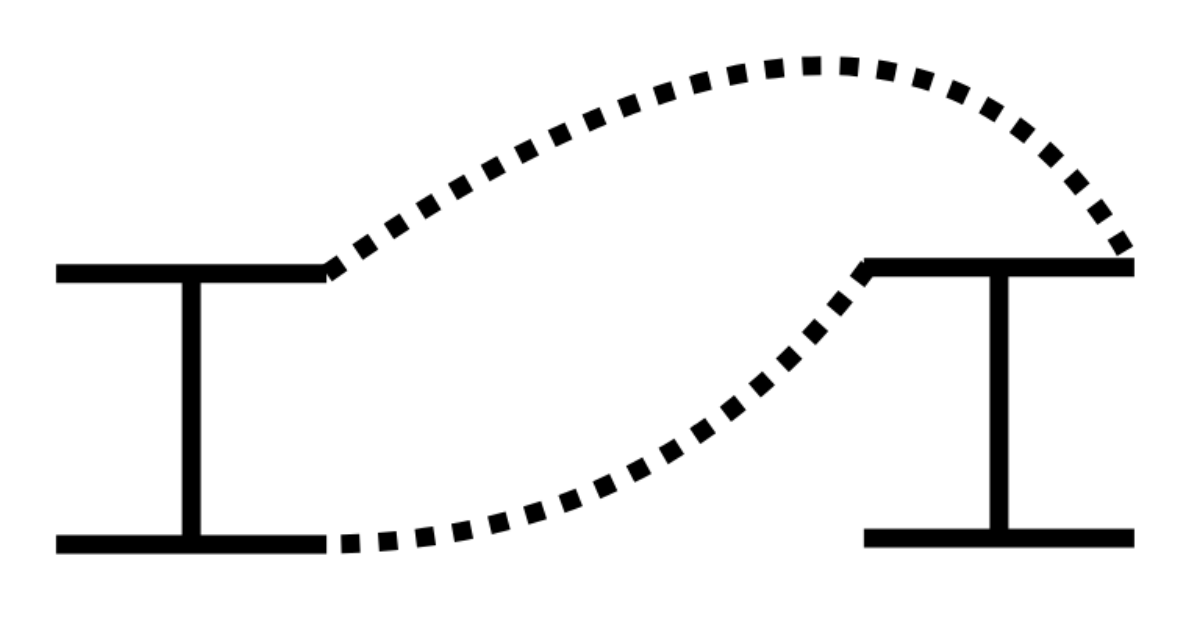}
    \caption{}
    \label{fig:Fig71}
\end{figure}

(VII)-(i) Assume that all the remaining four edges correspond to BSCC.
Denote that those four edges are $B_1$,$B_2$,$B_3$ and $B_4$ as shown in Figure \ref{fig:Fig72a}.
Then, $\Gamma_b$ is as shown in Figure \ref{fig:Fig72b} where $\widehat{B_k}=f|_{\Gamma_a}(B_k)$ for $k=1,2,3,4$.
Note that $B_3$ and $B_4$ share a vertex, but $\widehat{B_3}$ and $\widehat{B_4}$ do not.
Such $f$ does not exist.

\begin{figure}[H]
 \begin{minipage}[ht]{0.45\linewidth}
    \centering
    \includegraphics[width=0.7\linewidth]{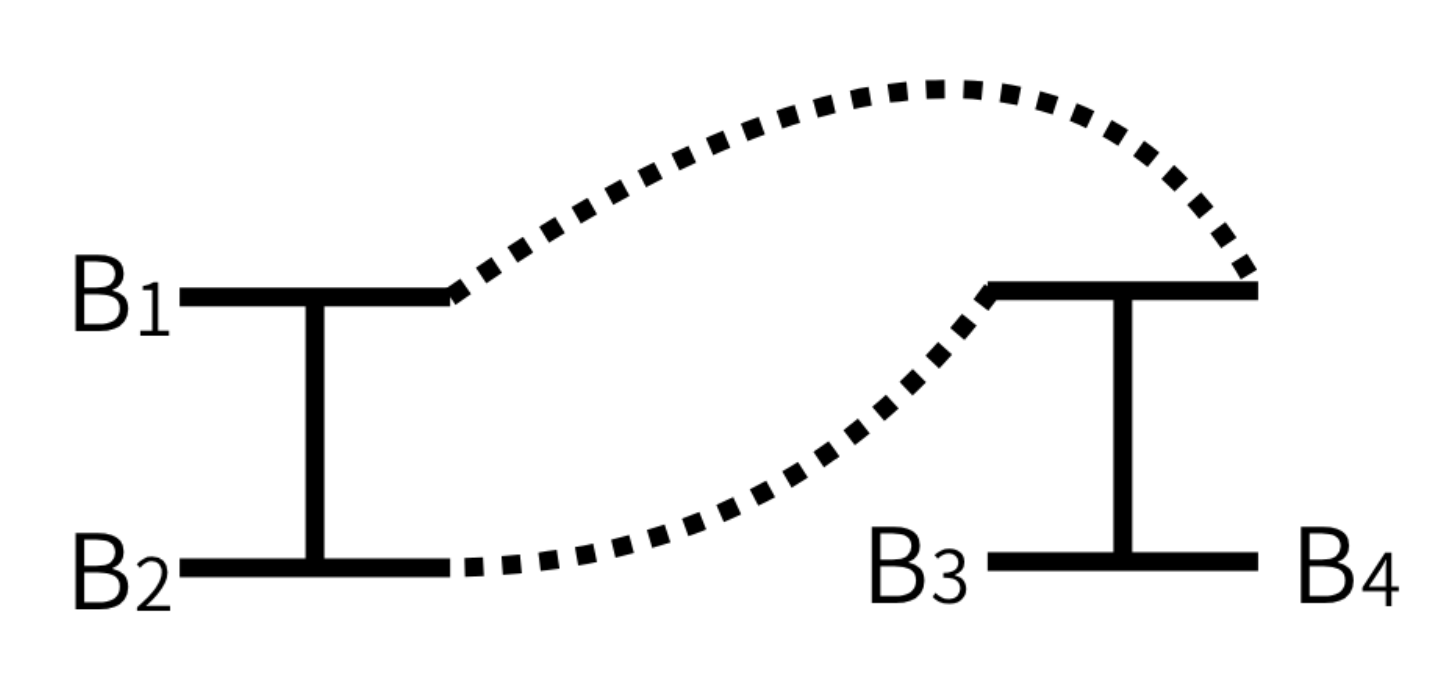}
    \caption{}
    \label{fig:Fig72a}
 \end{minipage}
 \begin{minipage}[ht]{0.45\linewidth}
    \centering
    \includegraphics[width=0.7\linewidth]{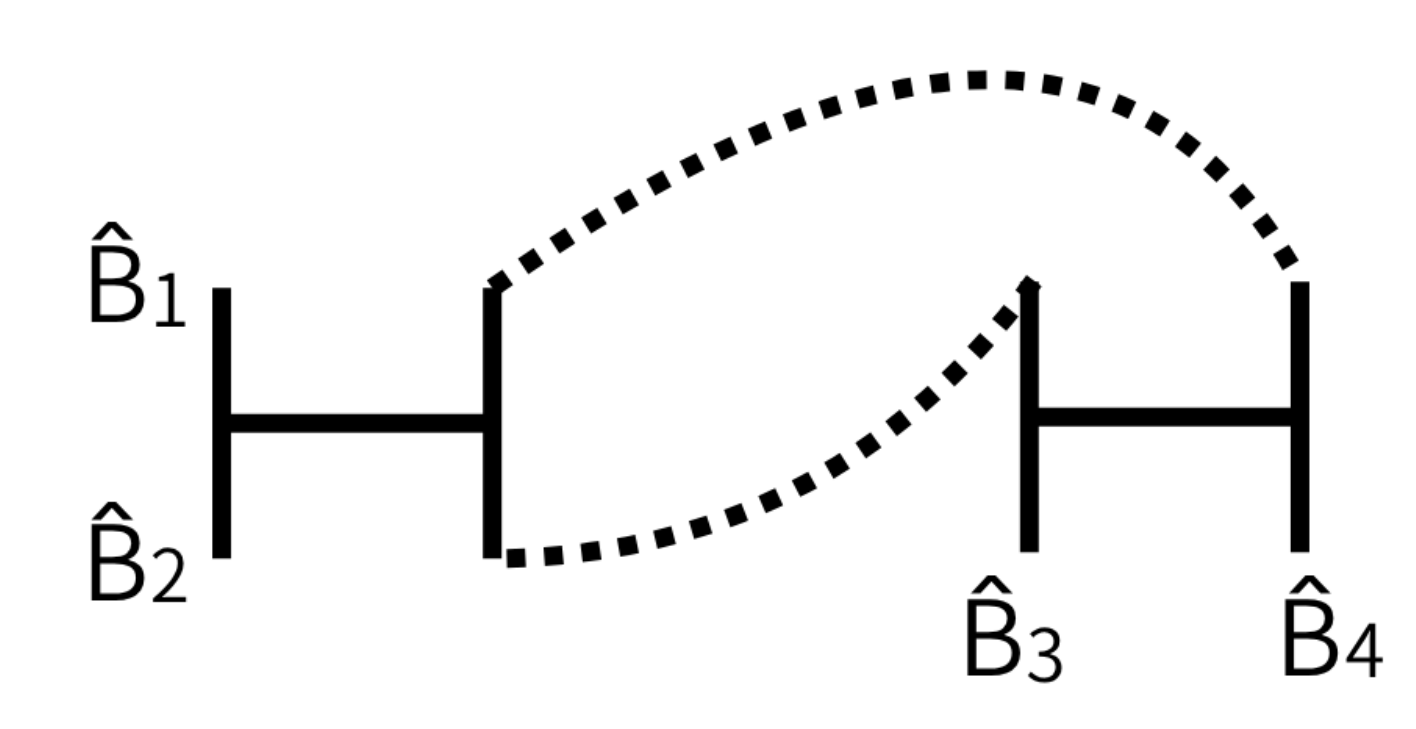}
    \caption{}
    \label{fig:Fig72b}
 \end{minipage}
\end{figure}

\begin{figure}[H]
 \begin{minipage}[ht]{0.45\linewidth}
    \centering
    \includegraphics[width=0.7\linewidth]{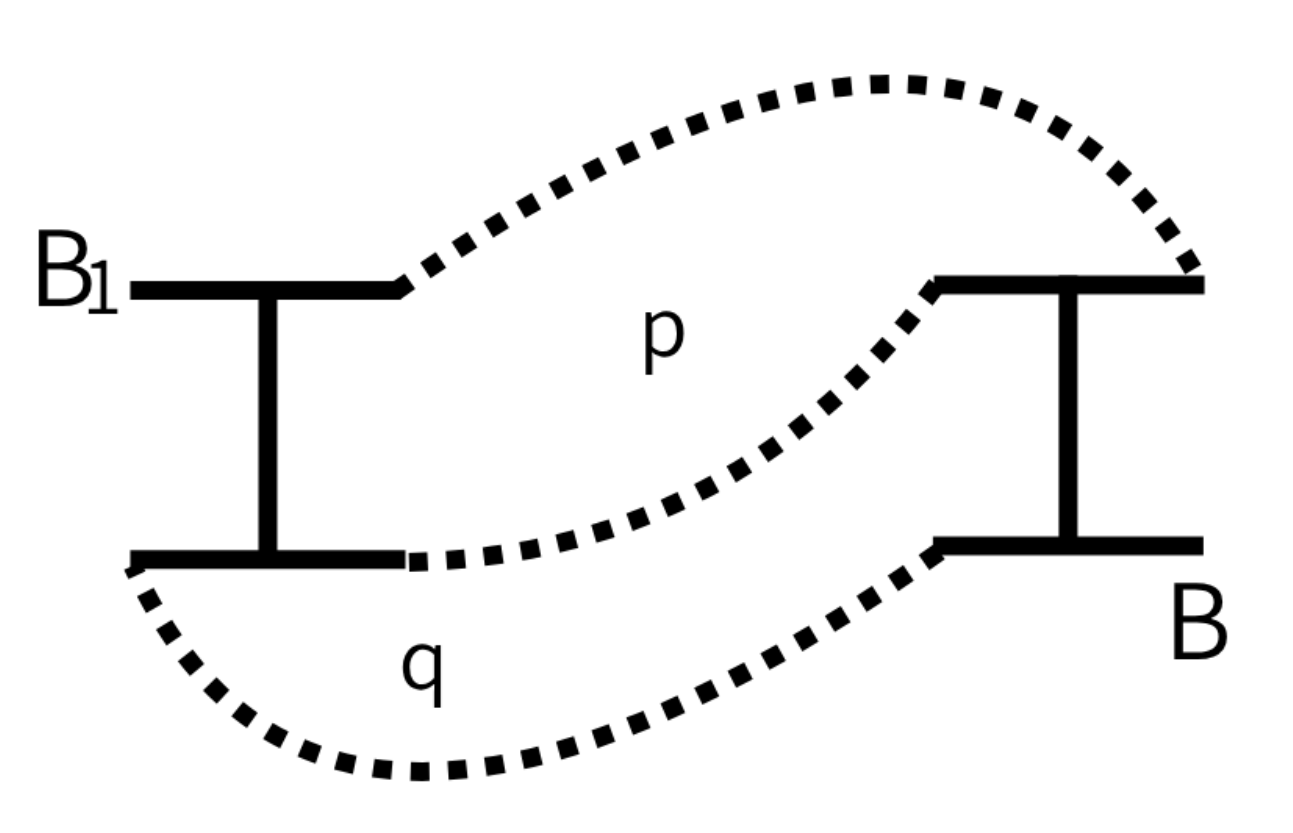}
    \caption{}
    \label{fig:Fig73a}
 \end{minipage}
 \begin{minipage}[ht]{0.45\linewidth}
    \centering
    \includegraphics[width=0.7\linewidth]{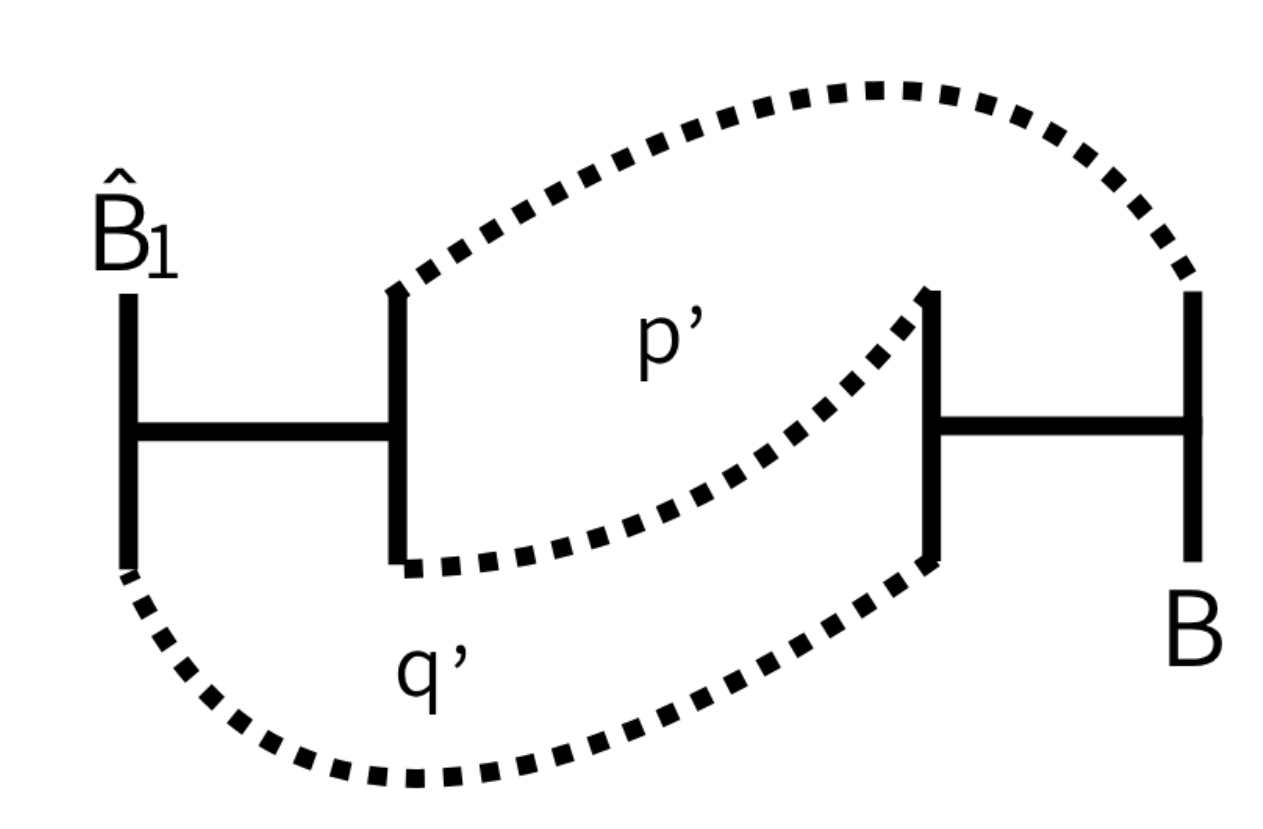}
    \caption{}
    \label{fig:Fig73b}
 \end{minipage}
\end{figure}

(VII)-(ii) Assume that $\Gamma_a$ is further connected as shown in Figure \ref{fig:Fig73a}.
That is, there exists another closed path $q$.
Then, $\Gamma_b$ has closed paths $p'$ and $q'$. See Figure \ref{fig:Fig73b}.
Note that $B_1$ and $p$ share a vertex, but $\widehat{B_1}$ and $p'$ do not.
Such $f$ does not exist.

(VII)-(iii) Assume that $\Gamma_a$ is further connected as shown in Figure \ref{fig:Fig74a}.
That is, there exists another closed path $r$.
Then, $\Gamma_b$ has another closed path $r'$.
Let $\widetilde{r}$ be 1-cycle of $r$, and let $\widetilde{r'}$ be 1-cycle of $r'$, such that $f$ maps $\widetilde{r}$ to $\widetilde{r'}$. 
In Figure \ref{fig:Fig74a}, let $\alpha$,$\beta$ and $\gamma$ be $1$, $m$ and $n$ edges ($m\geq 1$, $n\geq 1$).
Then, $f(\alpha)$, $f(\beta)$ and $f(\gamma)$ are as shown in Figure \ref{fig:Fig74b}.
Remark that $r$ share one edge $\alpha$ with $p$.
In this case, the cyclic order of $\alpha$, $\beta$, $\gamma$ induced by the orientation of $\widetilde{r}$ is different
from the cyclic order of $f(\alpha)$, $f(\beta)$, $f(\gamma)$
induced by the orientation of $\widetilde{r'}$.
Since $f|_{{\Gamma}_a}$ is an isomorphism of the graph,
such $f$ does not exist.

\begin{figure}[H]
 \begin{minipage}[ht]{0.45\linewidth}
    \centering
    \includegraphics[width=0.8\linewidth]{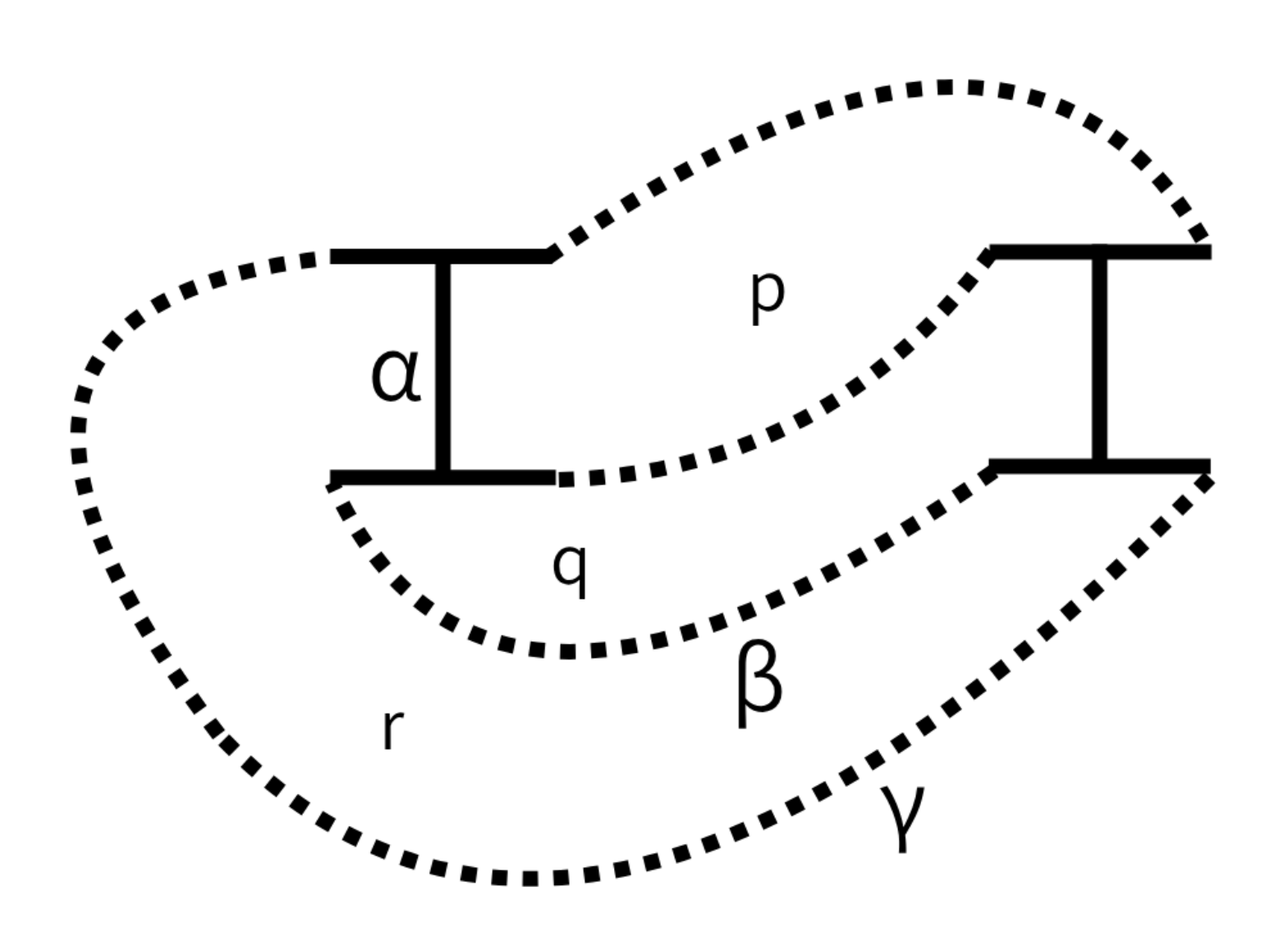}
    \caption{}
    \label{fig:Fig74a}
 \end{minipage}
 \begin{minipage}[ht]{0.45\linewidth}
    \centering
    \includegraphics[width=0.8\linewidth]{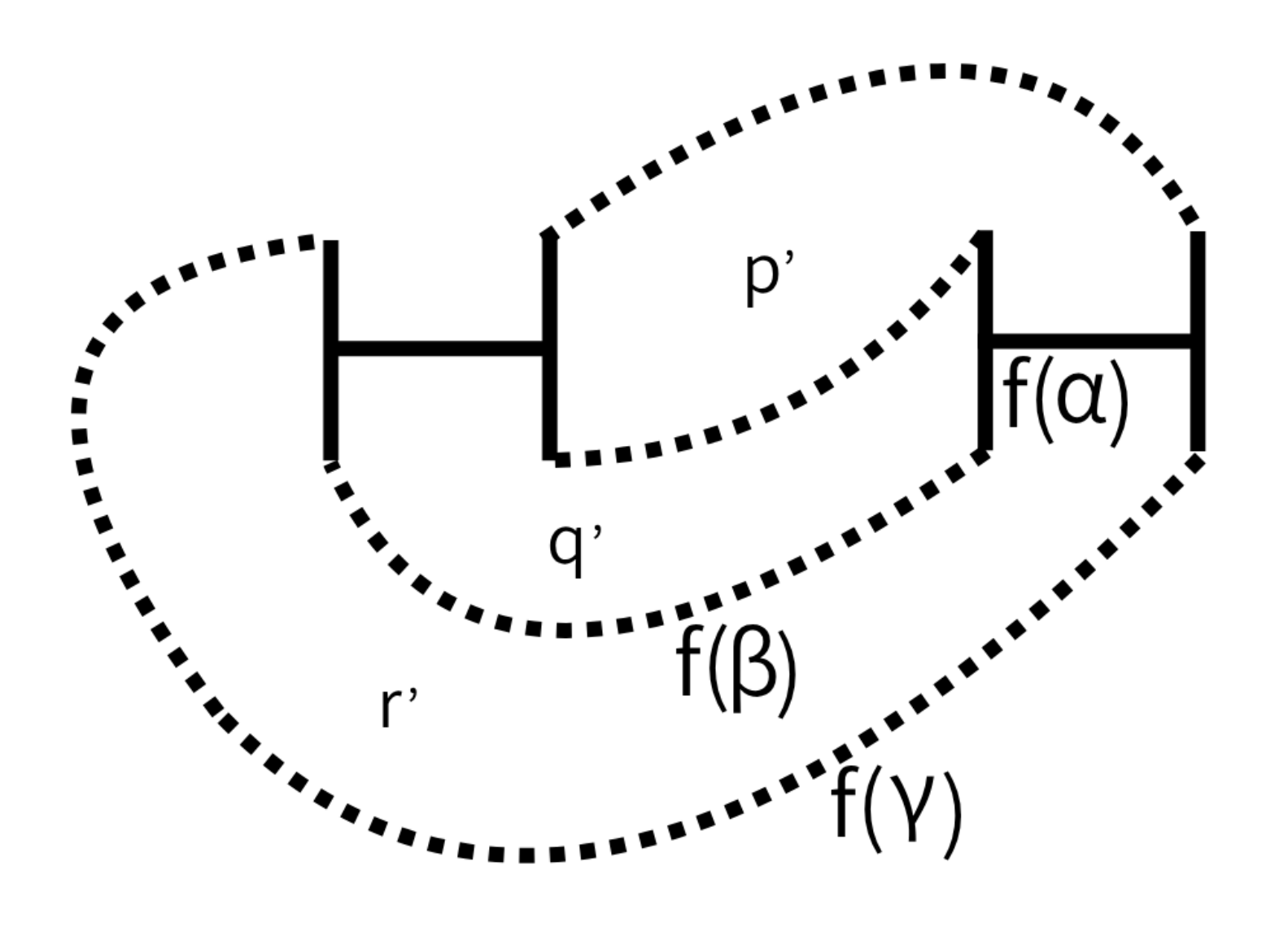}
    \caption{}
    \label{fig:Fig74b}
 \end{minipage}
\end{figure}


(VIII) Assume that $\Gamma_a$ is connected as shown in Figure \ref{fig:Fig75a}.
From the discussion so far, the remaining four edges of $\Gamma_a$ are not connected to anything, hence,
all the four edges correspond to BSCCs. Let them be $B_1$, $B_2$, $B_3$, and $B_4$.
See Figure \ref{fig:Fig75b}.
Let $p$ be the closed path of $\Gamma_a$, and let $\widetilde{p}$ be 1-cycle of $p$. 
Set $\widehat{B_k}=f(B_k)$ for $k=1,2,3,4$. Let $p'$ be the closed path of $\Gamma_b$. 
Let $\widetilde{p'}$ be 1-cycle of $p'$ such that $f$ maps $\widetilde{p}$ to $\widetilde{p'}$.
Then, it is as shown in Figure \ref{fig:Fig75c}.
In this case, the cyclic order of $B_1$,$B_2$,$B_3$,$B_4$ induced by the orientation of $\widetilde{p}$ is different from the cyclic order of
 $\widehat{B_1}$,$\widehat{B_2}$,$\widehat{B_3}$,$\widehat{B_4}$ induced by the orientation of $\widetilde{p'}$.
Such $f$ does not exist.

\begin{figure}[H]
 \begin{minipage}[ht]{0.3\linewidth}
    \centering
    \includegraphics[width=0.8\linewidth]{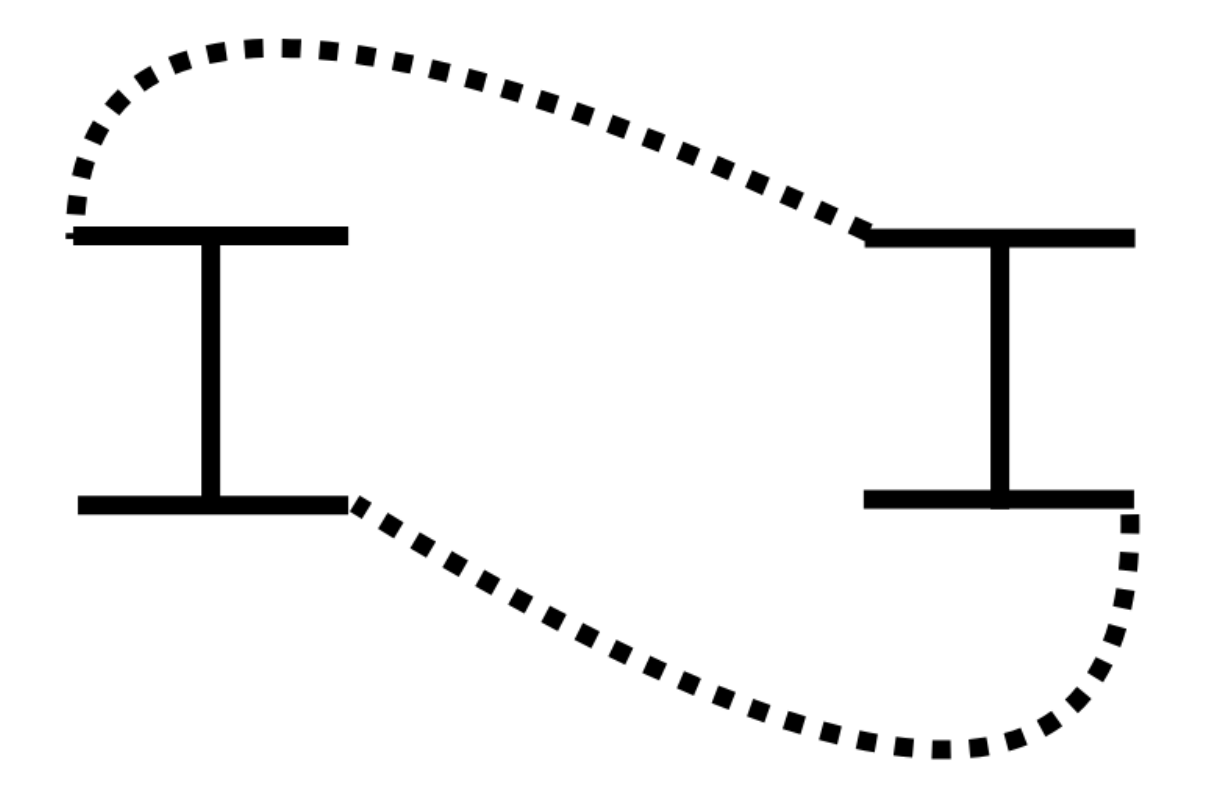}
    \caption{}
    \label{fig:Fig75a}
 \end{minipage}
 \begin{minipage}[ht]{0.3\linewidth}
    \centering
    \includegraphics[width=0.8\linewidth]{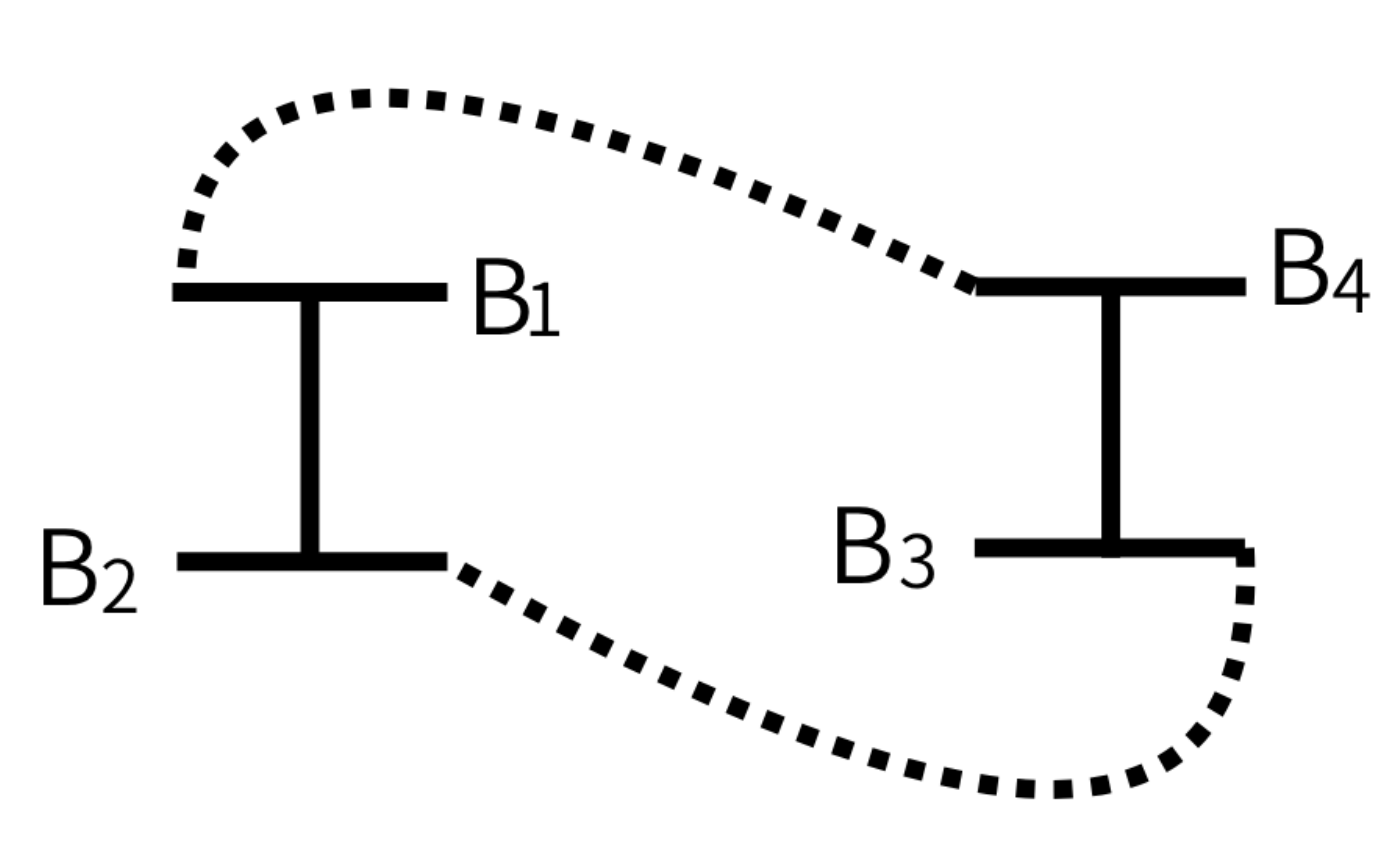}
    \caption{}
    \label{fig:Fig75b}
 \end{minipage}
 \begin{minipage}[ht]{0.3\linewidth}
    \centering
    \includegraphics[width=0.8\linewidth]{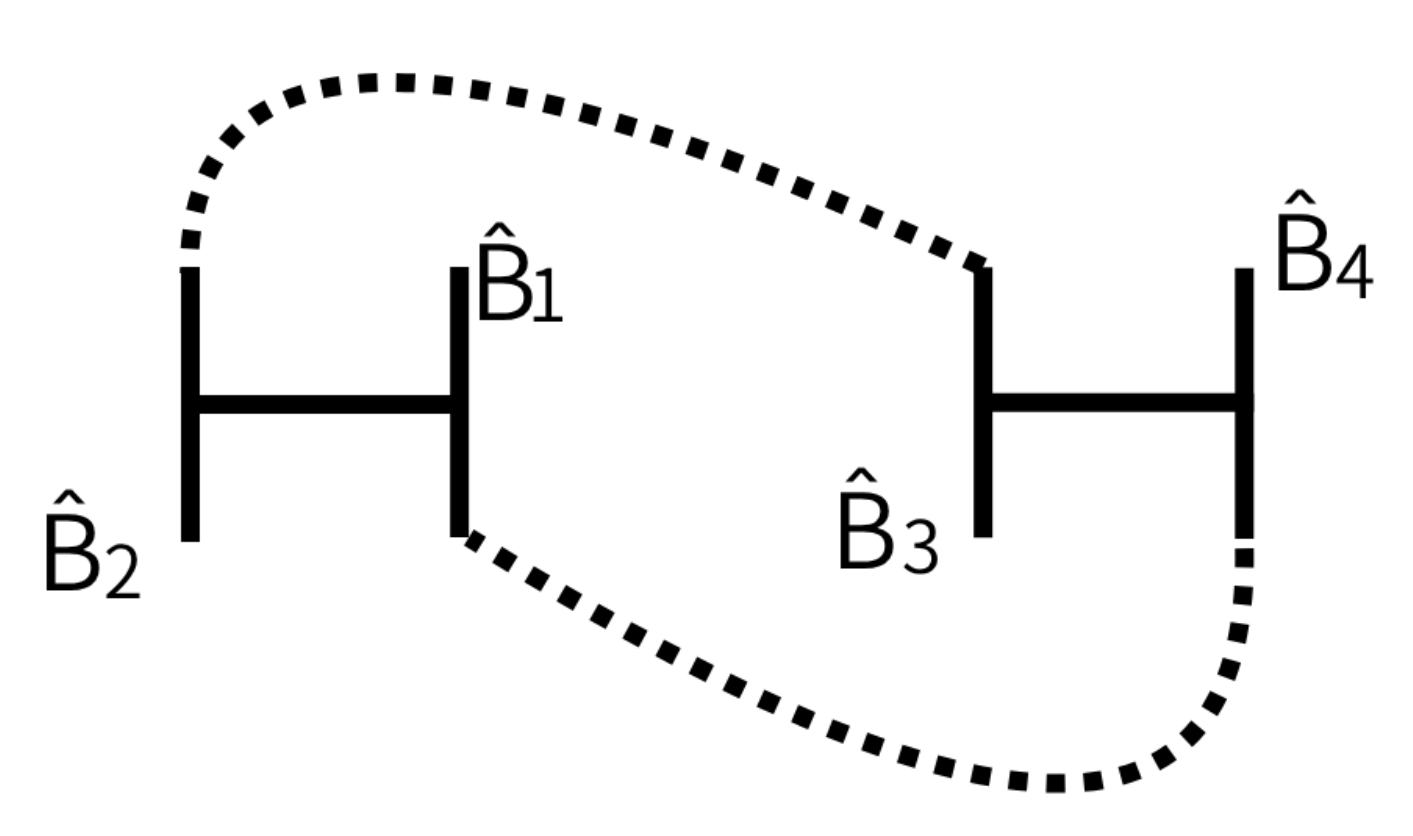}
    \caption{}
    \label{fig:Fig75c}
 \end{minipage}
\end{figure}

 (IX) Assume that $\Gamma_a$ is connected as shown in Figure \ref{fig:Fig76}.
Such $f$ does not exist.

(X) Assume that $\Gamma_a$ is connected as shown in Figure \ref{fig:Fig77}.
Such $f$ does not exist.

\begin{figure}[H]
 \begin{minipage}[ht]{0.45\linewidth}
    \centering
    \includegraphics[width=0.5\linewidth]{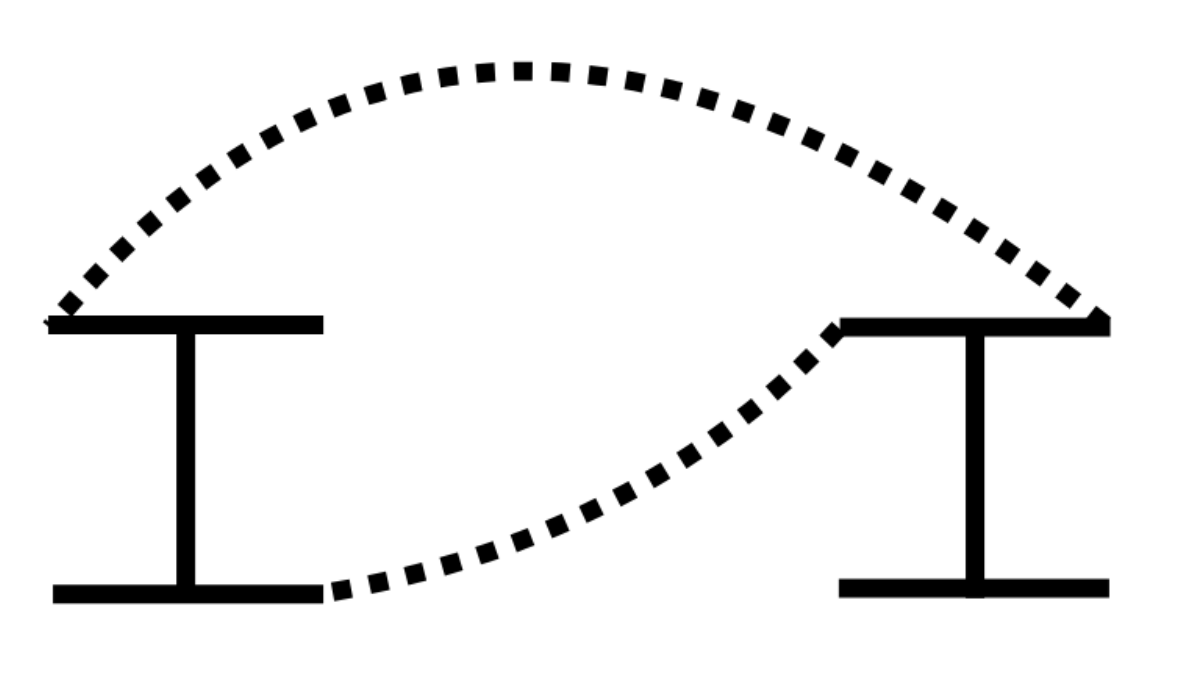}
    \caption{}
    \label{fig:Fig76}
 \end{minipage}
 \begin{minipage}[ht]{0.45\linewidth}
    \centering
    \includegraphics[width=0.5\linewidth]{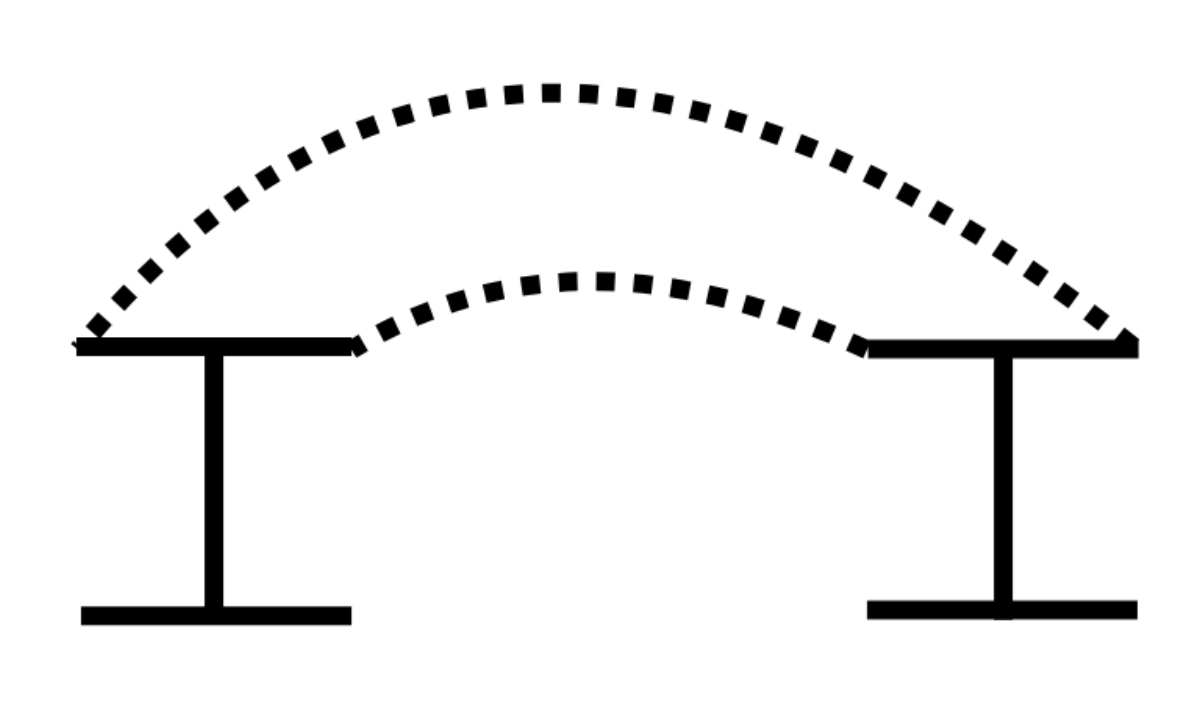}
    \caption{}
    \label{fig:Fig77}
 \end{minipage}
\end{figure}

We have checked all the cases.
This completes the proof of lemma \ref{thm:claim4}.

\end{proof}

Because this $DC$ is $2/4$ rotated by $f$, $f^2(a)=a$. From (1) of proposition \ref{thm:mainprop}, $f^2=1$.
From lemma \ref{thm:claim4}, $f^2\ne 1$. This is a contradiction. 
This proves (3), and hence proposition \ref{thm:mainprop}.
Note that we do not use the fact that the Torelli group is torsion-free.

\end{proof}

Now, let us turn to the two questions that naturally arise from our result.

In [N], Nakamura defined a two-dimensional cell complex called extended Hatcher complex
such that each 0-cell is a quilt decomposition of a surface (see definition 1 again),
each 1-cell is a move of three types between two quilt decompositions, and each 2-cell is 
one of six types of certain 2-cells.
The reference [N] is a handwritten note taken by one of the audience who attended the lecture.
For a precise definition of the extended Hatcher complex, see definition 4 of [F2].
Nakamura proved that the extended Hatcher complex is simply-connected.
For a proof, see theorem 5 of [F2].
We proved that a torsion-free subgroup of the mapping class group (in particular, the Torelli group)
 acts freely on the extended Hatcher complex
 using the fact that the group is torsion-free
(see (2) of theorem 18 of [F2]).

\begin{question}\label{q1}
    Is it possible to prove that the Torelli group acts freely on an extended Hatcher complex without using the fact
    that the Torelli group is torsion-free?
\end{question}

The following question also arises naturally.

\begin{question}\label{q2}
   Is it possible to prove that the Torelli group is torsion-free from our result?
\end{question}

\noindent
{\it Background, Acknowledgments and Tool Use}

The reference [F2] is based on [F1], which was submitted as a master’s thesis 
 to the Graduate School of Mathematical Sciences, the University of Tokyo, in March 2001. 
 Shortly after submitting my master’s thesis [F1], I fell ill and was unable to publish it.
 
Many years later, in 2025, I had the good fortune to meet Professor Koji Hasegawa.
 Since November 2025, Professor Hasegawa has encouraged me to publish [F1] and 
 helped me contact Professor Toshitake Kohno again, who had been my academic supervisor.
 I completed the writing of [F2] in September 2026.

The research presented in this paper is based on a manuscript prepared for 
 a student seminar under Professor Kohno at the University of Tokyo 
 on October 20, 2000, when I was a master’s student. 
 I had thought that the manuscript had been lost, however, it was found in early 2026 
 among some old belongings.
 
Professor Nariya Kawazumi, one of the teachers who helped me when I was a student, 
 provided me with valuable information about the Torelli groups. I am deeply grateful to 
 Professor Kohno, Professor Hasegawa and Professor Kawazumi.

I partly used the free version of Google Translate to translate Japanese to English. 
 I also used the free version of ChatGPT to find better expressions for English, to determine 
 the lettering of the notation, and to place the 96 figures on this paper, 
 while I did not use AI for mathematical reasoning or for drawing the figures.

\end{document}